\input amstex
\documentstyle{amsppt}
\magnification=\magstep1
\vsize =21 true cm 
\hsize =16 true cm 
\loadmsbm
\topmatter
\centerline{\bf Product formulas on a unitary group in 
                three variables}
\author{\smc Lei Yang}\endauthor
\endtopmatter
\document

\centerline{\bf 1. Introduction}
\vskip 0.5 cm

  In the present paper, we consider the eigenvalue problems which 
concern a differential equation
$$L_{k} f(z_1, z_2)=\lambda f(z_1, z_2). \tag 1.1$$
Here $f$ is a function on the Siegel domain of type II
$${\frak S}_{2}=\{ (z_1, z_2) \in {\Bbb C}^{2}: z_1+\overline{z_1}-
  z_2 \overline{z_2} > 0 \}, \tag 1.2$$
and $L_{k}$ is a differential operator given by
$$\aligned
  L_{k}
&=(z_1+\overline{z_1}-z_2 \overline{z_2})\\
&\times \left[ (z_1+\overline{z_1}) \frac{\partial^2}{\partial z_1
 \partial \overline{z_1}}+\frac{\partial^2}{\partial z_2 \partial
 \overline{z_2}}+z_2 \frac{\partial^2}{\partial \overline{z_1}
 \partial z_2}+\overline{z_2} \frac{\partial^2}{\partial z_1
 \partial \overline{z_2}}-k(\frac{\partial}{\partial z_1}-
 \frac{\partial}{\partial \overline{z_1}}) \right],  
\endaligned\tag 1.3$$
where $k$ is an integer. Let $\Gamma \subset U(2, 1)$ 
be an arithmetic subgroup or a convex cocompact subgroup, where 
$$U(2, 1)=\{ g \in GL(3, {\Bbb C}): g^{*} J g=J \}, \quad
          J=\left(\matrix
		      &   & -1\\
			  & 1 &   \\
		   -1 &   & 	   
		    \endmatrix\right),\tag 1.4$$ 
and assume that
$f$ is a $\Gamma$-automorphic form of weight $k$ in the sense that 
it is invariant under
$$f \mapsto f(\gamma(z_1, z_2)) (a_1 z_1+a_2 z_2+a_3)^{-k}
  \overline{(a_1 z_1+a_2 z_2+a_3)}^{k} \tag 1.5$$ 
for every $\gamma=\left(\matrix
                   *   & *   & *  \\
				   *   & *   & *  \\
				   a_1 & a_2 & a_3
				  \endmatrix\right) \in \Gamma$.
In fact, $L_{k}$ commutes with (1.5).

  Now, let us recall some basic facts about complex hyperbolic 
geometry (see \cite{Ap} and \cite{Go}). Geometry of complex 
hyperbolic space ${\bold H}_{\Bbb C}^{n}$ is the geometry of 
the unit ball ${\Bbb B}_{\Bbb C}^{n}$ in ${\Bbb C}^{n}$ with 
the K\"{a}hler structure given by the Bergman metric whose 
automorphisms are biholomorphic automorphisms of the ball, 
i.e., elements of $PU(n, 1)$. Any complex hyperbolic manifold
can be represented as the quotient 
$M={\bold H}_{\Bbb C}^{n}/\Gamma$ by a discrete torsion free
isometric action of the fundamental group of $M$,
$\pi_{1}(M) \cong \Gamma \subset PU(n, 1)$, its boundary
at infinity $\partial_{\infty} M$ is naturally identified
as the quotient $\Omega(\Gamma)/\Gamma$ of the discontinuity
set of $\Gamma$ at infinity. Here the discontinuity set
$\Omega(\Gamma)$ is the maximal subset of 
$\partial {\bold H}_{\Bbb C}^{n}$ where $\Gamma$ acts 
discretely, its complement $\Lambda(\Gamma)=\partial
{\bold H}_{\Bbb C}^{n} \backslash \Omega(\Gamma)$ is the
limit set of $\Gamma$, $\Lambda(\Gamma)=\overline{\Gamma(x)}
\cap \partial {\bold H}_{\Bbb C}^{n}$ for any 
$x \in {\bold H}_{\Bbb C}^{n}$.  

  In general, let $G$ be a connected, linear, real simple Lie
group of rank one, $G=KAN$ be an Iwasawa decomposition of $G$,   
${\frak g}={\frak k} \oplus {\frak a} \oplus {\frak n}$ be the
corresponding Iwasawa decomposition of the Lie algebra ${\frak g}$,
and $P=MAN$ be a minimal parabolic subgroup. The group $G$ acts
isometrically on the rank-one symmetric space $X=G/K$. Let 
$\partial X=G/P=K/M$ be its geodesic boundary. We regard 
$\overline{X}:=X \cup \partial X$ as a compact manifold with
boundary.

  By the classification of symmetric spaces with strictly
negative sectional curvature, we know that $X$ is one of the
following spaces: a real hyperbolic space ${\bold H}_{\Bbb R}^{n}$  
($n \geq 1$), a complex hyperbolic space ${\bold H}_{\Bbb C}^{n}$
($n \geq 2$), a quaternionic hyperbolic space ${\bold H}_{\Bbb H}^{n}$
($n \geq 2$) or the Cayley hyperbolic plane ${\bold H}_{\Bbb O}^{2}$,
and $G$ is a linear group finitely covering the 
orientation-preserving isometric group of $X$.

  Following \cite{BO}, we consider a torsion-free discrete
subgroup $\Gamma \subset G$ such that $\partial X$ admits a
$\Gamma$-invariant partition $\partial X=\Omega(\Gamma) \cup
\Lambda(\Gamma)$, where $\Omega(\Gamma) \neq \emptyset$ is
open and $\Gamma$ acts freely and cocompactly on $X \cup
\Omega(\Gamma)$. The closed subset $\Lambda(\Gamma)$ is
called the limit set of $\Gamma$. The locally symmetric
space $Y:=\Gamma \backslash X$ is a complete Riemannian
manifold of infinite volume without cusps. It can be
compactified by adjoining the geodesic boundary
$B:=\Gamma \backslash \Omega(\Gamma)$.  

  A subgroup $\Gamma$ satisfying this assumption is called
convex cocompact or geometrically cocompact since it acts
cocompactly on the convex hull of the limit set. The quotient
$Y$ is called a Kleinian manifold.  

  For $H \in {\frak a}$, we define $\rho \in {\frak a}^{*}$
by $\rho(H)=\frac{1}{2} \text{tr}(\text{ad}(H)_{|{\frak n}})$. 
Then $\rho({\bold H}_{\Bbb R}^{n})=\frac{n-1}{2}$ and
$\rho({\bold H}_{\Bbb C}^{n})=n$.

  Now, we need the following definition (see \cite{Co}):

{\it Definition}. For any discrete group $\Gamma \subset G$, the 
critical exponent $\delta(\Gamma)$ is the infimum of all $s$  
such that
$$\sum_{\gamma \in \Gamma} \exp(-s d(x, \gamma y))$$
converges for some (or any) $x, y \in X$, where $d(x, y)$ is the
Riemannian distance from $x$ to $y$. 

  The critical exponent $\delta(\Gamma)$ has been extensively
studied (see \cite{Pa1}, \cite{Su}, \cite{Co} and \cite{CoI}). 
It is known that $\delta(\Gamma) \in [0, 2 \rho]$ if $\Gamma$ 
is nontrivial. In fact, if $\Omega(\Gamma) \neq \emptyset$, 
then $\delta(\Gamma)$ is equal to the 
Hausdorff dimension of the limit set with respect to the 
natural class of sub-Riemannian metrics on $\partial X$.

  In the case of $U(2, 1)$, $\delta(\Gamma) \in [0, 4]$ for
$X={\bold H}_{\Bbb C}^{2}$.  

  The theory of Eisenstein series and the spectral theory for
Kleinian groups on real hyperbolic spaces ${\bold H}_{\Bbb R}^{n}$
has been extensively studied, in particular, by Mandouvalos
\cite{Ma1}, \cite{Ma2}, Patterson \cite{Pa3}, and Perry \cite{Pe1},
\cite{Pe2}, \cite{FHP}. A  Kleinian group $\Gamma$ is a discrete 
subgroup acting on ${\bold H}_{\Bbb R}^{n}$, which acts 
discontinuously on ${\bold H}_{\Bbb R}^{n}$ and has a fundamental
domain in ${\bold H}_{\Bbb R}^{n}$ of infinite hyperbolic volume.
This in particular implies that it is a non-arithmetic subgroup.
Their theory extended previous work of Roelcke \cite{Ro}, Elstrodt
\cite{El}, Patterson \cite{Pa2} and Fay \cite{Fa} on the spectral
theory and Eisenstein series for Fuchsian groups of the second kind.

  Now, we give the corresponding results on the complex hyperbolic
space ${\bold H}_{\Bbb C}^{2}$. We give the Eisenstein series for 
the Picard modular group, which is an arithmetic subgroup of $U(2, 1)$.
We also give the Eisenstein series for any discrete subgroup (either
arithmetic or non-arithmetic). In particular, for a convex cocompact
subgroup $\Gamma$ which satisfies that $\delta(\Gamma)<1$, we give 
the product formulas. Now, we state the results of this paper. 

  The Poisson kernel is given by
$$P(Z, W)=\frac{\rho(Z)}{|\rho(Z, W)|^2}, \quad \text{for} 
  \quad Z=(z_1, z_2) \in {\frak S}_2, W=(w_1, w_2) \in 
  \partial {\frak S}_2,\tag 1.6$$
where 
$$\rho(Z, W):=\overline{z_1}+w_1-\overline{z_2} w_2, \quad
  \rho(Z)=\rho(Z, Z).\tag 1.7$$    
Let 
$$\sigma(Z, Z^{\prime}):=\frac{|\rho(Z, Z^{\prime})|^{2}}
  {\rho(Z) \rho(Z^{\prime})}.\tag 1.8$$
It is a point-pair invariant under the action of $U(2, 1)$. 
The critical exponent is given by
$$\delta(\Gamma):=\text{inf} \{ s>0: \sum_{\gamma \in \Gamma}
  \sigma(Z, \gamma(Z^{\prime}))^{-s}<\infty \}.\tag 1.9$$

  For a discrete subgroup $\Gamma \subset U(2, 1)$, the Eisenstein 
series is defined by 			
$$E(Z, W; s):=\sum_{\gamma \in \Gamma} P(\gamma(Z), W)^s, \quad 
  \text{for} \quad \text{Re}(s)>\delta(\Gamma). \tag 1.10$$
The automorphic Green function is given by
$$G(Z, Z^{\prime}; s)=\sum_{\gamma \in \Gamma} 
  r(Z, \gamma(Z^{\prime}); s), \quad \text{for} \quad
  \text{Re}(s)>\delta(\Gamma), \tag 1.11$$
with
$$r(Z, Z^{\prime}; s)=\frac{\Gamma(s)}{\sqrt{\pi} \Gamma(s-
  \frac{1}{2})} (\frac{1}{4})^{s-1} u^{-s} F(s, s-1; 2s-1; -u^{-1}).
  \tag 1.12$$	 
Here
$$u=u(Z, Z^{\prime})=\sigma(Z, Z^{\prime})-1.\tag 1.13$$

  Now, the first main result of this paper is stated as follows:

{\smc Theorem 1.1}. {\it Assume that $\Gamma$ is convex cocompact
and $\delta(\Gamma)<1$. Then the following product formula holds:
$$\aligned
 &\int_{\Gamma \backslash \Omega(\Gamma)} E(Z, W; s) E(Z^{\prime}, W;
  2-s) dm(W)\\
=&-\frac{\pi}{s-1} c(s) c(2-s) [G(Z, Z^{\prime}; s)-
  G(Z, Z^{\prime}; 2-s)],
\endaligned\tag 1.14$$
where $\delta(\Gamma)<\text{Re}(s)<2-\delta(\Gamma)$ and
$c(s)=\sqrt{\pi} \Gamma(s-\frac{1}{2}) \Gamma(s)^{-1}$.}  

  For $\gamma=\left(\matrix
                * &   * &   *\\
			    * &   * &   *\\
			  a_1 & a_2 & a_3
              \endmatrix\right) \in U(2, 1)$, set
$$j(\gamma, Z)=a_1 z_1+a_2 z_2+a_3 z_3.\tag 1.15$$			  
The S-matrix is given by
$$S(W, W^{\prime}; s):=\sum_{\gamma \in \Gamma} |\rho(\gamma(W), 
  W^{\prime})|^{-2s} |j(\gamma, W)|^{-2s} \tag 1.16$$  
for $W, W^{\prime} \in \Omega(\Gamma)$ and $\text{Re}(s)>\delta(\Gamma)$.
By the S-matrix, we give the functional equation of Eisenstein series:

{\smc Theorem 1.2}. {\it Assume that $\Gamma$ is convex cocompact
and $\delta(\Gamma)<1$. Then the following functional equation for 
Eisenstein series holds:
$$\int_{\Gamma \backslash \Omega(\Gamma)} E(Z, W; s) S(W, W^{\prime}; 
  2-s) dm(W)=\frac{4^{s-1} \pi}{s-1} \frac{\sqrt{\pi} 
  \Gamma(s-\frac{1}{2})}{\Gamma(s)} E(Z, W^{\prime}; 2-s),
  \tag 1.17$$
for $1<\text{Re}(s)<2-\delta(\Gamma)$.}

  The Poisson kernel of weight $k$ is given by
$$P_{k}(Z, W; s)=\rho(Z)^{s} \rho(Z, W)^{k-s} 
  \rho(\overline{Z}, \overline{W})^{-k-s},\tag 1.18$$
where $Z \in {\frak S}_2$, $W \in \partial {\frak S}_2$.
The Eisenstein series of weight $k$ is defined by
$$E(Z, W; k, s)=\sum_{\gamma \in \Gamma} j(\gamma, Z)^{-k}
  j(\overline{\gamma}, \overline{Z})^{k} P_{k}(\gamma(Z), W; s),
  \quad \text{for} \quad \text{Re}(s)>\delta(\Gamma).\tag 1.19$$
Set
$$\aligned
  K(Z, W; k, s)
=&\frac{\pi^{\frac{3}{2}} \Gamma(\frac{3}{2}-s)}
  {(|k|+1-s) \Gamma(2-s)} (\frac{1}{4})^{s-1} \rho(Z, W)^{k} 
  \rho(\overline{Z}, \overline{W})^{-k}\\
 &\times \sigma^{-|k|} (\sigma-1)^{|k|-s} 
  F(s-|k|, s-1-|k|; 2s-1; -\frac{1}{\sigma-1}).
\endaligned\tag 1.20$$
The automorphic Green function of weight $k$ is given by
$$G(Z, Z^{\prime}; k, s)=\sum_{\gamma \in \Gamma}
  K(Z, \gamma(Z^{\prime}); k, s), \quad \text{for} \quad
  \text{Re}(s)>\delta(\Gamma). \tag 1.21$$

  Now, we state the main theorem of this paper: 

{\smc Theorem 1.3 (Main Theorem)}. {\it Assume that $\Gamma$ is
convex cocompact and $\delta(\Gamma)<1$. Then the following 
product formula holds:
$$\aligned
 &\int_{\Gamma \backslash \Omega(\Gamma)} E(Z, W; k, s)
  E(Z^{\prime}, W; -k, 2-s) dm(W)\\
=&G(Z, Z^{\prime}; k, s)+G(Z, Z^{\prime}; k, 2-s)\\
=&G(Z^{\prime}, Z; -k, s)+G(Z^{\prime}, Z; -k, 2-s),
\endaligned\tag 1.22$$
for $\delta(\Gamma)<\text{Re(s)}<2-\delta(\Gamma)$
and $k=0$ or $k=\pm 1$.}

  The S-matrix of weight $k$ is defined by
$$S(W, W^{\prime}; k, s)
:=\sum_{\gamma \in \Gamma} j(\gamma, W)^{k-s}
  j(\overline{\gamma}, \overline{W})^{-k-s} 
  \rho(\gamma(W), W^{\prime})^{-k-s} 
  \rho(\overline{\gamma(W)}, \overline{W^{\prime}})^{k-s}
  \tag 1.23$$  
for $W, W^{\prime} \in \Omega(\Gamma)$ and 
$\text{Re}(s)>\delta(\Gamma)$.
By the S-matrix of weight $k$, we obtain the functional equation
of Eisenstein series of weight $k$.

{\smc Theorem 1.4}. {\it Assume that $\Gamma$ is convex cocompact 
and $\delta(\Gamma)<1$. Then the following functional equation for 
Eisenstein series of weight $k$ holds:
$$\aligned
 &\int_{\Gamma \backslash \Omega(\Gamma)} E(Z, W; k, s) 
  S(W, W^{\prime}; -k, 2-s) dm(W)\\
=&\frac{\pi^{\frac{3}{2}} \Gamma(s-\frac{1}{2})}
  {(|k|+s-1) \Gamma(s)} 4^{s-1} E(Z, W^{\prime}; k, 2-s),
\endaligned\tag 1.24$$
for $1<\text{Re}(s)<2-\delta(\Gamma)$ and $k=0, \pm 1$.}

  In fact, when $k=0$, Theorem 1.3 and Theorem 1.4 reduces to
Theorem 1.1 and Theorem 1.2.  

  In the appendix, we give the product formulas on $SL(2, {\Bbb R})$.
  
{\smc Theorem 1.5}. {\it Assume that $\Gamma$ is convex cocompact
and $\delta(\Gamma)<\frac{1}{2}$. Then the following product 
formula on $SL(2, {\Bbb R})$ holds:
$$\int_{\Gamma \backslash \Omega(\Gamma)} E(z, \zeta; k, s)
  E(z^{\prime}, \zeta; -k, 1-s) dm(\zeta)
 =G(z, z^{\prime}; k, s)+G(z, z^{\prime}; k, 1-s),\tag 1.25$$
for $\delta(\Gamma)<\text{Re}(s)<1-\delta(\Gamma)$ and $k=0, \pm 1$.}

\vskip 0.5 cm
\centerline{\bf 2. Epstein zeta function and Eisenstein series 
                   on Picard modular groups}
\vskip 0.5 cm

  Let ${\frak S}_2$ be the Siegel domain
$${\frak S}_2=\{ (z_1, z_2) \in {\Bbb C}^2: z_1+
  \overline{z_1}>|z_2|^2 \}. \tag 2.1$$

  The imaginary quadratic field $K={\Bbb Q}(\sqrt{-3})$ is 
called the field of Eisenstein numbers. Its ring of integers
$${\Cal O}={\Cal O}_K=\{ a/2+b \sqrt{-3}/2: a, b \in {\Bbb Z}, 
  a \equiv b (\text{mod} 2) \}={\Bbb Z}+{\Bbb Z} [\omega] \tag 2.2$$
with $\omega=\frac{-1+\sqrt{-3}}{2}$ is called the ring of
Eisenstein integers.

  In their paper \cite{KW}, Kor\'{a}nyi and Wolf studied the general 
Cayley transform which carries the bounded domain of the Harish-Chandra
realization into a generalized half-plane. In \cite{KR}, Kor\'{a}nyi
and Reimann gave the Cayley transform on $U(2, 1)$. Now, we define 
the Cayley transform
$C: {\Bbb B}^2:=\{(w_1, w_2) \in {\Bbb C}^2: |w_1|^2+|w_2|^2<1 \}
\to {\frak S}_2$, where
$$C=\left(\matrix
     -\overline{\omega} &   & -\omega\\
                        & 1 &        \\
                     -1 &   & 1
    \endmatrix\right) \in GL(3, {\Cal O}). \tag 2.3$$ 

  It is known that
$$U(2, 1)=\{ g \in GL(3, {\Bbb C}): g I_{2, 1} g^{*}=I_{2, 1} \},
  \tag 2.4$$
where $I_{2, 1}=\text{diag}\{ 1, 1, -1 \}$.

  Let $E_{i, j}$ be the $3 \times 3$ matrix with a $1$ at the $i, j$-th
component, zeros elsewhere. Using Kronecker's delta, we have
$$[E_{p, q}, E_{r, s}]=\delta_{q, r} E_{p, s}-\delta_{s, p} E_{r, q}.$$
A basis of ${\frak g}=u(2, 1)$ is as follows:
$$X_1=E_{2, 3}+E_{3, 2}, \quad X_2=E_{1, 3}+E_{3, 1}, \quad
  X_3=E_{1, 2}-E_{2, 1},$$
$$X_4=i E_{2, 3}-i E_{3, 2}, \quad X_5=i E_{1, 3}-i E_{3, 1}, \quad
  X_6=i E_{1, 2}+i E_{2, 1},$$
$$X_7=i E_{1, 1}, \quad X_8=i E_{2, 2}, \quad X_9=i E_{3, 3}.$$
Thus,
$$E_{2, 3}=\frac{1}{2}(X_1-i X_4), \quad 
  E_{3, 2}=\frac{1}{2}(X_1+i X_4),$$
$$E_{1, 3}=\frac{1}{2}(X_2-i X_5), \quad
  E_{3, 1}=\frac{1}{2}(X_2+i X_5),$$
$$E_{1, 2}=\frac{1}{2}(X_3-i X_6), \quad
  E_{2, 1}=-\frac{1}{2}(X_3+i X_6),$$
$$E_{1, 1}=-i X_7, \quad E_{2, 2}=-i X_8, \quad
  E_{3, 3}=-i X_9.$$

 The action of $U(2, 1)$ on ${\Bbb B}^2$ is as follows:
$$\left(\matrix
  a_1 &a_2 &a_3\\
  b_1 &b_2 &b_3\\
  c_1 &c_2 &c_3
\endmatrix\right) (w_1, w_2)=\left(\frac{a_1 w_1+a_2 w_2+a_3}
{c_1 w_1+c_2 w_2+c_3}, \frac{b_1 w_1+b_2 w_2+b_3}{c_1 w_1+c_2 w_2+c_3}
\right).$$

  We have
$$\aligned
 &X_3 f(w_1, w_2, \overline{w_1}, \overline{w_2}) \\
=&\frac{d}{dt}[f(\exp(t \left(\matrix
		  0 & 1 & 0\\
		 -1 & 0 & 0\\
		  0 & 0 & 0
   \endmatrix\right)) (w_1, w_2))+
    f(\exp(t \left(\matrix
		  0 & 1 & 0\\
		 -1 & 0 & 0\\
		  0 & 0 & 0
   \endmatrix\right)) (\overline{w_1}, 
   \overline{w_2}))]|_{t=0}\\
=&\frac{d}{dt}[f(\left(\matrix
		  \cos t & \sin t & 0\\
		 -\sin t & \cos t & 0\\ 
		  0      &  0     & 1
    \endmatrix\right) (w_1, w_2))+
    f(\left(\matrix
		  \cos t & \sin t & 0\\
		 -\sin t & \cos t & 0\\ 
		  0      &  0     & 1
		 \endmatrix\right) (\overline{w_1}, 
  \overline{w_2}))]|_{t=0}\\
=&\frac{d}{dt}f(w_1 \cos t+w_2 \sin t, -w_1 \sin t+w_2 
  \cos t)|_{t=0}\\
 &+\frac{d}{dt}f(\overline{w_1} \cos t+\overline{w_2} \sin t, 
   -\overline{w_1} \sin t+\overline{w_2} \cos t)|_{t=0}\\
=&w_2 \frac{\partial f}{\partial w_1}-w_1 \frac{\partial f}
  {\partial w_2}+\overline{w_2} \frac{\partial f}{\partial 
  \overline{w_1}}-\overline{w_1} \frac{\partial f}{\partial 
  \overline{w_2}}.
\endaligned$$
So, the corresponding operator is
$$X_3=w_2 \frac{\partial}{\partial w_1}-w_1 \frac{\partial}
      {\partial w_2}+\overline{w_2} \frac{\partial}{\partial 
	  \overline{w_1}}-\overline{w_1} \frac{\partial}{\partial 
	  \overline{w_2}}.$$
  
  Similarly, we have
$$X_2=(1-w_1^2) \frac{\partial}{\partial w_1}-w_1 w_2 
      \frac{\partial}{\partial w_2}+(1-\overline{w_1}^2) 
	  \frac{\partial}{\partial \overline{w_1}}-\overline{w_1} 
	  \overline{w_2} \frac{\partial}{\partial \overline{w_2}}.$$
$$X_1=-w_1 w_2 \frac{\partial}{\partial w_1}+(1-w_2^2) 
      \frac{\partial}{\partial w_2}-\overline{w_1} \overline{w_2} 
	  \frac{\partial}{\partial \overline{w_1}}+(1-\overline{w_2}^2) 
	  \frac{\partial}{\partial \overline{w_2}}.$$
$$X_6=iw_2 \frac{\partial}{\partial w_1}
     +iw_1 \frac{\partial}{\partial w_2}
     -i \overline{w_2} \frac{\partial}{\partial \overline{w_1}}
     -i \overline{w_1} \frac{\partial}{\partial \overline{w_2}}.$$
$$X_5=i(1+w_1^2) \frac{\partial}{\partial w_1}
     +i w_1 w_2 \frac{\partial}{\partial w_2}
     -i(1+\overline{w_1}^2) \frac{\partial}{\partial \overline{w_1}}
     -i \overline{w_1} \overline{w_2} \frac{\partial}
     {\partial \overline{w_2}}.$$
$$X_4=i w_1 w_2 \frac{\partial}{\partial w_1}
     +i(1+w_2^2)\frac{\partial}{\partial w_2}
     -i \overline{w_1} \overline{w_2}\frac{\partial}
	 {\partial \overline{w_1}}-i(1+\overline{w_2}^2)
	 \frac{\partial}{\partial \overline{w_2}}.$$
$$X_7=iw_1 \frac{\partial}{\partial w_1}
     -i\overline{w_1} \frac{\partial}{\partial \overline{w_1}}.$$
$$X_8=i w_2 \frac{\partial}{\partial w_2}
     -i \overline{w_2} \frac{\partial}{\partial \overline{w_2}}.$$
$$X_9=-i w_1 \frac{\partial}{\partial w_1}
      -i w_2 \frac{\partial}{\partial w_2}
      +i \overline{w_1} \frac{\partial}{\partial \overline{w_1}}
      +i \overline{w_2} \frac{\partial}{\partial \overline{w_2}}.$$

  We denote $X_{ij}$ as the differential operator corresponding 
to $E_{i, j}$, then  
$$X_{23}=\frac{1}{2}(X_1-i X_4)=\frac{\partial}{\partial w_2}-
  \overline{w_1} \overline{w_2} \frac{\partial}{\partial 
  \overline{w_1}}-\overline{w_2}^2 \frac{\partial}{\partial
  \overline{w_2}}, \quad
  X_{32}=\frac{1}{2}(X_1+i X_4)=\overline{X_{23}},$$
$$X_{13}=\frac{1}{2}(X_2-i X_5)=\frac{\partial}{\partial w_1}
  -\overline{w_1}^2 \frac{\partial}{\partial \overline{w_1}}-
  \overline{w_1} \overline{w_2} \frac{\partial}{\partial
  \overline{w_2}}, \quad
  X_{31}=\frac{1}{2}(X_2+i X_5)=\overline{X_{13}},$$
$$X_{12}=\frac{1}{2}(X_3-i X_6)=w_2 \frac{\partial}{\partial w_1}
  -\overline{w_1} \frac{\partial}{\partial \overline{w_2}}, \quad
  X_{21}=-\frac{1}{2}(X_3+i X_6)=-\overline{X_{12}},$$
$$X_{11}=-i X_7=w_1 \frac{\partial}{\partial w_1}-
  \overline{w_1} \frac{\partial}{\partial \overline{w_1}},$$
$$X_{22}=-i X_8=w_2 \frac{\partial}{\partial w_2}-
  \overline{w_2} \frac{\partial}{\partial \overline{w_2}},$$
$$X_{33}=-i X_9=-w_1 \frac{\partial}{\partial w_1}-w_2
  \frac{\partial}{\partial w_2}+\overline{w_1} \frac{\partial}
  {\partial \overline{w_1}}+\overline{w_2} \frac{\partial}
  {\partial \overline{w_2}}.$$

  Now, if $n$ is any positive integer, we claim that
$$\sum_{i_1=1}^{3} \sum_{i_2=1}^{3} \cdots \sum_{i_n=1}^{3}
  X_{i_1 i_2} X_{i_2 i_3} \cdots X_{i_n i_1}$$
lies in the center of $U({\frak g})$. To prove this, without loss
of generality, let us assume $n=2$. Then 
$$[X_{pq}, \sum_{i_1} \sum_{i_2} X_{i_1 i_2} X_{i_2 i_1}]=
  \sum_{i_1} \sum_{i_2}([X_{pq}, X_{i_1 i_2}] X_{i_2 i_1}+
  X_{i_1 i_2} [X_{pq}, X_{i_2 i_1}]).$$
By $[X_{pq}, X_{rs}]=\delta_{qr} X_{ps}-\delta_{sp} X_{rq}$, 
this equals
$$\aligned
  &\sum_{i_1} \sum_{i_2} (\delta_{q i_1} X_{p i_2} X_{i_2 i_1}-
   \delta_{i_2 p} X_{i_1 q} X_{i_2 i_1}+X_{i_1 i_2} \delta_{q i_2}
   X_{p i_1}-X_{i_1 i_2} \delta_{i_1 p} X_{i_2 q})\\
 =&\sum_{i_2} X_{p i_2} X_{i_2 q}-\sum_{i_1} X_{i_1 q} X_{p i_1}+
   \sum_{i_1} X_{i_1 q} X_{p i_1}-\sum_{i_2} X_{p i_2} X_{i_2 q}=0.
\endaligned$$
Thus the above sum commutes with the generators $X_i$ of 
$U({\frak g})$, and lies in the center of $U({\frak g})$. 
The proof is similar if $n \ne 2$.

  Set 
$$D_1=\sum_{i_1=1}^{3} \sum_{i_2=1}^{3} E_{i_1 i_2} E_{i_2 i_1},$$
$$D_2=\sum_{i_1=1}^{3} \sum_{i_2=1}^{3} \sum_{i_3=1}^{3}
  E_{i_1 i_2} E_{i_2 i_3} E_{i_3 i_1}.$$
In the category of differential operators,
$$[X_{pq}, X_{rs}]=\delta_{sp} X_{rq}-\delta_{qr} X_{ps}.$$
Therefore, we have
$$D_1=2(X_{11}^2+X_{22}^2+X_{11}X_{22}+X_{12}X_{21}+X_{23}X_{32}+
        X_{31}X_{13})=2 \Delta,$$
$$\aligned
D_2=&3[X_{11}(X_{12}X_{21}-X_{32}X_{23}-X_{22}^2)+X_{22}(X_{21}
       X_{12}-X_{31}X_{13}-X_{11}^2)\\
	&+X_{23}X_{31}X_{12}+X_{13}X_{32}X_{21}-X_{31}X_{13}-
	  X_{32}X_{23}+X_{11}+X_{22}]=-3 \Delta,
\endaligned$$
where
$$\aligned
  \Delta
 &=(1-|w_1|^2-|w_2|^2)\\
 &\times \left[(1-|w_1|^2) \frac{\partial^2}{\partial
  w_1 \partial \overline{w_1}}+(1-|w_2|^2) \frac{\partial^2}
  {\partial w_2 \partial \overline{w_2}}-w_1 \overline{w_2} 
  \frac{\partial^2}{\partial w_1 \partial \overline{w_2}}-
  \overline{w_1} w_2 \frac{\partial^2}{\partial \overline{w_1} 
  \partial w_2}\right] 
\endaligned\tag 2.5$$
is the Laplace-Beltrami operator of $U(2, 1)$ on ${\Bbb B}^2$.

  Let $G=C U(2, 1) C^{-1}$, $K=C (U(2) \times U(1)) C^{-1}$ and
$G({\Cal O})=C U(2, 1; {\Cal O}) C^{-1}$, then $G({\Cal O})$
is an arithmetic subgroup of $G$ and $S:=G({\Cal O}) \backslash G/K$ 
is called a Picard modular surface (see \cite{Pic1}, \cite{Pic2},
\cite{Gir}, \cite{Ho1}, \cite{Ho2}, \cite{LR} and \cite{BB} 
for the details).

  The complex hyperbolic Gauss-Bonnet formula states that 
(see \cite{HP} and \cite{P})
$$\text{Vol} ({\Bbb H}_{\Bbb C}^2 / \Gamma)=\frac{8 \pi^2}{3}
  \chi({\Bbb H}_{\Bbb C}^2 / \Gamma),$$
where $\chi({\Bbb H}_{\Bbb C}^2/ \Gamma)$ is the Euler number of
${\Bbb H}_{\Bbb C}^2 / \Gamma$.
Holzapfel showed that
$$\chi({\Bbb B}^{2}/PSU(2, 1; {\Cal O}))=\frac{1}{24}.$$
As $PSU(2, 1; {\Cal O})$ has index $3$ in $PU(2, 1; {\Cal O})$, 
we see that $\chi({\Bbb B}^{2}/PU(2, 1; {\Cal O}))=\frac{1}{72}$, 
and so
$$\text{Vol}_{{\Bbb B}^{2}}({\Bbb B}^{2}/PU(2, 1; {\Cal O}))
 =\frac{8 \pi^2}{3} \cdot \frac{1}{72}=\frac{\pi^2}{27}.$$

  Denote $J=\left(\matrix
            0 & 0 & -1\\
			0 & 1 &  0\\
		   -1 & 0 &  0
		   \endmatrix\right)$, then
$C I_{2, 1} C^{*}=J$. It follows that
$$G=\{ g \in GL(3, {\Bbb C}): g J g^{*}=J \}.\tag 2.6$$
This is the other realization of $U(2, 1)$. For simplicity, 
from now on, we use the same symbol $U(2, 1)$ to denote as $G$. 
In fact, $G=\{ g \in GL(3, {\Bbb C}): g^{*} J g=J \}$.
   
  Let $P$ be the parabolic subgroup of $U(2, 1)$ which consists 
of upper triangular matrices. This group $P$ has the Langlands 
decomposition $P=NAM$, where  
$$N=\left\{n=[z, t]=\left(\matrix 
    1 & z & \frac{1}{2} |z|^2+it\\
	0 & 1 & \overline{z}\\
	0 & 0 & 1
	\endmatrix\right): z \in {\Bbb C}, t \in {\Bbb R} \right\}.$$
$$A=\left\{a=\left(\matrix
    \rho^{-1} &   &     \\
	          & 1 &     \\
			  &   & \rho
    \endmatrix\right): \rho > 0 \right\},$$
$$M=\left\{m=\left(\matrix
    \beta &            &      \\
	      & \beta^{-2} &      \\
	      &            & \beta
	\endmatrix\right): \beta \in {\Bbb C}, |\beta|^2=1 \right\}.$$	
Then $G=NAMK=NAK$ is the Iwasawa decomposition.
Set $\Gamma=G({\Cal O})$ , then
$$\Gamma \cap N=\left\{ [\alpha, \beta]=\left(\matrix
                         1 & \alpha & \beta\\
                           &      1 & \overline{\alpha}\\
                           &        &     1
                         \endmatrix\right):
  |\alpha|^2=\beta+\overline{\beta}, \alpha, \beta \in 
  {\Cal O} \right\}.\tag 2.7$$
If $k \in K$, then 
$$k=C \left(\matrix
      \alpha & \beta  & 0\\
      \gamma & \delta & 0\\
	       0 &      0 & e^{i \theta}
      \endmatrix\right) C^{-1},$$
where $\left(\matrix
       \alpha & \beta\\
	   \gamma & \delta
	   \endmatrix\right) \in U(2)$ and $\theta \in [0, 2 \pi)$.
Write explicitly,
$$k=\left(\matrix
    -\overline{\omega} \alpha-\omega e^{i \theta}
   &-\overline{\omega} \beta & -\alpha+e^{i \theta}\\	 
    \gamma & \delta & \omega \gamma\\
    -\alpha+e^{i \theta} & -\beta 
   &-\omega \alpha-\overline{\omega} e^{i \theta}
   \endmatrix\right).$$
We have $k(-\omega, 0)=(-\omega, 0)$.
Note that $\rho(-\omega, 0)=1$. 
Therefore, $g(-\omega, 0)=(-\omega, 0)$ if and only if $g \in K$.
Define $\pi: G \to {\frak S}_2$ by $\pi(g)=g(-\omega, 0)$. It 
induces a bijection $G/K \to {\frak S}_2$. For $g=nak$,
$$\pi(g)=nak(-\omega, 0)=na(-\omega, 0)=(-\omega \rho^{-2}+
  \frac{1}{2} |z|^2+it, \overline{z}).$$
Thus, if we set $(z_1, z_2)=\pi(g)$, then $z_1+\overline{z_1}-z_2
\overline{z_2}=\rho^{-2}$.

  Following \cite{He}, let $\Delta=\Delta(\Gamma)=\left\{ 
\left(\matrix
       1 &   & i x\\
	     & 1 &    \\
		 &   & 1
\endmatrix\right): x \in {\Bbb R} \right\} \cap \Gamma$ and
$q=q(\Gamma)$ be the positive integer such that
$\left(\matrix
        1 &   & q \sqrt{-3}\\
		  & 1 &            \\
		  &   & 1
\endmatrix\right)$ generates $\Delta$. For $[\alpha, \beta], 
[\alpha^{\prime}, \beta^{\prime}] \in \Gamma \cap N$,
$$[\alpha^{\prime}, \beta^{\prime}] [\alpha, \beta]
  ([\alpha, \beta] [\alpha^{\prime}, \beta^{\prime}])^{-1}=
  [0, 2 \sqrt{-3} \text{Im}(\alpha^{\prime} \overline{\alpha})]
  \in \Delta,$$ 
where $\text{Im} \alpha:=b$ for $\alpha=a+b \sqrt{-3}$, 
$a, b \in {\Bbb Z}$. Hence, $2 q^{-1} \text{Im}(\alpha^{\prime}
\overline{\alpha}) \in {\Bbb Z}$.

  By Cayley transform, we have

{\smc Theorem 2.1}. {\it The Laplace-Beltrami operator of 
$U(2, 1)$ on ${\frak S}_2$ is
$$L=(z_1+\overline{z_1}-|z_2|^2) \left[(z_1+\overline{z_1}) 
    \frac{\partial^2}{\partial z_1 \partial \overline{z_1}}+
	\frac{\partial^2}{\partial z_2 \partial \overline{z_2}}+
	z_2 \frac{\partial^2}{\partial \overline{z_1} \partial z_2}+
	\overline{z_2} \frac{\partial^2}{\partial z_1 \partial
    \overline{z_2}}\right]. \tag 2.8$$}

  If $g=\left(\matrix
         1 & \alpha & \beta\\
           &      1 & \gamma\\
           &        &      1
        \endmatrix\right) \in \Gamma$, then
$\alpha=\overline{\gamma}$, $|\alpha|^2=\beta+\overline{\beta}$.
Denote $Z=(z_1, z_2)$, set $\rho(Z)=z_1+\overline{z_1}-|z_2|^2$, 
we have
$$\rho(g, Z):=\rho(g(Z))=(z_1+\alpha z_2+\beta)+(\overline{z_1}
  +\overline{\alpha} \overline{z_2}+\overline{\beta})-(z_2+
  \gamma)(\overline{z_2}+\overline{\gamma})=\rho(Z).$$

  If $g=\left(\matrix
          * & * & *\\
          * & * & *\\
          a_1 & a_2 & a_3
         \endmatrix\right) \in U(2, 1)$, then
$a_1 \overline{a_3}+\overline{a_1} a_3=a_2 \overline{a_2}$.
By a straightforward calculation, we have
$$\rho(g, Z)=\frac{\rho(Z)}{|j(\gamma, Z)|^2}, \quad \text{with}
  \quad j(\gamma, Z):=a_1 z_1+a_2 z_2+a_3. \tag 2.9$$

  For the Picard modular group $\Gamma=U(2, 1; {\Cal O})$, the 
Eisenstein series is defined as
$$E(Z; s)=\sum_{\gamma \in \Gamma \cap N \backslash
  \Gamma} \rho(\gamma,Z)^s, \quad \text{for} \quad
  \text{Re}(s)>2. \tag 2.10$$

  The unit group of ${\Cal O}$ is $\{ \pm 1, \pm \omega,
\pm \overline{\omega} \}$, its order is $6$.
$$E(Z; s)=\frac{1}{6} \sum_{\aligned
 &a_1, a_2, a_3 \in {\Cal O}, (a_1, a_2, a_3)={\Cal O},\\
 &a_1 \overline{a_3}+\overline{a_1} a_3=a_2 \overline{a_2}
 \endaligned} \frac{\rho^s}{|a_1 z_1+a_2 z_2+a_3|^{2s}},$$
where $(a_1, a_2, a_3)$ denotes the ideal generated by $a_1$, 
$a_2$ and $a_3$.
               
{\smc Theorem 2.2}. {\it The Eisenstein series satisfies the
following equation:
$$L E(Z; s)=s(s-2) E(Z; s). \tag 2.11$$}

{\it Proof}. For $j(\gamma, Z)=a_1 z_1+a_2 z_2+a_3$,
$$L \rho(\gamma, Z)^s=s(s-2) \rho(\gamma, Z)^s+(-\overline{a_1} 
  a_3-a_1 \overline{a_3}+a_2 \overline{a_2}) s^2 \frac{\rho(Z)^s}
  {|j(\gamma, Z)|^{2(s+1)}}.$$
$\qquad \qquad \qquad \qquad \qquad \qquad \qquad \qquad \qquad
 \qquad \qquad \qquad \qquad \qquad \qquad \qquad \qquad \qquad
 \quad \boxed{}$

  Let
$${\Cal S} {\Cal P}_3({\Bbb C})=\{ Y \in M_3({\Bbb C}):
  Y \quad \text{positive definite}, \quad Y^{*}=Y, 
  \text{det}Y=1 \},$$  
${\Cal H}_3({\Bbb C}):={\Cal S}{\Cal P}_3({\Bbb C}) 
\cap U(2, 1)$ and denote $Y[a]=a^{*} Y a$, for 
$Y \in {\Cal S} {\Cal P}_3({\Bbb C})$.

  For $a \in A, n \in N$ as above,
$$a[n]:=n^{*} a n=\left(\matrix
  \frac{1}{\rho} & \frac{z}{\rho} & \frac{\frac{1}{2} |z|^2+it}{\rho}\\
  \frac{\overline{z}}{\rho} & \frac{|z|^2}{\rho}+1 
& \overline{z} (\frac{\frac{1}{2} |z|^2+it}{\rho}+1)\\
  \frac{\frac{1}{2}|z|^2-it}{\rho} & z (\frac{\frac{1}{2}|z|^2-it}
  {\rho}+1) & \frac{\frac{1}{4}|z|^4+t^2}{\rho}+|z|^2+\rho
\endmatrix\right).$$

  We can identify ${\frak S}_2$ with ${\Cal H}_3({\Bbb C})$ by
$(z_1, z_2) \mapsto W_{(\rho, t, z)}=a[n]$, where
$$z_1=\rho+\frac{1}{2} |z|^2-it, z_2=-\overline{z}, \rho=\frac{1}{2}
      (z_1+\overline{z_1}-|z_2|^2).$$  

{\smc Lemma 2.3}. {\it If $a=(a_3, -a_2, a_1) \in {\Cal O}^3$, 
where $a_1 \overline{a_3}+\overline{a_1} a_3=a_2 \overline{a_2}$, then
$$W_{(\rho, t, z)}[a^{*}]=\rho^{-1} |a_1 z_1+a_2 z_2+a_3|^2. 
  \tag 2.12$$}

{\it Proof}. In fact,
$$\rho^{-1} |a_1 z_1+a_2 z_2+a_3|^2-W_{(\rho, t, z)}[a^{*}]=
  a_1 \overline{a_3}+\overline{a_1} a_3-a_2 \overline{a_2}=0.$$
$\qquad \qquad \qquad \qquad \qquad \qquad \qquad \qquad \qquad 
 \qquad \qquad \qquad \qquad \qquad \qquad \qquad \qquad \qquad
 \quad \boxed{}$ 

{\it Definition} {\smc 2.4}. For the Picard modular group 
$\Gamma=U(2, 1; {\Cal O})$ the Epstein zeta function is given by
$$Z(Y, s):=\frac{1}{6} \sum_{a \in {\Cal O}^3-\{ 0 \},
           a_1 \overline{a_3}+\overline{a_1} a_3=a_2 \overline{a_2}}
		   Y[a]^{-s}, \tag 2.13$$
where $\text{Re}(s)>2.$

{\smc Lemma 2.5}. {\it $Z(W_{(\rho, t, z)}, s)=\zeta_{K}(2s) 
E(z_1, z_2; s)$, where $\zeta_{K}(s)$ is the Dedekind zeta function 
of the number field $K={\Bbb Q}(\sqrt{-3})$.}

  Because $E(z_1, z_2; s)=Z(W_{(\rho, t, z)}, s)/\zeta_{K}(2s)$, 
any zero of $\zeta_{K}(2s)$ will give a pole of $E(z_1, z_2; s)$ 
unless $Z(W, s)$ vanishes at that value of $s$. The trivial zeros 
of $Z(W, s)$ and $\zeta_{K}(2s)$ are both $s=-1, -2, -3, \cdots$, 
and both have order $1$. So the trivial zeros cancel out in 
$E(z_1, z_2; s)$.

  Suppose $a_1, a_2, a_3 \in {\Cal O}$ and $\overline{a_1} a_3+a_1
\overline{a_3}=a_2 \overline{a_2}$, then
$J[(a_1, a_2, a_3)^{*}]=0$ and $J[(a_3, -a_2, a_1)^{*}]=0$.

  For $n$, $a \in U(2, 1)$ as above, then $n^{*} \in U(2, 1)$, thus 
$W \in U(2, 1)$. If $g \in {\Cal H}_3({\Bbb C})$, $g J g=J$, thus 
$g^{-1}=J g J$. Suppose $a=(a_1, a_2, a_3) \in {\Cal O}^3-\{ 0 \}$, 
we have the following lemma:

{\smc Lemma 2.6}. {\it For $g \in {\Cal H}_3({\Bbb C})$, 
$g [(a_3, -a_2, a_1)^{*}]=g^{-1} [(a_1, a_2, a_3)^{*}]$.}
 
{\it Proof}. In fact,
$$g^{-1}[(a_1, a_2, a_3)^{*}]=JgJ[(a_1, a_2, a_3)^{*}]=
 (-a_3, a_2, -a_1)g \left(\matrix
  -\overline{a_3}\\ \overline{a_2}\\ -\overline{a_1}
  \endmatrix\right)=g[(a_3, -a_2, a_1)^{*}].$$
$\qquad \qquad \qquad \qquad \qquad \qquad \qquad \qquad \qquad
 \qquad \qquad \qquad \qquad \qquad \qquad \qquad \qquad \qquad
 \quad \boxed{}$ 
  
  By Lemma 2.6, it follows that
    
{\smc Theorem 2.7}. {\it The functional equation of Epstein 
zeta functions: 
$$Z(Y^{-1}, s)=Z(Y, s) \tag 2.14$$ 
if $Y \in {\Cal H}_3({\Bbb C})$.}

  Following \cite{Ho2}, we define the Picard modular forms on 
${\frak S}_2$ as follows:

{\it Definition} {\smc 2.8}. A holomorphic function 
$f: {\frak S}_2 \to {\Bbb C}$ is a Picard modular form of 
$K={\Bbb Q}(\sqrt{-3})$ and of weight $k$, if it satisfies:
$$f(\gamma(Z))=\text{Jac}(\gamma,Z)^{-k} f(Z) \quad 
  \text{for all} \quad \gamma \in \Gamma, \tag 2.15$$
where $\text{Jac}(\gamma, Z)$ is the Jacobi determinant of 
$\gamma: {\frak S}_2 \to {\frak S}_2$ at $Z=(z_1, z_2)$.

  Now, for the Picard modular group $\Gamma=U(2, 1; {\Cal O})$, 
we define the Eisenstein series of weight $k$:
$$E_k(Z)=\sum_{\gamma \in \Gamma \cap N \backslash \Gamma}
  \text{Jac}(\gamma, Z)^k, \quad \text{for} \quad k>1. 
  \tag 2.16$$
In fact, in \cite{BB} and \cite{Sh1}, the authors studied 
the Poincar\'{e}-Eisenstein series.  
    
  The following proposition is important.
  	
{\smc Proposition 2.9}. {\it $j(g, Z)$ satisfies the following
identity:
$$j(g_1 g_2, Z)=j(g_1, g_2(Z)) j(g_2, Z).\tag 2.17$$}

{\it Proof}. For
$g_1=\left(\matrix
      a_1 & a_2 & a_3\\
	  b_1 & b_2 & b_3\\  
      c_1 & c_2 & c_3
      \endmatrix\right) \in U(2, 1)$ and
$g_2=\left(\matrix
      a_1^{\prime} & a_2^{\prime} & a_3^{\prime}\\
      b_1^{\prime} & b_2^{\prime} & b_3^{\prime}\\ 
	  c_1^{\prime} & c_2^{\prime} & c_3^{\prime}
     \endmatrix\right) \in U(2, 1)$, we have
$$g_1 g_2=\left(\matrix
           *  &  *  &  *\\               
           *  &  *  &  *\\  
		a_1^{\prime} c_1+b_1^{\prime} c_2+c_1^{\prime} c_3
	   &a_2^{\prime} c_1+b_2^{\prime} c_2+c_2^{\prime} c_3
	   &a_3^{\prime} c_1+b_3^{\prime} c_2+c_3^{\prime} c_3		
       \endmatrix\right).$$	 
Thus,
$$j(g_1 g_2, Z)
 =(a_1^{\prime} c_1+b_1^{\prime} c_2+c_1^{\prime} c_3) z_1
 +(a_2^{\prime} c_1+b_2^{\prime} c_2+c_2^{\prime} c_3) z_2
 +(a_3^{\prime} c_1+b_3^{\prime} c_2+c_3^{\prime} c_3).$$
On the other hand.
$$g_2(Z)=\left(\frac{a_1^{\prime} z_1+a_2^{\prime} z_2+a_3^{\prime}}
  {c_1^{\prime} z_1+c_2^{\prime} z_2+c_3^{\prime}},
  \frac{b_1^{\prime} z_1+b_2^{\prime} z_2+b_3^{\prime}}
  {c_1^{\prime} z_1+c_2^{\prime} z_2+c_3^{\prime}}\right).$$
$$j(g_1, g_2(Z))=c_1 \frac{a_1^{\prime} z_1+a_2^{\prime} z_2
  +a_3^{\prime}}{c_1^{\prime} z_1+c_2^{\prime} z_2+c_3^{\prime}}
  +c_2 \frac{b_1^{\prime} z_1+b_2^{\prime} z_2+b_3^{\prime}}
  {c_1^{\prime} z_1+c_2^{\prime} z_2+c_3^{\prime}}+c_3.$$
$$j(g_2, Z)=c_1^{\prime} z_1+c_2^{\prime} z_2+c_3^{\prime}.$$
So,
$$j(g_1, g_2(Z)) j(g_2, Z)=j(g_1 g_2, Z).$$
$\qquad \qquad \qquad \qquad \qquad \qquad \qquad \qquad \qquad
 \qquad \qquad \qquad \qquad \qquad \qquad \qquad \qquad \qquad
 \quad \boxed{}$

\vskip 0.5 cm
\centerline{\bf 3. Product formulas on $U(2, 1)$}
\vskip 0.5 cm

  In his paper \cite{Ko}, the Poisson kernel $P: B \times D 
\to {\Bbb R}$ of $D$ is defined by
$$P(u, z)=\frac{|S(u, z)|^{2}}{S(z, z)}$$
where $D=\{ (z_1, z_2) | \text{Im} z_1-\Phi(z_2, z_2) \in \Omega \}$ 
and $B$ is the distinguished boundary of $D$:
$$B=\{ (z_1, z_2) | \text{Im} z_1-\Phi(z_2, z_2)=0 \}.$$ 
The Szeg\"{o} kernel function $S$ is given by
$$S(z, w)=\frac{1}{(2 \pi)^{n_1}} \int_{\Omega^{\prime}}
  e^{-<\alpha, \rho(z, w)>} L(\alpha)^{-1} d \alpha,$$
where
$$\rho(z, w)=\frac{z_1-\overline{w_1}}{i}-2 \Phi(z_2, w_2)$$
for $z=(z_1, z_2)$, $w=(w_1, w_2)$, and
$$L(\alpha)=(\frac{i}{2})^{n_2} \int_{V_2}
  e^{-2 <\alpha, \Phi(z_2, z_2)>} dz_2 d\overline{z_2}.$$

  Now, we define the Poisson kernel as follows:
$$P(Z, W)=\frac{\rho(Z)}{|\rho(Z, W)|^2}, \quad \text{for} 
  \quad Z=(z_1, z_2) \in {\frak S}_2, W=(w_1, w_2) \in 
  \partial {\frak S}_2, \tag 3.1$$
where 
$$\rho(Z, W):=\overline{z_1}+w_1-\overline{z_2} w_2.$$    
In fact, for $Z=(\frac{\rho+|z|^2}{2}+it, z)$ and 
$Z^{\prime}=(\frac{\rho^{\prime}+|z^{\prime}|^2}{2}+it^{\prime},
z^{\prime})$, where $\rho>0$, $\rho^{\prime}>0$, 
$t, t^{\prime} \in {\Bbb R}$ and $z, z^{\prime} \in {\Bbb C}$, we have
$$|\rho(Z, Z^{\prime})|^2=\frac{1}{4}(\rho+\rho^{\prime}
 +|z-z^{\prime}|^2)^2+(t^{\prime}-t+\text{Im} z 
 \overline{z^{\prime}})^2 \geq 0.$$
$\rho(Z, Z^{\prime})=0$ if and only if 
$Z=Z^{\prime} \in \partial {\frak S}_2$.

  In their papers \cite{Ko}, \cite{KW}, Kor\'{a}nyi and Wolf 
studied the realization of Siegel domain of type II and the 
Poisson integral. Our definition has a little difference 
from theirs. We have 
$$L P(Z, W)^s=s(s-2) P(Z, W)^s-s^2 \rho(W) P(Z, W)^{s+1}
             =s(s-2) P(Z, W)^s. \tag 3.2$$ 
In fact, by Helgason's conjecture, which was proved by Kashiwara 
et al. in \cite{K}, that the eigenfunctions on Riemannian
symmetric spaces can be represented as Poisson integrals
of their hyperfunction boundary values.

  The boundary of ${\frak S}_2$:
$$\partial {\frak S}_2=\{ (z_1, z_2) \in {\Bbb C}^2:
  z_1+\overline{z_1}=z_2 \overline{z_2} \}. \tag 3.3$$
In the homogeneous coordinates, one has
$$\overline{\partial {\frak S}_2}=\{ (z_1, z_2, z_3) \in {\Bbb C}^3:
  z_1 \overline{z_3}+\overline{z_1} z_3=z_2 \overline{z_2} \}.$$
In particular,
$$\overline{\partial {\frak S}_2}({\Cal O})=\{ (a_1, a_2, a_3) \in
  {\Cal O}^3: a_1 \overline{a_3}+\overline{a_1} a_3=a_2 
  \overline{a_2} \}.$$

  The Green function associated to $L f=0$ on ${\frak S}_2$ 
is defined as
$$G(Z, W)=\log |\rho(Z, W)| \quad \text{for} \quad Z, W 
  \in {\frak S}_2. \tag 3.4$$
It satisfies that $L G(Z, W)=0$.

  If $g=\left(\matrix
     a_1 & a_2 & a_3\\
     b_1 & b_2 & b_3\\
     c_1 & c_2 & c_3 
     \endmatrix\right) \in U(2, 1)$, 
then $g^{*} J g=J$, i.e.,
$$\left\{\aligned
 -\overline{c_1} a_1+\overline{b_1} b_1-\overline{a_1} c_1 &=0,\\
 -\overline{c_2} a_1+\overline{b_2} b_1-\overline{a_2} c_1 &=0,\\
 -\overline{c_2} a_2+\overline{b_2} b_2-\overline{a_2} c_2 &=1,\\
 -\overline{c_3} a_1+\overline{b_3} b_1-\overline{a_3} c_1 &=-1,\\
 -\overline{c_3} a_2+\overline{b_3} b_2-\overline{a_3} c_2 &=0,\\
 -\overline{c_3} a_3+\overline{b_3} b_3-\overline{a_3} c_3 &=0.
\endaligned\right.$$
By a straightforward calculation, we have
$$\rho(g(Z), g(W))=\frac{\overline{z_1}+w_1-
  \overline{z_2} w_2}{\overline{(c_1 z_1+c_2 z_2+c_3)}
  (c_1 w_1+c_2 w_2+c_3)}.$$
Hence,
$$\rho(g(Z), g(W)) \overline{j(g, Z)} j(g, W)
 =\rho(Z, W), \quad \text{for} \quad g \in G. \tag 3.5$$			
Therefore,
$$P(g(Z), g(W))=|j(g, W)|^2 P(Z, W). \tag 3.6$$
By Proposition 2.9, we have
$$P(g(Z), W)=|j(g, g^{-1}(W))|^2 P(Z, g^{-1}(W))
 =|j(g^{-1}, W)|^{-2} P(Z, g^{-1}(W)). \tag 3.7$$

  The point-pair invariant on ${\frak S}_2$ is defined by
$$\sigma(Z, Z^{\prime}):=\frac{|\rho(Z, Z^{\prime})|^2}{\rho(Z)
  \rho(Z^{\prime})}, \tag 3.8$$
which is invariant under the action of $U(2, 1)$. The critical 
exponent is given by
$$\delta(\Gamma)=\text{inf} \{ s>0: \sum_{\gamma \in \Gamma}
  \sigma(Z, \gamma(Z^{\prime}))^{-s}<\infty \}. \tag 3.9$$
The number $\delta(\Gamma)$ does not depend on the choice of $Z$,
$Z^{\prime}$. It is an invariant of the group $\Gamma$. 

  The Eisenstein series is given by 			
$$E(Z, W; s):=\sum_{\gamma \in \Gamma} P(\gamma(Z), W)^s,
  \quad \text{for} \quad \text{Re}(s)>\delta(\Gamma),
  \tag 3.10$$
where $\Gamma$ is a discrete subgroup of $U(2, 1)$.
We have the following lemma:

{\smc Lemma 3.1}. {\it The properties of Eisenstein series: 
\roster
\item $L E(Z, W; s)=s(s-2) E(Z, W; s)$. 
\item
$E(\gamma(Z), W; s)=E(Z, W; s)$.
\item
$E(Z, \gamma(W); s)=|j(\gamma, W)|^{2s} E(Z, W; s)$ 
 for $\gamma \in \Gamma$.
\endroster}

  We have the following lemma:

{\smc Lemma 3.2}. {\it Let $f=f(u)$ with 
$u=u(Z, Z^{\prime})=\sigma(Z, Z^{\prime})-1$, then
$$Lf=u(u+1) f^{\prime \prime}(u)+(3u+2) f^{\prime}(u). \tag 3.11$$}

  By $Lf=s(s-2)f$, we have	 
$$u(u+1) f^{\prime \prime}(u)+(3u+2) f^{\prime}(u)+s(2-s) f(u)=0. 
  \tag 3.12$$	 
A solution is	 
$$\varphi_s(u(Z, Z^{\prime}))=\frac{\Gamma(s)}{\sqrt{\pi} \Gamma(s-
  \frac{1}{2})} (\frac{1}{4})^{s-1} u^{-s} F(s, s-1; 2s-1; -u^{-1}).
  \tag 3.13$$	 
Here $F={}_{2}F_{1}$ is the Gauss hypergeometric function. 

  Let
$$G(Z, Z^{\prime}; s)=\sum_{\gamma \in \Gamma} 
  r(Z, \gamma(Z^{\prime}); s), \quad \text{for} \quad
  \text{Re}(s)>\delta(\Gamma), \tag 3.14$$
where $r(Z, Z^{\prime}; s)=\varphi_s(u(Z, Z^{\prime}))$. 
By the integral representation of hypergeometric functions,
we can only consider the series
$$\sum_{\gamma \in \Gamma} \frac{1}{[1+u(Z, \gamma(Z^{\prime}))]^{s}}
 =\sum_{\gamma \in \Gamma} \sigma(Z, \gamma(Z^{\prime}))^{-s},$$
which is convergent for $\text{Re}(s)>\delta(\Gamma)$.
Thus, the right hand side of (3.14) is convergent when 
$\text{Re}(s)>\delta(\Gamma)$. We have
$$\lim_{\rho(Z) \to 0} \rho(Z)^{-s} r(Z, Z^{\prime}; s)=
  \frac{\Gamma(s)}{\sqrt{\pi} \Gamma(s-\frac{1}{2})} 
  (\frac{1}{4})^{s-1} P(Z^{\prime}, W)^{s},$$
where $Z \to W \in \partial {\frak S}_2$ as $\rho(Z) \to 0$.
Consequently,
$$\lim_{\rho(Z) \to 0} \rho(Z)^{-s} G(Z, Z^{\prime}; s)=
  \frac{\Gamma(s)}{\sqrt{\pi} \Gamma(s-\frac{1}{2})} (\frac{1}{4})^{s-1} 
  E(Z^{\prime}, W; s).$$

  In the real hyperbolic geometry, the double transitivity of
the action of the isometric groups on the real hyperbolic space
has been studied extensively (see \cite{ElGM} for three 
dimensional hyperbolic space). Now, we give the corresponding
results on the complex hyperbolic space.   

{\smc Proposition 3.3}. {\it The group $U(2, 1)$ acts in the 
following sense doubly transitively on ${\frak S}_{2}$: For all 
$P, P^{\prime}, Q, Q^{\prime} \in {\frak S}_{2}$, such that 
$d(P, P^{\prime})=d(Q, Q^{\prime})$, there exists an element 
$g \in U(2, 1)$ such that $g(P)=Q$, $g(P^{\prime})=Q^{\prime}$.}

{\it Proof}. Set $P=(z_1, z_2)$ with $z_1=-\omega \rho+\frac{|z|^2}{2}
+it$, $z_2=z$, where $\rho>0$, $t \in {\Bbb R}$ and $z \in {\Bbb C}$. 
Then $\rho(z_1, z_2)=\rho$.
$$T_{1}:=\left(\matrix
         \frac{1}{\sqrt{\rho}} &   &  \\
		                       & 1 &  \\  
                               &   & \sqrt{\rho} 
         \endmatrix\right)
		 \left(\matrix
		 1 & -\overline{z} & \frac{|z|^2}{2}-it\\
		   &             1 &                 -z\\
		   &               &                  1
		 \endmatrix\right) \in U(2, 1).$$
Moreover, $T_1(z_1, z_2)=(-\omega, 0)$, i.e. $T_1$ maps $P$ onto
$(-\omega, 0)$. Set $(z_1^{\prime}, z_2^{\prime})=(-\omega \rho^{\prime}
+\frac{|z^{\prime}|^2}{2}+i t^{\prime}, z^{\prime})$, then 
$$T_1(P^{\prime})=\left(-\omega \frac{\rho^{\prime}}{\rho}+\frac{1}{2}
  \left(\frac{|z^{\prime}-z|}{\sqrt{\rho}}\right)^{2}+i
  \frac{t^{\prime}-t+\text{Im}(z \overline{z^{\prime}})}{\rho},
  \frac{z^{\prime}-z}{\sqrt{\rho}}\right),$$
i.e., $T_1(P^{\prime})$ has the form of $(-\omega \rho+\frac{|z|^2}{2}
+it, z)$.

  If $k(-\omega \rho+\frac{|z|^2}{2}+it, z)=(-\omega \rho^{\prime}, 0)$,
where $k \in K$, then we obtain two complex equations about $\alpha,
\beta, \gamma, \delta$ and $\theta$. In fact, the number of real linear 
equations is four. While, the real dimension of $K$ is five. Hence, 
the existence of $k$ is proved. Applying a suitable element $T_2$ 
of the stabilizer $K$ of $(-\omega, 0)$ in $U(2, 1)$ to 
$T_1(P^{\prime})$ we get a transform $T=T_2 T_1$ such that 
$T(P)=(-\omega, 0)$, $T(P^{\prime})=(-\lambda \omega, 0)$ with 
$\lambda \geq 1$. Since the isometric group of ${\frak S}_{2}$ 
is $U(2, 1)$, we have
$$d(P, P^{\prime})=d(T(P), T(P^{\prime}))=d((-\omega, 0),
  (-\lambda \omega, 0))=\int_{1}^{\lambda} \frac{dr}{r}=
  \log \lambda.$$
So $\lambda=\exp(d(P, P^{\prime}))$. Similarly, there exists an 
element $S \in U(2, 1)$ such that $S(Q)=(-\omega, 0)$, 
$S(Q^{\prime})=(-\lambda \omega, 0)$ with the 
same $\lambda$, since $d(P, P^{\prime})=d(Q, Q^{\prime})$. Thus,
$g:=S^{-1} T$ has the desired properties.
\flushpar
$\qquad \qquad \qquad \qquad \qquad \qquad \qquad \qquad \qquad
 \qquad \qquad \qquad \qquad \qquad \qquad \qquad \qquad \qquad
 \quad \boxed{}$

  In \cite{ElGM}, the concept of real hyperbolic distance was
introduced and the basic properties was studied. Now, we give
the complex hyperbolic distance and obtain the triangle
inequality for our function $\delta$.  
  
{\smc Proposition 3.4}. {\it The complex hyperbolic distance 
$d(P, P^{\prime})$ is given by
$$\cosh d(P, P^{\prime})=\delta(P, P^{\prime}), \tag 3.15$$
where $\delta$ is defined by
$$\delta(P, P^{\prime})=\delta(Z, Z^{\prime})
 =\frac{1}{2}[\sigma(Z, Z^{\prime})+1],\tag 3.16$$
and $P=Z$, $P^{\prime}=Z^{\prime}$.}

{\it Proof}. If $P=(-\omega, 0)$ and $P^{\prime}=(-\lambda \omega, 0)$ 
$(\lambda \geq 1)$, we have
$$d((-\omega, 0), (-\lambda \omega, 0))=\int_{1}^{\lambda}
  \frac{dr}{r}=\log \lambda.$$
$$\delta((-\omega, 0), (-\lambda \omega, 0))=\frac{1}{2}(\lambda
  +\frac{1}{\lambda})=\cosh d((-\omega, 0), (-\lambda \omega, 0)).$$
Thus the proposition is true in the special case $P=(-\omega, 0)$,
$P^{\prime}=(-\lambda \omega, 0)$. Note that $\delta$ is a point-pair
invariant. Since $U(2, 1)$ acts doubly transitively on ${\frak S}_{2}$,
there exists for all $P, P^{\prime} \in {\frak S}_{2}$ an element
$g \in U(2, 1)$, such that $g(P)=(-\omega, 0)$, $g(P^{\prime})=
(-\lambda \omega, 0)$, $\lambda=\exp(d(P, P^{\prime}))$. Therefore,
by point-pair invariance, we have
$$\aligned
  \cosh d(P, P^{\prime})
=&\cosh d((-\omega, 0), (-\lambda \omega, 0))
 =\delta((-\omega, 0), (-\lambda \omega, 0))
 =\delta(g(P), g(P^{\prime}))\\
=&\delta(P, P^{\prime}).
\endaligned$$
$\qquad \qquad \qquad \qquad \qquad \qquad \qquad \qquad \qquad
 \qquad \qquad \qquad \qquad \qquad \qquad \qquad \qquad \qquad
 \quad \boxed{}$

{\smc Proposition 3.5}. {\it For $g=\left(\matrix
a_1 & a_2 & a_3\\
b_1 & b_2 & b_3\\
c_1 & c_2 & c_3
\endmatrix\right) \in U(2, 1)$, the following formulas hold:
\roster
\item
$$\aligned
 &\delta(P, g(P))+\delta(Q, g(Q))+\delta(P, g(Q))+\delta(Q, g(P))\\
=&\frac{5}{2}+2(|a_1|^2+|a_3|^2+|c_1|^2+|c_3|^2)+
  (|a_2|^2+|c_2|^2+|b_1|^2+|b_3|^2)+\frac{1}{2} |b_2|^2,
\endaligned\tag 3.17$$
where $P=(-\omega, 0)$ and $Q=(-\overline{\omega}, 0)$.
\item
$$\delta((1, 0), g(1, 0))=\frac{5}{8}+\frac{1}{8} \sum_{j=1}^{3}
  (|a_j|^2+|b_j|^2+|c_j|^2).\tag 3.18$$
\endroster}

{\it Proof}. It is obtained by a straightforward calculation.
\flushpar
$\qquad \qquad \qquad \qquad \qquad \qquad \qquad \qquad \qquad
 \qquad \qquad \qquad \qquad \qquad \qquad \qquad \qquad \qquad
 \quad \boxed{}$

{\smc Proposition 3.6}. {\it The triangle inequality for the metric 
function $\delta$ on the complex hyperbolic space: For all $P, Q,
R \in {\frak S}_{2}$, the following inequality holds:
$$\frac{1}{72} \frac{\delta(Q, R)}{\delta(P, Q)} \leq
  \delta(P, R) \leq 72 \delta(P, Q) \delta(Q, R).\tag 3.19$$}

{\it Proof}. Since $\delta$ is a point-pair invariant, it is
sufficient to prove (3.19) only in the special
case $Q=(\frac{1}{2}, 0)$. Set $P=(\frac{\rho_1+|z_1|^2}{2}+
it_1, z_1)$ and $R=(\frac{\rho_2+|z_2|^2}{2}+i t_2, z_2)$. Then
$$\delta(P, Q)=\frac{\frac{1}{4}(1+\rho_1+|z_1|^2)^2+t_1^2+\rho_1}
  {2 \rho_1}, \quad
  \delta(Q, R)=\frac{\frac{1}{4}(1+\rho_2+|z_2|^2)^2+t_2^2+\rho_2}
  {2 \rho_2},$$ 
and
$$\delta(P, R)=\frac{1}{8 \rho_1 \rho_2} [(\rho_1+\rho_2+|z_1-
  z_2|^2)^2+4(t_2-t_1+\text{Im} z_1 \overline{z_2})^2+4 \rho_1
  \rho_2].$$
Note that $|z_1-z_2|^2 \leq 2(|z_1|^2+|z_2|^2)$ and
$|\text{Im} z_1 \overline{z_2}| \leq |z_1| |z_2| \leq \frac{1}{2}
(|z_1|^2+|z_2|^2)$, we have
$$8 \rho_1 \rho_2 \delta(P, R) \leq (\rho_1+\rho_2+2 |z_1|^2+
  2 |z_2|^2)^2+(2 |t_1|+2 |t_2|+|z_1|^2+|z_2|^2)^2+4 \rho_1 \rho_2.$$
By $2 |t_1| |t_2| \leq t_1^2+t_2^2$ and $2 |t_i| |z_j|^2 \leq
t_i^2+|z_j|^4$, $i, j=1, 2$, we have
$$\aligned
     &8 \rho_1 \rho_2 \delta(P, R)\\
\leq &\rho_1^2+\rho_2^2+9 |z_1|^4+9 |z_2|^4+6 \rho_1 \rho_2
      +4(\rho_1+\rho_2)(|z_1|^2+|z_2|^2)+10 |z_1|^2 |z_2|^2\\
	 &+12 t_1^2+12 t_2^2\\
\leq &9(1+\rho_1^2+|z_1|^4+6 \rho_1+2 |z_1|^2+2 \rho_1 |z_1|^2
      +4 t_1^2)\\
	 &\times (1+\rho_2^2+|z_2|^4+6 \rho_2+2 |z_2|^2+2 \rho_2 
	  |z_2|^2+4 t_2^2)\\
    =&576 \delta(P, Q) \delta(Q, R).
\endaligned$$
Hence, $\delta(P, R) \leq 72 \delta(P, Q) \delta(Q, R)$.
Interchanging $P$ and $Q$, we obtain the other part of the
inequality.
\flushpar
$\qquad \qquad \qquad \qquad \qquad \qquad \qquad \qquad \qquad
 \qquad \qquad \qquad \qquad \qquad \qquad \qquad \qquad \qquad
 \quad \boxed{}$

{\smc Theorem 3.7}. {\it Let
$$\phi(Z, Z^{\prime}; s)=\int_{\partial {\frak S}_2} P(Z, W)^{s}
  P(Z^{\prime}, W)^{2-s} dm(W). \tag 3.20$$
Then one has the following formula:
$$\phi(Z, Z^{\prime}; s)=-\frac{\pi}{s-1} \frac{\sqrt{\pi} 
  \Gamma(\frac{3}{2}-s)}{\Gamma(2-s)} (\frac{1}{4})^{s-1} u^{-s} 
  F(s, s-1; 2s-1; -u^{-1})+(s \mapsto 2-s), \tag 3.21$$  
where $dm(W)$ is the Lebesque measure on $\partial {\frak S}_2$,
$u=u(Z, Z^{\prime})$ is the point-pair invariant, and $(s \mapsto 2-s)$
is equal to the first term on the right hand side with $s$ replaced
by $2-s$.}

{\it Proof}. $\phi(Z, Z^{\prime}; s)$ is a point-pair invariant. 
Thus, it suffices to calculate it for 
$Z=(z, t, \rho)=(0, 0, \rho)$ and $Z^{\prime}=(z^{\prime}, t^{\prime},
\rho^{\prime})=(0, 0, \rho^{\prime})$, i.e., $Z=(\frac{1}{2} \rho, 0)$ 
and $Z^{\prime}=(\frac{1}{2} \rho^{\prime}, 0)$. 
Let $\rho$ and $\rho^{\prime}$ stand for
$(\rho, 0)$ and $(\rho^{\prime}, 0)$, respectively. Now, 
$\rho(Z, Z^{\prime})=\frac{1}{2}(\rho+\rho^{\prime})$ and
$$\sigma(Z, Z^{\prime})=\frac{(\rho+\rho^{\prime})^2}{4 \rho 
  \rho^{\prime}}=\frac{(\lambda+1)^2}{4 \lambda}, 
  \quad u=\frac{(\lambda-1)^2}{4 \lambda},$$
where $\lambda=\rho/\rho^{\prime}$. 
For $W \in \partial {\frak S}_2$, $w_1=\frac{1}{2}
|w|^2+iv$ and $w_2=w$, where $w \in {\Bbb C}$, $v \in {\Bbb R}$.
$\rho(Z, W)=\frac{1}{2}(\rho+|w|^2)+iv$ and
$\rho(Z^{\prime}, W)=\frac{1}{2}(\rho^{\prime}+|w|^{2})+iv$. Hence,
$$\phi(\frac{1}{2} \rho, \frac{1}{2} \rho^{\prime}; s)
 =\int_{{\Bbb C}} \int_{{\Bbb R}} \frac{\rho^{s} 
  {\rho^{\prime}}^{2-s}}{|\frac{1}{2}(\rho+|w|^2)+iv|^{2s} 
  |\frac{1}{2}(\rho^{\prime}+|w|^2)+iv|^{2(2-s)}} dm(w) dv.$$
By transform $w=r e^{i \theta}$ with $r \geq 0, \theta \in [0, 2 \pi)$
and $u=r^2$, we have
$$\phi(\frac{1}{2} \rho, \frac{1}{2} \rho^{\prime}; s)=16 \pi \rho^{s} 
 {\rho^{\prime}}^{2-s} \int_{0}^{\infty} \int_{0}^{\infty} 
  \frac{du dv}{[(\rho+u)^2+v^2]^{s} [(\rho^{\prime}+u)^2+v^2]^{2-s}}.$$
Set $v=(\rho^{\prime}+u) \sqrt{w}$, then
$$\int_0^{\infty} \frac{dv}{[(\rho+u)^2+v^2]^{s} 
  [(\rho^{\prime}+u)^2+v^2]^{2-s}}=\frac{1}{2} (\rho^{\prime}+u)^{-3}
  \int_0^{\infty} w^{-\frac{1}{2}} (1+w)^{s-2} (l^2+w)^{-s} dw,$$
where $l=l(u)=\frac{\rho+u}{\rho^{\prime}+u}$.
 	 
  By (see \cite{Er}, p.115, (5))
$$F(a, b; c; 1-z)=\frac{\Gamma(c)}{\Gamma(b) \Gamma(c-b)}
  \int_0^{\infty} s^{b-1} (1+s)^{a-c} (1+sz)^{-a} ds$$
for $\text{Re}(c)>\text{Re}(b)>0$, $|\text{arg}(z)|<\pi$, we have
$$\int_0^{\infty} w^{-\frac{1}{2}} (1+w)^{s-2} (1+w l^{-2})^{-s} dw
 =\frac{\pi}{2} F(s, \frac{1}{2}; 2; 1-l^{-2}).$$
Now, we need the following two formulas:

(1) $z=1-l^{-2}$, $\frac{z}{z-1}=1-l^2$.

By (see \cite{Er}, p.105, 2.9.(3))
$$F(a, b; c; z)=(1-z)^{-a} F(a, c-b; c; \frac{z}{z-1}),$$
one has
$$F(s, \frac{1}{2}; 2; 1-l^{-2})=l^{2s} F(s, \frac{3}{2}; 2; 1-l^2).$$
  
(2) $z=1-l^2$, $1-z=l^2$.

By (see \cite{Er}, p.108, 2.10.(1))  
$$\aligned
  F(a, b; c; z)
 &=\frac{\Gamma(c) \Gamma(c-a-b)}{\Gamma(c-a)
  \Gamma(c-b)} F(a, b; a+b-c+1; 1-z)\\
 &+\frac{\Gamma(c) \Gamma(a+b-c)}{\Gamma(a) \Gamma(b)} 
  (1-z)^{c-a-b} F(c-a, c-b; c-a-b+1; 1-z),
\endaligned$$
where $|\arg(1-z)|<\pi$, one has
$$\aligned
  F(s, \frac{3}{2}; 2; 1-l^2)
 &=\frac{\Gamma(\frac{1}{2}-s)}{\sqrt{\pi} \Gamma(2-s)} 
  F(s, \frac{3}{2}; s+\frac{1}{2}; l^2)\\
 &+\frac{2 \Gamma(s-\frac{1}{2})}{\sqrt{\pi} \Gamma(s)} l^{1-2s}
  F(2-s, \frac{1}{2}; \frac{3}{2}-s; l^2).
\endaligned$$
Hence, $\phi(\frac{1}{2} \rho, \frac{1}{2} \rho^{\prime}; s)=
\phi_{I}+\phi_{II}$, where
$$\phi_{I}=4 \pi^2 \lambda^{s} {\rho^{\prime}}^{2}
  \frac{\Gamma(\frac{1}{2}-s)}{\sqrt{\pi} \Gamma(2-s)} 
  \int_0^{\infty} (\rho^{\prime}+u)^{-3} F(s, \frac{3}{2};
  s+\frac{1}{2}; l^2) du,$$
$$\phi_{II}=4 \pi^2 \lambda^{s} {\rho^{\prime}}^{2}
  \frac{2 \Gamma(s-\frac{1}{2})}{\sqrt{\pi} \Gamma(s)}
  \int_0^{\infty} (\rho^{\prime}+u)^{-3} l^{1-2s} 
  F(2-s, \frac{1}{2}; \frac{3}{2}-s; l^2) du.$$

  $l=\frac{\rho+u}{\rho^{\prime}+u}$ implies that 
$u=\frac{\rho-\rho^{\prime} l}{l-1}$. Thus, 
$du=\frac{\rho^{\prime}-\rho}{(l-1)^2} dl$ and 
$\rho^{\prime}+u=\frac{\rho-\rho^{\prime}}{l-1}$.
Without loss of generality, we can assume that $\rho<\rho^{\prime}$.
Then
$$\phi_{I}=-\frac{4 \pi^2 \lambda^s}{(1-\lambda)^2}
           \frac{\Gamma(\frac{1}{2}-s)}{\sqrt{\pi} \Gamma(2-s)}
		   \int_{\lambda}^{1} (l-1) F(s, \frac{3}{2};
		   s+\frac{1}{2}; l^2) dl,$$
$$\phi_{II}=-\frac{4 \pi^2 \lambda^s}{(1-\lambda)^2}
            \frac{2 \Gamma(s-\frac{1}{2})}{\sqrt{\pi} \Gamma(s)}
			\int_{\lambda}^{1} (l-1) l^{1-2s} F(2-s,
			\frac{1}{2}; \frac{3}{2}-s; l^2) dl.$$
We have			
$$\aligned
 &\int_{\lambda}^{1} (l-1) F(s, \frac{3}{2}; s+\frac{1}{2};
  l^2) dl\\
=&\sum_{k=0}^{\infty} \frac{(s)_{k} (\frac{3}{2})_{k}}{k! (s+
  \frac{1}{2})_{k}}(\frac{\lambda^{2k+1}}{2k+1}-\frac{\lambda^{2k+2}}
  {2k+2}-\frac{1}{(2k+1)(2k+2)})\\   
=&\lambda F(s, \frac{1}{2}; s+\frac{1}{2}; \lambda^2)-
  \frac{2s-1}{2s-2} F(s-1, \frac{1}{2}; s-\frac{1}{2};
  \lambda^2)+\frac{1}{s-1} \frac{\Gamma(s+\frac{1}{2})}{\sqrt{\pi} 
  \Gamma(s)}. 
\endaligned$$
Here we use the following formula:
$$F(s-1, -\frac{1}{2}; s-\frac{1}{2}; 1)=\frac{\Gamma(s-
  \frac{1}{2})}{\sqrt{\pi} \Gamma(s)}, \quad \text{for} \quad
  \text{Re}(s)>0.$$
Similarly,
$$\aligned
  &\int_{\lambda}^{1} (l-1) l^{1-2s} F(2-s, \frac{1}{2};
   \frac{3}{2}-s; l^2) dl\\
= &\frac{\lambda^{2-2s}}{2-2s} F(1-s, \frac{1}{2}; 
   \frac{3}{2}-s; \lambda^2)-\frac{\lambda^{3-2s}}{3-2s}
   F(2-s, \frac{1}{2}; \frac{5}{2}-s; \lambda^2)-
   \frac{1}{2-2s} \frac{\Gamma(\frac{3}{2}-s)}{\sqrt{\pi} \Gamma(2-s)},
\endaligned$$
for $\text{Re}(s)<2$.

  Set $z=-\lambda$, then $\frac{4z}{(1+z)^2}=-u^{-1}$.
By (see \cite{Er}, p.111, 2.11.(5))
$$F(a, b; 2b; \frac{4z}{(1+z)^2})=(1+z)^{2a}
  F(a, a+\frac{1}{2}-b; b+\frac{1}{2}; z^2),$$
we have
$$F(s, \frac{1}{2}; s+\frac{1}{2}; \lambda^2)=(1-\lambda)
 ^{-2s} F(s, s; 2s; -u^{-1}),$$
$$F(s-1, \frac{1}{2}; s-\frac{1}{2}; \lambda^2)=(1-\lambda)
 ^{2-2s} F(s-1, s-1; 2s-2; -u^{-1}),$$
$$F(1-s, \frac{1}{2}; \frac{3}{2}-s; \lambda^2)=(1-\lambda)
 ^{2s-2} F(1-s, 1-s; 2-2s; -u^{-1}),$$
$$F(2-s, \frac{1}{2}; \frac{5}{2}-s; \lambda^2)=(1-\lambda)
 ^{2s-4} F(2-s, 2-s; 4-2s; -u^{-1}).$$
Thus, 
$$\aligned
  \phi_{I}=
 &-4 \pi^2 \frac{\Gamma(\frac{1}{2}-s)}{\sqrt{\pi} \Gamma(2-s)}
  (4u)^{-s-1} F(s, s; 2s; -u^{-1})\\
 &-\frac{4 \pi^2}{s-1} \frac{\Gamma(\frac{3}{2}-s)}{\sqrt{\pi} 
  \Gamma(2-s)} (4u)^{-s} F(s-1, s-1; 2s-2; -u^{-1})\\
 &-\frac{4 \pi}{s-1} \frac{\Gamma(\frac{1}{2}-s) \Gamma(s+\frac{1}{2})}
  {\Gamma(s) \Gamma(2-s)} \lambda^{s} (1-\lambda)^{-2}, 
\endaligned$$
and
$$\aligned
  \phi_{II}=
 &-\frac{4 \pi^2}{1-s} \frac{\Gamma(s-\frac{1}{2})}{\sqrt{\pi} 
  \Gamma(s)} (4u)^{s-2} F(1-s, 1-s; 2-2s; -u^{-1})\\
 &-4 \pi^2 \frac{\Gamma(s-\frac{3}{2})}{\sqrt{\pi} \Gamma(s)}
  (4u)^{s-3} F(2-s, 2-s; 4-2s; -u^{-1})\\
 &-\frac{4 \pi}{s-1} \frac{\Gamma(s-\frac{1}{2}) \Gamma(\frac{3}{2}-s)} 
  {\Gamma(s) \Gamma(2-s)} \lambda^s (1-\lambda)^{-2}.
\endaligned$$
Note that
$$\Gamma(\frac{1}{2}-s) \Gamma(s+\frac{1}{2})+\Gamma(s-\frac{1}{2})
  \Gamma(\frac{3}{2}-s)=(s-\frac{1}{2}) \Gamma(\frac{1}{2}-s)
  \Gamma(s-\frac{1}{2})+(\frac{1}{2}-s) \Gamma(s-\frac{1}{2})
  \Gamma(\frac{1}{2}-s)=0.$$
Moreover, we need the following lemma.

{\smc Lemma 3.8}. {\it The following formula about hypergeometric
functions holds:
$$\gamma F-\beta z F(\beta+1, \gamma+1)-\gamma F(\alpha-1)=0,$$
where $F$ denotes $F(\alpha, \beta; \gamma; z)$, 
$F(\beta+1, \gamma+1)$ stands for $F(\alpha, \beta+1; \gamma+1; z)$
and $F(\alpha-1)$ stands for $F(\alpha-1, \beta; \gamma; z)$.}

{\it Proof}. It is obtained by the integral representation of
hypergeometric functions.
\flushpar
$\qquad \qquad \qquad \qquad \qquad \qquad \qquad \qquad \qquad
 \qquad \qquad \qquad \qquad \qquad \qquad \qquad \qquad \qquad
 \quad \boxed{}$

  By Lemma 3.8 and the following formula (see \cite{Er}, p.103, 
(35)): 
$$(\gamma-1) F(\gamma-1)-\alpha F(\alpha+1)-(\gamma-\alpha-1) F=0.$$
we get
$$(2s-1) F(s, s-1; 2s-1; z)-(s-1) z F(s, s; 2s; z)-(2s-1) F(s-1, s-1;
  2s-1; z)=0,$$
$$(2s-2) F(s-1, s-1; 2s-2; z)-(s-1) F(s, s-1; 2s-1; z)-(s-1) F(s-1,
  s-1; 2s-1; z)=0.$$
Thus, we have
$$\frac{4s-2}{s-1} F(s, s-1; 2s-1; z)=z F(s, s; 2s; z)+\frac{4s-2}{s-1}
  F(s-1, s-1; 2s-2; z).$$
Substituting this identity to the expression of 
$\phi(\frac{1}{2} \rho, \frac{1}{2} \rho^{\prime}; s)$, We get 
the desired results.
\flushpar
$\qquad \qquad \qquad \qquad \qquad \qquad \qquad \qquad \qquad
 \qquad \qquad \qquad \qquad \qquad \qquad \qquad \qquad \qquad
 \quad \boxed{}$

  Consequently, we have 
  
{\smc Corollary 3.9}. {\it The following formula holds:
$$\phi(Z, Z^{\prime}; s)=-\frac{\pi}{s-1} c(s) c(2-s)
  [r(Z, Z^{\prime}; s)-r(Z, Z^{\prime}; 2-s)], \tag 3.22$$
where 
$$c(s)=\frac{\sqrt{\pi} \Gamma(s-\frac{1}{2})}{\Gamma(s)}.$$}
  
{\smc Corollary 3.10}. {\it The functional equation for the trivial 
group on ${\frak S}_2$: If $\text{Re}(s)>1$, then the following 
identity holds:
$$\int_{\partial {\frak S}_2} P(Z, W)^{s} |\rho(W, W^{\prime})|^
  {-2(2-s)} dm(W)=\frac{\pi}{s-1} \frac{\sqrt{\pi} 
  \Gamma(s-\frac{1}{2})}{\Gamma(s)} 4^{s-1} P(Z, W^{\prime})^{2-s},
  \tag 3.23$$
where $W^{\prime} \in \partial {\frak S}_2$.}

{\it Proof}. The integral on the left-hand side is absolutely
convergent. We set
$$P(Z^{\prime}, W)=\frac{\rho(Z^{\prime})}{|\rho(Z^{\prime}, W)|^2}$$
in the formula of Theorem 3.7. Next, we multiply both sides of
the formula by $\rho(Z^{\prime})^{s-2}$ and take the limit as
$Z^{\prime} \to W^{\prime} \in \partial {\frak S}_2$. Thus we get
the desired results.
\flushpar
$\qquad \qquad \qquad \qquad \qquad \qquad \qquad \qquad \qquad
 \qquad \qquad \qquad \qquad \qquad \qquad \qquad \qquad \qquad
 \quad \boxed{}$

{\smc Proposition 3.11}. {\it The integral formula for the Poisson 
kernel:
$$\int_{\partial {\frak S}_2} P(Z, W)^{s} dm(W)=\frac{\pi}{s-1}
  \frac{\sqrt{\pi} \Gamma(s-\frac{1}{2})}{\Gamma(s)} 4^{s-1}
  \rho^{2-s}, \quad \text{for} \quad \text{Re}(s)>1.
  \tag 3.24$$}

{\it Proof}. Denote the left hand side of (3.24) as $f(Z)$,
then 
$$f(Z)=\rho^s \int_{\partial {\frak S}_2} |\overline{z_1}+w_1-
  \overline{z_2} w_2|^{-2 s} dm(W).$$
Set $z_1=\frac{\rho+|z|^2}{2}+ir$, $z_2=z$, $w_1=\frac{|w|^2}{2}+it$ and
$w_2=w$, where $z, w \in {\Bbb C}$, $t, r \in {\Bbb R}$, we have
$$\aligned
  f(Z)
 &=\rho^s \int_{{\Bbb C}} \int_{{\Bbb R}} |\frac{\rho+|w-z|^2}{2}+
   i(t-r+\text{Im}(z \overline{w}))|^{-2s} dm(w) dt\\
 &=\rho^s \int_{{\Bbb C}} \int_{{\Bbb R}} |\frac{\rho+|w|^2}{2}+
   it|^{-2s} dm(w) dt\\ 
 &=2 \pi \rho^{s} \int_{0}^{\infty} \int_{-\infty}^{\infty}
   \frac{r dr dt}{[\frac{1}{4}(\rho+r^2)^{2}+t^2]^{s}}\\
 &=2^{2s-1} \pi \rho^{s} \int_{0}^{\infty} \int_{-\infty}^{\infty}
   \frac{dr du}{[(\rho+r)^{2}+u^{2}]^{s}}.
\endaligned$$
By the well known formula (see \cite{Ku}, p.15):
$$\int_{-\infty}^{\infty} \frac{dt}{(1+t^2)^{s}}=\sqrt{\pi}
  \frac{\Gamma(s-\frac{1}{2})}{\Gamma(s)}, \quad \text{for} \quad
  \text{Re}(s)>1,$$
we have
$$f(Z)=\frac{\pi}{s-1} \frac{\sqrt{\pi} \Gamma(s-\frac{1}{2})}
       {\Gamma(s)} 4^{s-1} \rho^{2-s}.$$
$\qquad \qquad \qquad \qquad \qquad \qquad \qquad \qquad \qquad
 \qquad \qquad \qquad \qquad \qquad \qquad \qquad \qquad \qquad
 \quad \boxed{}$

  We define the S-matrix as follows:
$$S(W, W^{\prime}; s):=\sum_{\gamma \in \Gamma} |\rho(\gamma(W), 
  W^{\prime})|^{-2s} |j(\gamma, W)|^{-2s} \tag 3.25$$  
for $W, W^{\prime} \in \Omega(\Gamma)$ and 
$\text{Re}(s)>\delta(\Gamma)$.

{\smc Theorem 3.12}. {\it Assume that $\Gamma$ is convex cocompact
and $\delta(\Gamma)<1$. Then the following functional equation for 
Eisenstein series holds:
$$\int_{\Gamma \backslash \Omega(\Gamma)} E(Z, W; s) S(W, W^{\prime}; 
  2-s) dm(W)=\frac{4^{s-1} \pi}{s-1} \frac{\sqrt{\pi} 
  \Gamma(s-\frac{1}{2})}{\Gamma(s)} E(Z, W^{\prime}; 2-s),
  \tag 3.26$$
for $1<\text{Re}(s)<2-\delta(\Gamma)$.}

{\it Proof}. By Corollary 3.10,
$$\int_{\partial {\frak S}_2} E(Z, W; s) |\rho(W, W^{\prime})|^{-2(2-s)}
  dm(W)=\frac{4^{s-1} \pi}{s-1} \frac{\sqrt{\pi} \Gamma(s-\frac{1}{2})}
  {\Gamma(s)} E(Z, W^{\prime}; 2-s).\tag 3.27$$
On the other hand,
$$\aligned
  \int_{\partial {\frak S}_2} f(W) dW
=&\int_{\Gamma \backslash \Omega(\Gamma)} \sum_{\gamma \in \Gamma} 
  f(\gamma(W)) dm(\gamma(W))\\
=&\int_{\Gamma \backslash \Omega(\Gamma)} \sum_{\gamma \in \Gamma} 
  |J(\gamma(W))| f(\gamma(W)) dm(W).
\endaligned$$

  Hence, the left hand side of (3.27):
$$\aligned
  LHS
=&\int_{\Gamma \backslash \Omega(\Gamma)} \sum_{\gamma \in \Gamma}
  |J(\gamma(W))| E(Z, \gamma(W); s) |\rho(\gamma(W), W^{\prime})|^
  {-2(2-s)} dm(W)\\
=&\int_{\Gamma \backslash \Omega(\Gamma)} \sum_{\gamma \in \Gamma}
  |j(\gamma, W)|^{-4} |j(\gamma, W)|^{2s} E(Z, W; s) |\rho(\gamma(W),
  W^{\prime})|^{-2(2-s)} dm(W)\\
=&\int_{\Gamma \backslash \Omega(\Gamma)} E(Z, W; s) S(W, W^{\prime};
  2-s) dm(W).
\endaligned$$  
$\qquad \qquad \qquad \qquad \qquad \qquad \qquad \qquad \qquad
 \qquad \qquad \qquad \qquad \qquad \qquad \qquad \qquad \qquad
 \quad \boxed{}$

{\smc Proposition 3.13}. {\it Assume that $\Gamma$ is convex
cocompact and $\delta(\Gamma)<1$. Then the following product 
formula holds:
$$\aligned
 &\int_{\Gamma \backslash \Omega(\Gamma)} E(Z, W; s) E(Z^{\prime}, W;
  2-s) dm(W)\\
=&-\frac{\pi}{s-1} c(s) c(2-s) [G(Z, Z^{\prime}; s)-
  G(Z, Z^{\prime}; 2-s)],
\endaligned\tag 3.28$$
where $\delta(\Gamma)<\text{Re}(s)<2-\delta(\Gamma)$.}  

  For $s \in {\Bbb C}$, similar as in \cite{Ma1} and \cite{Ma2}, we set 
$$L_{s}^{2}(\Gamma):=\{ f(s, \cdot): f(s, \cdot) \in C^{\infty}(\Omega
  (\Gamma)), f(s, \gamma(W)) j(\gamma, W)^{2-s}=f(s, W), \gamma \in 
  \Gamma \}.$$
If $f(s, \cdot) \in L_{s}^{2}(\Gamma)$, we define the Eisenstein
integral:
$$E(Z, \cdot, s): L_{s}^{2}(\Gamma) \to C^{\infty}(\Gamma \backslash
  {\frak S}_2),$$
$$E(Z, f, s)=\int_{\Gamma \backslash \Omega(\Gamma)} E(Z, W; s)
  f(s, W) dm(W), \quad \text{for} \quad \text{Re}(s)>\delta(\Gamma),
  \tag 3.29$$
and the scattering operator:
$$S_{s}(f)(z, t):=\lim_{Z \to (z, t)} \rho(Z)^{-s} [E(Z, f, s)-
  \frac{4^{s-1} \pi}{s-1} \frac{\sqrt{\pi} \Gamma(s-\frac{1}{2})}
  {\Gamma(s)} f(s, (z, t)) \rho(Z)^{2-s}],\tag 3.30$$
where $\delta(\Gamma)<\text{Re}(s)<2$.

  Now, the functional equation for Eisenstein series can be written in
terms of Eisenstein integrals as
$$E(Z, S_{2-s}(f), s)=\frac{4^{s-1} \pi}{s-1} \frac{\sqrt{\pi}
  \Gamma(s-\frac{1}{2})}{\Gamma(s)} E(Z, f, 2-s),\tag 3.31$$  
and in terms involving the scattering operator as:
$$S_{s} S_{2-s}=-\frac{\pi^3}{(s-1)^2} \frac{\Gamma(s-\frac{1}{2})
  \Gamma(\frac{3}{2}-s)}{\Gamma(s) \Gamma(2-s)} I.\tag 3.32$$

\vskip 0.5 cm
\centerline{\bf 4. Product formulas of weight $k$ on $U(2, 1)$}
\vskip 0.5 cm

  In the case of $SL(2, {\Bbb R})$, many people, especially Maass 
\cite{M}, Selberg \cite{Se}, Roelcke \cite{Ro}, Fay \cite{Fa}, 
Elstrodt \cite{El}, Hejhal \cite{H}, and Shimura \cite{Sh2} studied the 
eigenvalue problems which concern a homogeneous equation
$$(\Delta_{k}-\lambda) f(z)=0.$$
Here $f$ is a function on the complex upper plane
$${\Bbb H}=\{ z \in {\Bbb C}: \text{Im}(z)>0 \},$$  
and $L_{k}$ is a differential operator given by
$$\Delta_{k}=y^{2}\left(\frac{\partial^2}{\partial x^2}+
  \frac{\partial^2}{\partial y^2}\right)-iky 
  \frac{\partial}{\partial x},$$
where $z=x+iy$ and $k$ is an integer. Let $\Gamma$ be a
congruence subgroup of $SL(2, {\Bbb Z})$ and assume that
$f$ is a $\Gamma$-automorphic form of weight $k$ in the 
sense that it is invariant under
$$f \mapsto f(\gamma(z)) (cz+d)^{-k} |cz+d|^{k}$$ 
for any $\gamma=\left(\matrix
                       * & *\\
					   c & d 
                      \endmatrix\right) \in \Gamma$.
In fact, $L_{k}$ commutes with the above map.
One of the most remarkable facts about such $\lambda$ 
and $f$ is the existence of Selberg zeta function (see 
\cite{Se}), whose set of zeros coincides essentially with the 
set $S(\Gamma, k)$ of $s \in {\Bbb C}$ such that $\lambda=s(1-s)$
occurs as an eigenvalue of $L_{k}$ as in the above equation
with a cusp form $f$ (see \cite{Sh2}). 

  Now, we study the corresponding problems in the case 
of $U(2, 1)$.
     
{\smc Proposition 4.1}. {\it Let
$$\aligned
  L_{(\alpha, \beta)}
&=(z_1+\overline{z_1}-z_2 \overline{z_2})\\
&\times \left[ (z_1+\overline{z_1}) \frac{\partial^2}{\partial z_1
  \partial \overline{z_1}}+\frac{\partial^2}{\partial z_2
  \partial \overline{z_2}}+z_2 \frac{\partial^2}{\partial 
  \overline{z_1} \partial z_2}+\overline{z_2} \frac{\partial^2}
  {\partial z_1 \partial \overline{z_2}}+\beta \frac{\partial}
  {\partial z_1}+\alpha \frac{\partial}{\partial \overline{z_1}}
  \right]
\endaligned\tag 4.1$$
and
$$E(Z, \overline{Z}; \alpha, \beta)=\sum_{\gamma \in \Gamma \cap N
  \backslash \Gamma} j(\gamma, Z)^{-\alpha} \overline{j(\gamma, Z)}
  ^{-\beta}, \tag 4.2$$
where $\Gamma=U(2, 1; {\Cal O})$ is the Picard modular group.
Then the following formula holds:
$$L_{(\alpha, \beta)} E(Z, \overline{Z}; \alpha, \beta)=0.\tag 4.3$$}

{\it Proof}. We have
$$\aligned
 &\left[(z_1+\overline{z_1}) \frac{\partial^2}{\partial z_1 \partial
  \overline{z_1}}+\frac{\partial^2}{\partial z_2 \partial 
  \overline{z_2}}+z_2 \frac{\partial^2}{\partial \overline{z_1}
  \partial z_2}+\overline{z_2} \frac{\partial^2}{\partial z_1
  \partial \overline{z_2}} \right] [j(\gamma, Z)^{-\alpha}
  \overline{j(\gamma, Z)}^{-\beta}]\\
=&\alpha \beta [a_1 j(\gamma, Z)^{-\alpha-1} \overline{j(\gamma, 
  Z)}^{-\beta}+\overline{a_1} j(\gamma, Z)^{-\alpha} \overline{
  j(\gamma, Z)}^{-\beta-1}]\\
=&-\left(\beta \frac{\partial}{\partial z_1}+\alpha \frac{\partial}
  {\partial \overline{z_1}} \right) [j(\gamma, Z)^{-\alpha}
  \overline{j(\gamma, Z)}^{-\beta}].
\endaligned$$
Hence, $L_{(\alpha, \beta)} E(Z, \overline{Z}; \alpha, \beta)=0$.
\flushpar
$\qquad \qquad \qquad \qquad \qquad \qquad \qquad \qquad \qquad
 \qquad \qquad \qquad \qquad \qquad \qquad \qquad \qquad \qquad
 \quad \boxed{}$

  Set
$$E_{3k}(Z, \overline{Z}; \beta)=\rho(Z)^{\beta} E(Z, \overline{Z};
  \alpha, \beta)=\sum_{\gamma \in \Gamma \cap N \backslash \Gamma}
  \frac{\rho(Z)^{\beta}}{j(\gamma, Z)^{3k} |j(\gamma, Z)|^{2 \beta}}
  \tag 4.4$$
with $3k=\alpha-\beta$. 

  In particular, when $\beta=-k$ and $\alpha=k$, we have 
$$\aligned
  L_{k}
 &=(z_1+\overline{z_1}-z_2 \overline{z_2})\\
 &\times \left[(z_1+\overline{z_1}) \frac{\partial^2}{\partial z_1 
  \partial \overline{z_1}}+\frac{\partial^2}{\partial z_2 \partial
  \overline{z_2}}+z_2 \frac{\partial^2}{\partial \overline{z_1}
  \partial z_2}+\overline{z_2} \frac{\partial^2}{\partial z_1
  \partial \overline{z_2}}-k (\frac{\partial}{\partial z_1}
  -\frac{\partial}{\partial \overline{z_1}}) \right].
\endaligned$$  

{\smc Proposition 4.2}. {\it The operator $L_{k}$ has the following
automorphic property:
$$L_{k}[f(\gamma(Z)) j(\gamma, Z)^{-k} \overline{j(\gamma, Z)}^{k}]=
  j(\gamma, Z)^{-k} \overline{j(\gamma, Z)}^{k} (L_{k}f)(\gamma(Z)),
  \tag 4.5$$
for $\gamma \in U(2, 1)$ and $Z=(z_1, z_2) \in {\Bbb C}^{2}$, where 
$f=f(Z, \overline{Z})$ is a real analytic function. In other words,
$L_k$ commutes with the following map
$$f(Z) \mapsto f(\gamma(Z)) j(\gamma, Z)^{-k} 
  \overline{j(\gamma, Z)}^{k}.\tag 4.6$$}

{\it Proof}. 
For $\gamma=\left(\matrix
             a_1 & a_2 & a_3\\
			 b_1 & b_2 & b_3\\
			 c_1 & c_2 & c_3
			\endmatrix\right) \in U(2, 1)$, one has
$\gamma J \gamma^{*}=J$, i.e.,
$$\left\{\aligned
  -\overline{a_1} a_3+a_2 \overline{a_2}-a_1 \overline{a_3} &=0,\\
  -a_3 \overline{b_1}+a_2 \overline{b_2}-a_1 \overline{b_3} &=0,\\
  -a_3 \overline{c_1}+a_2 \overline{c_2}-a_1 \overline{c_3} &=-1,\\
  -b_3 \overline{b_1}+b_2 \overline{b_2}-b_1 \overline{b_3} &=1,\\
  -b_3 \overline{c_1}+b_2 \overline{c_2}-b_1 \overline{c_3} &=0,\\
  -c_3 \overline{c_1}+c_2 \overline{c_2}-c_1 \overline{c_3} &=0.
\endaligned\right.$$
Denote
$$\gamma(z_1, z_2)=(u(z_1, z_2), v(z_1, z_2))=\left( 
  \frac{a_1 z_1+a_2 z_2+a_3}{c_1 z_1+c_2 z_2+c_3},
  \frac{b_1 z_1+b_2 z_2+b_3}{c_1 z_1+c_2 z_2+c_3} \right).$$
Set
$$A_{1, 1}:=\frac{\partial}{\partial z_1} \left(
  \frac{a_1 z_1+a_2 z_2+a_3}{c_1 z_1+c_2 z_2+c_3} \right)
 =[(a_1 c_2-a_2 c_1) z_2+(a_1 c_3-a_3 c_1)] j(\gamma, Z)^{-2}.$$
$$A_{1, 2}:=\frac{\partial}{\partial z_2} \left( 
  \frac{a_1 z_1+a_2 z_2+a_3}{c_1 z_1+c_2 z_2+c_3} \right)
 =[(a_2 c_1-a_1 c_2) z_1+(a_2 c_3-a_3 c_2)] j(\gamma, Z)^{-2}.$$
$$A_{2, 1}:=\frac{\partial}{\partial z_1} \left( 
  \frac{b_1 z_1+b_2 z_2+b_3}{c_1 z_1+c_2 z_2+c_3} \right)
 =[(b_1 c_2-b_2 c_1) z_2+(b_1 c_3-b_3 c_1)] j(\gamma, Z)^{-2}.$$
$$A_{2, 2}:=\frac{\partial}{\partial z_2} \left(
  \frac{b_1 z_1+b_2 z_2+b_3}{c_1 z_1+c_2 z_2+c_3} \right)
 =[(b_2 c_1-b_1 c_2) z_1+(b_2 c_3-b_3 c_2)] j(\gamma, Z)^{-2}.$$
We have
$$\aligned
  L_k(\gamma(Z))
=&(u+\overline{u}-v \overline{v})\\
 &\times \left[(u+\overline{u}) \frac{\partial^2}{\partial u \partial
  \overline{u}}+\frac{\partial^2}{\partial v \partial \overline{v}}
  +v \frac{\partial^2}{\partial \overline{u} \partial v}+
  \overline{v} \frac{\partial^2}{\partial u \partial 
  \overline{v}}-k \left(\frac{\partial}{\partial u}-
  \frac{\partial}{\partial \overline{u}} \right) \right],
\endaligned$$
while $u+\overline{u}-v \overline{v}=(z_1+\overline{z_1}-z_2
\overline{z_2}) |j(\gamma, Z)|^{-2}$.

 On the other hand,
$$\aligned
 &L_k[f(\gamma(Z)) j(\gamma, Z)^{-k} \overline{j(\gamma, Z)}^{k}]\\
=&j(\gamma, Z)^{-k} \overline{j(\gamma, Z)}^{k} (z_1+\overline{z_1}-
  z_2 \overline{z_2})[F_{1, 1} \frac{\partial^2 f}{\partial u
  \partial \overline{u}}+F_{1, 2} \frac{\partial^2 f}{\partial u 
  \partial \overline{v}}+F_{2, 1} \frac{\partial^2 f}{\partial 
  \overline{u} \partial v}\\
 &+F_{2, 2} \frac{\partial^2 f}{\partial v \partial \overline{v}}+
  H_1 \frac{\partial f}{\partial u}+H_2 \frac{\partial f}{\partial v}
  +G_1 \frac{\partial f}{\partial \overline{u}}+G_2 \frac{\partial f}
  {\partial \overline{v}}](\gamma(Z)),
\endaligned$$
where
$$\left\{\aligned
 F_{1, 1} &=(z_1+\overline{z_1}) A_{1, 1} \overline{A_{1, 1}}+
 A_{1, 2} \overline{A_{1, 2}}+z_2 A_{1, 2} \overline{A_{1, 1}}+
 \overline{z_2} A_{1, 1} \overline{A_{1, 2}},\\
 F_{1, 2} &=(z_1+\overline{z_1}) A_{1, 1} \overline{A_{2, 1}}+
 A_{1, 2} \overline{A_{2, 2}}+z_2 A_{1, 2} \overline{A_{2, 1}}+
 \overline{z_2} A_{1, 1} \overline{A_{2, 2}},\\
 F_{2, 2} &=(z_1+\overline{z_1}) A_{2, 1} \overline{A_{2, 1}}+
 A_{2, 2} \overline{A_{2, 2}}+z_2 A_{2, 2} \overline{A_{2, 1}}+
 \overline{z_2} A_{2, 1} \overline{A_{2, 2}},\\
 H_1 &=-k A_{1, 1}+k[(z_1+\overline{z_1}) \overline{c_1} A_{1, 1}
 +\overline{c_2} A_{1, 2}+z_2 \overline{c_1} A_{1, 2}+\overline{z_2}
 \overline{c_2} A_{1, 1}] j(\overline{\gamma}, \overline{Z})^{-1},\\
 H_2 &=-k A_{2, 1}+k[(z_1+\overline{z_1}) \overline{c_1} A_{2, 1}
 +\overline{c_2} A_{2, 2}+z_2 \overline{c_1} A_{2, 2}+\overline{z_2}
 \overline{c_2} A_{2, 1}] j(\overline{\gamma}, \overline{Z})^{-1},\\
 F_{2, 1} &=\overline{F_{1, 2}}, \quad G_1=-\overline{H_1}, \quad 
 G_2=-\overline{H_2}.
\endaligned\right.$$
Now, we will prove that $F_{1, 1}=(u+\overline{u}) |j(\gamma, Z)|^{-2}$.
Note that
$$\aligned
 &F_{1, 1} |j(\gamma, Z)|^{4}\\
 &=(z_1+\overline{z_1})[(a_1 c_2-a_2 c_1) z_2+(a_1 c_3-a_3 c_1)]
  [(\overline{a_1} \overline{c_2}-\overline{a_2} \overline{c_1})
  \overline{z_2}+(\overline{a_1} \overline{c_3}-\overline{a_3}
  \overline{c_1})]\\
 &+[(a_2 c_1-a_1 c_2) z_1+(a_2 c_3-a_3 c_2)]
  [(\overline{a_2} \overline{c_1}-\overline{a_1} \overline{c_2})
  \overline{z_1}+(\overline{a_2} \overline{c_3}-\overline{a_3}
  \overline{c_2})]\\
 &+z_2 [(a_2 c_1-a_1 c_2) z_1+(a_2 c_3-a_3 c_2)]
  [(\overline{a_1} \overline{c_2}-\overline{a_2} \overline{c_1})
  \overline{z_2}+(\overline{a_1} \overline{c_3}-\overline{a_3}
  \overline{c_1})]\\
 &+\overline{z_2}[(a_1 c_2-a_2 c_1)z_2+(a_1 c_3-a_3 c_1)]
  [(\overline{a_2} \overline{c_1}-\overline{a_1} \overline{c_2})
  \overline{z_1}+(\overline{a_2} \overline{c_3}-\overline{a_3}
  \overline{c_2})]. 
\endaligned$$
On the other hand,
$$u+\overline{u}=[(a_1 z_1+a_2 z_2+a_3)(\overline{c_1} \overline{z_1}
  +\overline{c_2} \overline{z_2}+\overline{c_3})+(\overline{a_1}
  \overline{z_1}+\overline{a_2} \overline{z_2}+\overline{a_3})
  (c_1 z_1+c_2 z_2+c_3)] |j(\gamma, Z)|^{-2}.$$
Now, we need the following six identities:
$$\aligned
 &(a_2 c_1-a_1 c_2)(\overline{a_2} \overline{c_1}-\overline{a_1}
  \overline{c_2})\\
=&(\overline{a_1} a_3+a_1 \overline{a_3}) c_1 \overline{c_1}-a_1
  \overline{a_2} \overline{c_1} c_2-\overline{a_1} a_2 c_1
  \overline{c_2}+a_1 \overline{a_1} (c_1 \overline{c_3}+c_3
  \overline{c_1})\\
=&\overline{a_1} c_1(a_3 \overline{c_1}-a_2 \overline{c_2}+a_1
  \overline{c_3})+a_1 \overline{c_1}(\overline{a_3} c_1-
  \overline{a_2} c_2+\overline{a_1} c_3)\\  
=&\overline{a_1} c_1+a_1 \overline{c_1}.
\endaligned$$
$$\aligned
 &(a_1 c_2-a_2 c_1)(\overline{a_1} \overline{c_3}-\overline{a_3}
  \overline{c_1})\\
=&\overline{a_1} c_2(1-a_3 \overline{c_1}+a_2 \overline{c_2})+a_2
  \overline{c_1} (1+\overline{a_2} c_2-\overline{a_1} c_3)-
  \overline{a_1} a_2 c_1 \overline{c_3}-a_1 \overline{a_3}
  \overline{c_1} c_2\\
=&\overline{a_1} c_2+a_2 \overline{c_1}+\overline{c_1} c_2 (-
  \overline{a_1} a_3+a_2 \overline{a_2}-a_1 \overline{a_3})+
  \overline{a_1} a_2 (c_2 \overline{c_2}-\overline{c_1} c_3-
  c_1 \overline{c_3})\\
=&\overline{a_1} c_2+a_2 \overline{c_1}.      
\endaligned$$ 
$$\aligned
 &(a_2 c_3-a_3 c_2)(\overline{a_1} \overline{c_2}-\overline{a_2}
  \overline{c_1})+(a_1 c_2-a_2 c_1)(\overline{a_2} \overline{c_3}-
  \overline{a_3} \overline{c_2})\\
=&a_2 \overline{c_2} (\overline{a_1} c_3+\overline{a_3} c_1)+
  \overline{a_2} c_2 (a_3 \overline{c_1}+a_1 \overline{c_3})-
  \overline{a_1} a_3 c_2 \overline{c_2}-a_2 \overline{a_2}
  \overline{c_1} c_3-a_2 \overline{a_2} c_1 \overline{c_3}-
  a_1 \overline{a_3} c_2 \overline{c_2}\\
=&a_2 \overline{c_2}+\overline{a_2} c_2+c_2 \overline{c_2}(a_2
  \overline{a_2}-\overline{a_1} a_3-a_1 \overline{a_3})+a_2
  \overline{a_2}(c_2 \overline{c_2}-\overline{c_1} c_3-c_1
  \overline{c_3})\\
=&a_2 \overline{c_2}+\overline{a_2} c_2.
\endaligned$$ 
$$\aligned
 &(a_1 c_3-a_3 c_1)(\overline{a_1} \overline{c_3}-\overline{a_3}
  \overline{c_1})+(a_2 c_1-a_1 c_2)(\overline{a_2} \overline{c_3}-
  \overline{a_3} \overline{c_2})\\
=&a_1 \overline{c_3}(\overline{a_1} c_3-\overline{a_2} c_2)+
  \overline{a_3} c_1 (a_3 \overline{c_1}-a_2 \overline{c_2})-
  \overline{a_1} a_3 c_1 \overline{c_3}-a_1 \overline{a_3}
  \overline{c_1} c_3+a_2 \overline{a_2} c_1 \overline{c_3}+
  a_1 \overline{a_3} c_2 \overline{c_2}\\
=&a_1 \overline{c_3}+\overline{a_3} c_1-a_1 \overline{a_3}(c_1
  \overline{c_3}+\overline{c_1} c_3-c_2 \overline{c_2})-c_1
  \overline{c_3}(a_1 \overline{a_3}+\overline{a_1} a_3-a_2
  \overline{a_2})\\
=&a_1 \overline{c_3}+\overline{a_3} c_1.    
\endaligned$$ 
$$\aligned
 &(a_2 c_3-a_3 c_2)(\overline{a_1} \overline{c_3}-\overline{a_3}
  \overline{c_1})\\
=&a_2 \overline{c_3}(1-\overline{a_3} c_1+\overline{a_2} c_2)+
  \overline{a_3} c_2(1+a_2 \overline{c_2}-a_1 \overline{c_3})-
  \overline{a_1} a_3 c_2 \overline{c_3}-a_2 \overline{a_3}
  \overline{c_1} c_3\\
=&a_2 \overline{c_3}+\overline{a_3} c_2-a_2 \overline{a_3}(c_1
  \overline{c_3}+\overline{c_1} c_3-c_2 \overline{c_2})+c_2
  \overline{c_3}(a_2 \overline{a_2}-a_1 \overline{a_3}-
  \overline{a_1} a_3)\\
=&a_2 \overline{c_3}+\overline{a_3} c_2.    
\endaligned$$ 
$$\aligned
 &(a_2 c_3-a_3 c_2)(\overline{a_2} \overline{c_3}-\overline{a_3}
  \overline{c_2})\\
=&a_3 \overline{c_3}(1-\overline{a_3} c_1-\overline{a_1} c_3)+
  \overline{a_3} c_3(1-a_3 \overline{c_1}-a_1 \overline{c_3})+
  a_2 \overline{a_2} c_3 \overline{c_3}+a_3 \overline{a_3} c_2
  \overline{c_2}\\  
=&a_3 \overline{c_3}+\overline{a_3} c_3-a_3 \overline{a_3}(c_1
  \overline{c_3}+\overline{c_1} c_3-c_2 \overline{c_2})-c_3
  \overline{c_3} (\overline{a_1} a_3+a_1 \overline{a_3}-a_2
  \overline{a_2})\\
=&a_3 \overline{c_3}+\overline{a_3} c_3.  
\endaligned$$ 
By these identities, we get the desired formula.
 
  Similarly, we have
$$F_{1, 2}=\overline{v} |j(\gamma, Z)|^{-2}, \quad
  F_{2, 1}=v |j(\gamma, Z)|^{-2}, \quad 
  F_{2, 2}=|j(\gamma, Z)|^{-2},$$ 
and
$$H_1=-k |j(\gamma, Z)|^{-2}, \quad H_2=0, \quad
  G_1=k |j(\gamma, Z)|^{-2}, \quad G_2=0.$$ 
Therefore,
$$\aligned
 &L_k(f(\gamma(Z)) j(\gamma, Z)^{-k} j(\overline{\gamma}, 
  \overline{Z})^{k})\\
=&j(\gamma, Z)^{-k} j(\overline{\gamma}, \overline{Z})^{k}
  (z_1+\overline{z_1}-z_2 \overline{z_2}) |j(\gamma, Z)|^{-2}\\
 &\times \left[(u+\overline{u}) \frac{\partial^2}{\partial u
  \partial \overline{u}}+\frac{\partial^2}{\partial v \partial
  \overline{v}}+v \frac{\partial^2}{\partial \overline{u} \partial
  v}+\overline{v} \frac{\partial^2}{\partial u \partial \overline{v}}
  -k\left( \frac{\partial}{\partial u}-\frac{\partial}{\partial
  \overline{u}}\right)\right] f(\gamma(Z))\\
=&j(\gamma, Z)^{-k} j(\overline{\gamma}, \overline{Z})^{k}
  (L_k f)(\gamma(Z)).   
\endaligned$$
$\qquad \qquad \qquad \qquad \qquad \qquad \qquad \qquad \qquad 
 \qquad \qquad \qquad \qquad \qquad \qquad \qquad \qquad \qquad
 \quad \boxed{}$

  For the Picard modular group $\Gamma=U(2, 1; {\Cal O})$, the 
Eisenstein series of weight $k$ is defined by
$$E(Z; k, s)=\sum_{\gamma \in \Gamma \cap N \backslash \Gamma}
  j(\overline{\gamma}, \overline{Z})^{k} j(\gamma, Z)^{-k}
  \rho(\gamma, Z)^{s}, \quad \text{for} \quad \text{Re}(s)>2,
  k \in {\Bbb Z}. \tag 4.7$$

{\smc Theorem 4.3}. {\it The Eisenstein series satisfies the 
following equation:
$$L_{k} E(Z; k, s)=s(s-2) E(Z; k, s).\tag 4.8$$}

{\it Proof}. First, we note that
$$L_{k} \rho^{s}=L \rho^{s}+\rho \cdot k \left(\frac{\partial}
  {\partial z_1}-\frac{\partial}{\partial \overline{z_1}} \right)
  \rho^{s}=s(s-2) \rho^{s}.$$
By Proposition 4.2, we have
$$\aligned
  L_{k} E(Z; k, s)
&=\sum_{\gamma \in \Gamma \cap N \backslash \Gamma}
  L_{k}[\rho(\gamma(Z))^{s} j(\gamma, Z)^{-k} j(\overline{\gamma},
  \overline{Z})^{k}]\\
&=\sum_{\gamma \in \Gamma \cap N \backslash \Gamma}
  j(\gamma, Z)^{-k} j(\overline{\gamma}, \overline{Z})^{k}
  (L_{k} \rho^{s})(\gamma(Z))\\
&=s(s-2) \sum_{\gamma \in \Gamma \cap N \backslash \Gamma}
  j(\gamma, Z)^{-k} j(\overline{\gamma}, \overline{Z})^{k}
  \rho^{s}(\gamma(Z))\\
&=s(s-2) E(Z; k, s).        
\endaligned$$
$\qquad \qquad \qquad \qquad \qquad \qquad \qquad \qquad \qquad
 \qquad \qquad \qquad \qquad \qquad \qquad \qquad \qquad \qquad
 \quad \boxed{}$

  Set
$$J(\gamma, Z):=\frac{j(\gamma, Z)}{|j(\gamma, Z)|}. \tag 4.9$$
Then $J(\gamma, Z)=e^{i \arg j(\gamma, Z)}$.

  A generalized automorphic form on the Picard modular group
$\Gamma=U(2, 1; {\Cal O})$ of weight $3k$ is given by
$$f(\gamma(Z))=\frac{j(\gamma, Z)^{3k}}{|j(\gamma, Z)|^
  {3k}} f(Z) \quad \text{for} \quad \gamma \in \Gamma.\tag 4.10$$
Put
$$H_{k}(Z, W)=\frac{\rho(Z, W)^{2k}}{|\rho(Z, W)|^{2k}} \quad \text{for}
  \quad Z, W \in {\frak S}_2.\tag 4.11$$

{\it Definition} {\smc 4.4}. For any function 
$\Phi \in C^{\infty}({\Bbb R})$, the associated real 
point-pair invariant is given by
$$k(Z, W)=\Phi \left[ \frac{|\rho(Z, W)|^2}{\rho(Z) \rho(W)} \right],
  \tag 4.12$$
where $Z, W \in {\frak S}_{2}$.

{\smc Proposition 4.5}. {\it The following formulas hold:
\roster
\item $k(Z, W)=k(g(Z), g(W))$, for any $g \in G$.
\item $k(Z, W)=k(W, Z)$.
\endroster}

{\it Definition} {\smc 4.6}. Let $\Phi \in C^{\infty}({\Bbb R})$ be a 
real-valued function. The complex point-pair invariant is given by
$$K_{k}(Z, W)=H_{k}(Z, W) \Phi \left[ \frac{|\rho(Z, W)|^{2}}{\rho(Z) 
  \rho(W)} \right],\tag 4.13$$  
where $Z, W \in {\frak S}_{2}$.

{\smc Proposition 4.7}. {\it The following formulas hold:
\roster
\item $H_{k}(g(Z), g(W))=J(g, Z)^{2k} H_{k}(Z, W) J(g, W)^{-2k}$, 
      for $g \in G$.
\item $K_{k}(g(W), g(W))=J(g, Z)^{2k} K_{k}(Z, W) J(g, W)^{-2k}$, 
      for $g \in G$.
\item $|H_{k}(Z, W)|=1$, and $H_{k}(Z, W)=\overline{H_{k}(W, Z)}$.
\item $K_{k}(Z, W)=\overline{K_{k}(W, Z)}$.
\endroster}

{\it Proof}. Since $\rho(g(Z), g(W)) \overline{j(g, Z)} j(g, W)=
\rho(Z, W)$ for $g \in G$. We have
$$H_{k}(g(Z), g(W))=\frac{|j(g, Z)|^{2k}}{\overline{j(g, Z)}^{2k}}
  H_{k}(Z, W) \frac{|j(g, W)|^{2k}}{j(g, W)^{2k}}=J(g, Z)^{2k} 
  H_{k}(Z, W) J(g, W)^{-2k}.$$
The others are trivial.
\flushpar
$\qquad \qquad \qquad \qquad \qquad \qquad \qquad \qquad \qquad
 \qquad \qquad \qquad \qquad \qquad \qquad \qquad \qquad \qquad
 \quad \boxed{}$

{\smc Lemma 4.8}. {\it The following formula holds: 
$$L_{k} K_{k}(Z, W)=H_{k}(Z, W) \left[(\sigma^2-\sigma) \Phi^{\prime 
  \prime} (\sigma)+(3 \sigma-1) \Phi^{\prime}(\sigma)+\frac{k^2}{\sigma}
  \Phi(\sigma) \right],\tag 4.14$$
where $\sigma(Z, W)=\frac{|\rho(Z, W)|^2}{\rho(Z) \rho(W)}$.}  

{\it Proof}. It is obtained by a straightforward calculation.
\flushpar
$\qquad \qquad \qquad \qquad \qquad \qquad \qquad \qquad \qquad 
 \qquad \qquad \qquad \qquad \qquad \qquad \qquad \qquad \qquad
 \quad \boxed{}$

  By $L_{k} K_{k}(Z, W)=s(s-2) K_{k}(Z, W)$, one has
$$\Phi^{\prime \prime}(\sigma)+\left(\frac{1}{\sigma}+\frac{2}{\sigma-1}
  \right) \Phi^{\prime}(\sigma)+\frac{1}{\sigma(\sigma-1)} \left[ 
  \frac{k^2}{\sigma}-s(s-2) \right] \Phi(\sigma)=0.\tag 4.15$$
It is known that the Fuchsian equation with three regular singularities
$a, b, c$ is given by
$$\aligned
 &\frac{d^2 w}{dz^2}+\left[ \frac{1-\alpha_1-\alpha_2}{z-a}+
  \frac{1-\beta_1-\beta_2}{z-b}+\frac{1-\gamma_1-\gamma_2}{z-c}
  \right] \frac{dw}{dz}+\frac{1}{(z-a)(z-b)(z-c)}\\
 &\times \left[ \frac{\alpha_1 \alpha_2 (a-b)(a-c)}
  {z-a}+\frac{\beta_1 \beta_2 (b-c)(b-a)}{z-b}+\frac{\gamma_1 \gamma_2 
  (c-a)(c-b)}{z-c} \right] w=0,
\endaligned$$
where $(\alpha_1, \alpha_2)$, $(\beta_1, \beta_2)$, $(\gamma_1, 
\gamma_2)$ are exponents belonging to $a, b, c$, respectively, 
and, they satisfy that 
$\alpha_1+\alpha_2+\beta_1+\beta_2+\gamma_1+\gamma_2=1$.
The solution of this equation is given in Riemann P-notation by
$$w(z)=P \left\{\matrix
         a & b & c\\
		 \alpha_1 & \beta_1 & \gamma_1\\
		 \alpha_2 & \beta_2 & \gamma_2
		 \endmatrix; z \right\}.$$
When $c=\infty$, it reduces to
$$\aligned
 &\frac{d^2 w}{d z^2}+\left[\frac{1-\alpha_1-\alpha_2}{z-a}+
  \frac{1-\beta_1-\beta_2}{z-b} \right] \frac{dw}{dz}+\frac{1}
  {(z-a)(z-b)}\\
 &\times \left[\frac{\alpha_1 \alpha_2 (a-b)}{z-a}+\frac{\beta_1 
  \beta_2 (b-a)}{z-b}+\gamma_1 \gamma_2 \right] w=0.
\endaligned$$
Now, in our case,
$$a=0, b=1, c=\infty, \alpha_1=-|k|, \alpha_2=|k|, \beta_1=0,  
  \beta_2=-1, \gamma_1=s, \gamma_2=2-s,$$
and
$$\aligned
 &\Phi(\sigma)=\Phi_{|k|}(\sigma)\\
=&P \left\{\matrix
         0 &  1 & \infty\\  
      -|k| &  0 & s\\
	   |k| & -1 & 2-s                        
        \endmatrix; \sigma \right\}\\
=&\sigma^{-|k|} P \left\{\matrix
         0 &  1 & \infty\\
         0 &  0 & s-|k|\\
      2|k| & -1 & 2-s-|k|
       \endmatrix; \sigma \right\}\\
=&\sigma^{-|k|} F(s-|k|, 2-s-|k|; 1-2|k|; \sigma)\\	   		
=&\sigma^{-|k|} (1-\sigma)^{|k|-s} F(s-|k|, s-1-|k|; 1-2|k|; 
  \frac{\sigma}{\sigma-1})\\
=&\sigma^{-|k|} (1-\sigma)^{|k|-s} F(s-|k|, s-1-|k|; 2s-1; 
  -\frac{1}{\sigma-1}).
\endaligned$$
Now, we set 
$$\aligned
  K(Z, W; k, s)
=&\frac{\pi^{\frac{3}{2}} \Gamma(\frac{3}{2}-s)}
  {(|k|+1-s) \Gamma(2-s)} (\frac{1}{4})^{s-1} \rho(Z, W)^{k} 
  \rho(\overline{Z}, \overline{W})^{-k}\\
 &\times \sigma^{-|k|} (\sigma-1)^{|k|-s} 
  F(s-|k|, s-1-|k|; 2s-1; -\frac{1}{\sigma-1}).
\endaligned\tag 4.16$$
It has the following property:
$$K(\gamma(Z), \gamma(Z^{\prime}); k, s)=j(\gamma, Z)^{k}
  j(\overline{\gamma}, \overline{Z})^{-k} j(\gamma, Z^{\prime})^{-k}
  j(\overline{\gamma}, \overline{Z^{\prime}})^{k} 
  K(Z, Z^{\prime}; k, s).\tag 4.17$$
Moreover, note that  
$\rho(Z^{\prime}, Z)=\rho(\overline{Z}, \overline{Z^{\prime}})$,  
we have
$$K(Z^{\prime}, Z; k, s)=K(Z, Z^{\prime}; -k, s).$$  
  
  The automorphic Green function of weight $k$ is given by
$$G(Z, Z^{\prime}; k, s)=\sum_{\gamma \in \Gamma}
  K(Z, \gamma(Z^{\prime}); k, s), \quad \text{for} \quad
  \text{Re}(s)>\delta(\Gamma). \tag 4.18$$
By the same method, we can prove that the right hand side of (4.18)
is convergent if $\text{Re}(s)>\delta(\Gamma)$.
 
  The automorphic Green function of weight $k$ has the following 
properties:
$$G(Z, Z^{\prime}; k, s)=G(Z^{\prime}, Z; -k, s), \quad
  L_{k}(Z) G(Z, Z^{\prime}; k, s)=s(s-2) G(Z, Z^{\prime}; k, s),$$
Furthermore,
$$\aligned
 &L_{-k}(Z^{\prime}) G(Z, Z^{\prime}; k, s)=
  L_{-k}(Z^{\prime}) G(Z^{\prime}, Z; -k, s)\\
=&s(s-2) G(Z^{\prime}, Z; -k, s)=
  s(s-2) G(Z, Z^{\prime}; k, s).
\endaligned$$

  The Poisson kernel of weight $k$ is given by
$$P_{k}(Z, W; s)=\rho(Z)^{s} \rho(Z, W)^{k-s} 
  \rho(\overline{Z}, \overline{W})^{-k-s},\tag 4.19$$
where $Z \in {\frak S}_2$, $W \in \partial {\frak S}_2$.
By 
$$\rho(\gamma(Z), \gamma(W)) \overline{j(\gamma, Z)} j(\gamma, W)
=\rho(Z, W), \quad j(\gamma, \gamma^{-1}(W))=j(\gamma^{-1}, W)^{-1},$$
we have
$$P_{k}(\gamma(Z), W; s)=P_{k}(Z, \gamma^{-1}(W); s) j(\gamma, Z)^{k}
  j(\overline{\gamma}, \overline{Z})^{-k} j(\gamma^{-1}, W)^{k-s}
  j(\overline{\gamma}^{-1}, \overline{W})^{-k-s}.$$

  The Eisenstein series with respect to $\partial {\frak S}_{2}$ 
of weight $k$ is defined by
$$E(Z, W; k, s)=\sum_{\gamma \in \Gamma} j(\gamma, Z)^{-k}
  j(\overline{\gamma}, \overline{Z})^{k} P_{k}(\gamma(Z), W; s),
  \quad \text{for} \quad \text{Re}(s)>\delta(\Gamma).
  \tag 4.20$$
Here $\Gamma$ is a discrete subgroup of $U(2, 1)$.
Now, we have the following formula:
$$\lim_{\rho(Z) \to 0} \rho(Z)^{-s} K_{k}(Z, Z^{\prime})
 =\frac{\pi^{\frac{3}{2}} \Gamma(\frac{3}{2}-s)}{(|k|+1-s)
  \Gamma(2-s)} (\frac{1}{4})^{s-1} P_{k}(Z^{\prime}, W; s),$$
where $Z \to W \in \partial {\frak S}_2$ as $\rho(Z) \to 0$.
Consequently,
$$\lim_{\rho(Z) \to 0} \rho(Z)^{-s} G(Z, Z^{\prime}; k, s)
 =\frac{\pi^{\frac{3}{2}} \Gamma(\frac{3}{2}-s)}{(|k|+1-s)
  \Gamma(2-s)} (\frac{1}{4})^{s-1} E(Z^{\prime}, W; k, s).$$

{\smc Lemma 4.9}. {\it The Eisenstein series of weight $k$ satisfies 
the following three equations:
\roster
\item $L_{k} E(Z, W; k, s)=s(s-2) E(Z, W; k, s)$.
\item $E(\gamma(Z), W; k, s)=E(Z, W; k, s)$, for $\gamma \in \Gamma$.
\item $E(Z, g(W); k, s)=j(g, W)^{s-k} j(\overline{g},\overline{W})^{s+k}
       E(Z, W; k, s)$, where $g \in \Gamma$.
\endroster}		

{\it Proof}. In fact,
$$L_{k} P_{k}(Z, W; s)=s(s-2) P_{k}(Z, W; s)+(k^2-s^2)
  \rho(W) P_{k}(Z, W; s+1).$$
Since $\rho(W)=0$, $L_{k} P_{k}(Z, W; s)=s(s-2) P_{k}(Z, W; s)$.
By Proposition 4.2,
$$\aligned
  L_{k} E(Z, W; k, s)
=&\sum_{\gamma \in \Gamma} L_{k}[P_{k}(\gamma(Z), W; s)
  j(\gamma, Z)^{-k} j(\overline{\gamma}, \overline{Z})^{k}]\\
=&\sum_{\gamma \in \Gamma} j(\gamma, Z)^{-k} j(\overline{\gamma},
  \overline{Z})^{k} (L_{k} P_{k})(\gamma(Z), W; s)\\   
=&s(s-2) \sum_{\gamma \in \Gamma} j(\gamma, Z)^{-k} j(\overline{\gamma},
  \overline{Z})^{k} P_{k}(\gamma(Z), W; s)\\
=&s(s-2) E(Z, W; k, s).
\endaligned$$

  To prove the third identity, we note that
$$P_{k}(\gamma(Z), \gamma(W); s)=j(\gamma, Z)^{k} j(\overline{\gamma},
  \overline{Z})^{-k} j(\gamma, W)^{s-k} j(\overline{\gamma},
  \overline{W})^{s+k} P_{k}(Z, W; s).$$   
By Proposition 2.9, we have
$$\aligned
  E(Z, g(W); k, s)
=&\sum_{\gamma \in \Gamma} j(\gamma, Z)^{-k} j(\overline{\gamma},
  \overline{Z})^{k} P_{k}(\gamma(Z), g(W); s)\\
=&\sum_{\gamma \in \Gamma} j(g, \gamma(Z))^{-k} j(\overline{g},
  \overline{\gamma}(\overline{Z}))^{k} j(\gamma, Z)^{-k} 
  j(\overline{\gamma}, \overline{Z})^{k} 
  P_{k}(g(\gamma(Z)), g(W); s)\\
=&\sum_{\gamma \in \Gamma} j(\gamma, Z)^{-k} j(\overline{\gamma},
  \overline{Z})^{k} P_{k}(\gamma(Z), W; s) j(g, W)^{s-k}
  j(\overline{g}, \overline{W})^{s+k}\\
=&j(g, W)^{s-k} j(\overline{g}, \overline{W})^{s+k}
  E(Z, W; k, s).
\endaligned$$
$\qquad \qquad \qquad \qquad \qquad \qquad \qquad \qquad \qquad
 \qquad \qquad \qquad \qquad \qquad \qquad \qquad \qquad \qquad
 \quad \boxed{}$

{\smc Theorem 4.10}. {\it Set
$$\phi(Z, Z^{\prime}; k, s)=\int_{\partial {\frak S}_2}
  P_{k}(Z, W; s) P_{-k}(Z^{\prime}, W; 2-s) dm(W).\tag 4.21$$
Then the following formula holds:
$$\aligned
 &\phi(Z, Z^{\prime}; k, s)
=K(Z, Z^{\prime}; k, s)+K(Z, Z^{\prime}; k, 2-s)\\
=&\phi(Z^{\prime}, Z; -k, 2-s)
=K(Z^{\prime}, Z; -k, s)+K(Z^{\prime}, Z; -k, 2-s),
\endaligned\tag 4.22$$  
where $k=0$ or $k=\pm 1$.}
  
{\it Proof}. In fact, if $k=0$, this is Theorem 3.7.  
Now, we only consider the nontrivial case $k=\pm 1$.

  For $\gamma \in G$,
$$\aligned
 &\phi(\gamma(Z), \gamma(Z^{\prime}); k, s)
 =\int_{\partial {\frak S}_{2}} P_{k}(\gamma(Z), W; s)
  P_{-k}(\gamma(Z^{\prime}), W; 2-s) dm(W)\\
=&j(\gamma, Z)^{k} j(\overline{\gamma}, \overline{Z})^{-k}
  j(\gamma, Z^{\prime})^{-k} j(\overline{\gamma}, 
  \overline{Z^{\prime}})^{k} \int_{\partial {\frak S}_{2}} 
  P_{k}(Z, \gamma^{-1}(W); s)\\
 &\times P_{-k}(Z^{\prime}, \gamma^{-1}(W); 2-s) j(\gamma^{-1}, W)^{-2}
  j(\overline{\gamma}^{-1}, \overline{W})^{-2} dm(W).
\endaligned$$
Since $dm(\gamma^{-1}(W))=|j(\gamma^{-1}, W)|^{-4} dm(W)$, we have
$$\phi(\gamma(Z), \gamma(Z^{\prime}); k, s)
 =j(\gamma, Z)^{k} j(\overline{\gamma}, \overline{Z})^{-k}
  j(\gamma, Z^{\prime})^{-k} j(\overline{\gamma},
  \overline{Z^{\prime}})^{k} \phi(Z, Z^{\prime}; k, s).
  \tag 4.23$$
Thus, $\phi(Z, Z^{\prime}; k, s)$ is a point-pair invariant covariant
with respect to the weight $k$.

  By the double transitivity of $U(2, 1)$ on ${\frak S}_{2}$, 
there exists a $\gamma \in U(2, 1)$, such that
$$Z=\gamma(\frac{1}{2} \rho, 0), \quad
  Z^{\prime}=\gamma(\frac{1}{2} \rho^{\prime}, 0).$$ 
Let $\rho$ and $\rho^{\prime}$ stand for $(\rho, 0)$ and 
$(\rho^{\prime}, 0)$, respectively.
Set $\gamma=\left(\matrix
                  * & *   & *\\
				  * & *   & *\\
				a_1 & a_2 & a_3
		        \endmatrix\right)$, then
$$|j(\gamma, \frac{1}{2} \rho)|^{2}=|\frac{1}{2} \rho a_1+a_3|^{2}
 =\frac{1}{4} \rho^{2} |a_1|^2+\frac{1}{2} \rho |a_2|^2+|a_3|^2.$$
Denote
$$\mu=\frac{\frac{1}{2}(\rho+\rho^{\prime})}
  {\frac{1}{4} \rho^2 |a_1|^2+\frac{1}{2} \rho |a_2|^2+|a_3|^2},
  \tag 4.24$$
then $\mu>0$. 

  Since
$$\aligned
 &\phi(Z, Z^{\prime}; k, s)=\phi(\gamma(\frac{1}{2} \rho),
  \gamma(\frac{1}{2} \rho^{\prime}); k, s)\\
=&j(\gamma, \frac{1}{2} \rho)^{k} j(\overline{\gamma}, 
  \frac{1}{2} \rho)^{-k} j(\gamma, \frac{1}{2} \rho^{\prime})^{-k} 
  j(\overline{\gamma}, \frac{1}{2} \rho^{\prime})^{k} 
  \phi(\frac{1}{2} \rho, \frac{1}{2} \rho^{\prime}; k, s).
\endaligned$$
By (4.24), we have
$$\rho(Z, Z^{\prime})=\rho(\gamma(\frac{1}{2} \rho), 
  \gamma(\frac{1}{2} \rho^{\prime}))
 =\frac{\frac{1}{2}(\rho+\rho^{\prime})}{j(\overline{\gamma}, 
  \frac{1}{2} \rho) j(\gamma, \frac{1}{2} \rho^{\prime})}
 =\mu \frac{j(\gamma, \frac{1}{2} \rho)}{j(\gamma, \frac{1}{2} 
  \rho^{\prime})}.$$
Hence,
$$\aligned
  \phi(Z, Z^{\prime}; k, s)
=&[\mu^{-1} \rho(Z, Z^{\prime})]^{k} [\mu^{-1} \rho(\overline{Z},
  \overline{Z^{\prime}})]^{-k} \phi(\frac{1}{2} \rho, \frac{1}{2}
  \rho^{\prime}; k, s)\\
=&\rho(Z, Z^{\prime})^{k} \rho(\overline{Z}, \overline{Z^{\prime}})^
  {-k} \phi(\frac{1}{2} \rho, \frac{1}{2} \rho^{\prime}; k, s).
\endaligned\tag 4.25$$

  For $W \in \partial {\frak S}_{2}$, i.e., $w_1=\frac{|w|^2}{2}+i v$, 
$w_2=w$ with $w \in {\Bbb C}$, $v \in {\Bbb R}$. 
$$\rho(Z, W)=\frac{\rho+|w|^2}{2}+iv, \quad
  \rho(Z^{\prime}, W)=\frac{\rho^{\prime}+|w|^2}{2}+iv.$$
Therefore,
$$\aligned
  \phi(\frac{1}{2} \rho, \frac{1}{2} \rho^{\prime}; k, s)
=&16 \rho^{s} {\rho^{\prime}}^{2-s} \int_{{\Bbb R}} \int_{{\Bbb C}}
  \frac{1}{|\rho+|w|^2+2iv|^{2s}} \frac{1}{|\rho^{\prime}+|w|^2
  +2iv|^{4-2s}}\\
 &\times \left(\frac{\rho+|w|^2+2iv}{\rho+|w|^2-2iv}\right)^{k}
  \left(\frac{\rho^{\prime}+|w|^2+2iv}{\rho^{\prime}+|w|^2-2iv}
  \right)^{-k} dw d \overline{w} dv. 
\endaligned$$
Set $w=r e^{i \theta}$ and $u=r^2$, then
$$\aligned
  \phi(\frac{1}{2} \rho, \frac{1}{2} \rho^{\prime}; k, s)
=&16 \pi \rho^{s} {\rho^{\prime}}^{2-s} \int_{-\infty}^{\infty}
  \int_{0}^{\infty} \frac{1}{|\rho+u+2iv|^{2s}} 
  \frac{1}{|\rho^{\prime}+u+2iv|^{4-2s}}\\ 
 &\times \left(\frac{\rho+u+2iv}{\rho+u-2iv}\right)^{k}
  \left(\frac{\rho^{\prime}+u+2iv}{\rho^{\prime}+u-2iv}
  \right)^{-k} du dv\\
=&8 \pi \rho^{s} {\rho^{\prime}}^{2-s} \int_{-\infty}^{\infty}
  \int_{0}^{\infty} \frac{1}{|\rho+u+iv|^{2s}} \frac{1}
  {|\rho^{\prime}+u+iv|^{2(2-s)}}\\
 &\times \exp(2ki(\arctan \frac{v}{\rho+u}-\arctan \frac{v}
  {\rho^{\prime}+u})) du dv.   
\endaligned$$
By
$$\arctan \frac{v}{\rho+u}-\arctan \frac{v}{\rho^{\prime}+u}=
\arctan \frac{(\rho^{\prime}-\rho) v}{(\rho+u)(\rho^{\prime}+u)+v^2},$$
we have
$$\aligned
  \phi(\frac{1}{2} \rho, \frac{1}{2} \rho^{\prime}; k, s)
=&8 \pi \rho^{s} {\rho^{\prime}}^{2-s} \int_{-\infty}^{\infty}
  \int_{0}^{\infty} \frac{1}{[(\rho+u)^2+v^2]^{s}} \frac{1}
  {[(\rho^{\prime}+u)^2+v^2]^{2-s}}\\
 &\times \exp(2ki \arctan \frac{(\rho^{\prime}-\rho) v}{(\rho+u)
  (\rho^{\prime}+u)+v^2}) du dv\\
=&16 \pi \rho^{s} {\rho^{\prime}}^{2-s} \int_{0}^{\infty} 
  \int_{0}^{\infty} \frac{1}{[(\rho+u)^2+v^2]^{s}} \frac{1}
  {[(\rho^{\prime}+u)^2+v^2]^{2-s}}\\
 &\times \cos 2k \arctan \frac{(\rho^{\prime}-\rho) v}{(\rho+u)
  (\rho^{\prime}+u)+v^2} du dv.
\endaligned$$
Assume that $\rho<\rho^{\prime}$, set 
$\lambda=\frac{\rho}{\rho^{\prime}}$, $v=(\rho^{\prime}-\rho) w$ 
and $l=\frac{\rho+u}{\rho^{\prime}+u}$, then 
$u=\frac{l \rho^{\prime}-\rho}{1-l}$, 
$du=\frac{\rho^{\prime}-\rho}{(1-l)^2} dl$.
We have
$$\aligned
  \phi(\frac{1}{2} \rho, \frac{1}{2} \rho^{\prime}; k, s)
=&16 \pi \lambda^{s} (1-\lambda)^{-2} \int_{0}^{\infty} 
  \int_{\lambda}^{1} \left[w^2+\frac{l^2}{(1-l)^2}\right]^{-s}
  \left[w^2+\frac{1}{(1-l)^2} \right]^{s-2}\\
 &\times \cos 2k \arctan \left(\frac{w}{w^2+\frac{l}{(1-l)^2}}\right)
  (1-l)^{-2} dl dw.
\endaligned$$
Set $w=\frac{\sqrt{v}}{1-l}$, then
$$\aligned
  \phi(\frac{1}{2} \rho, \frac{1}{2} \rho^{\prime}; k, s)
=&8 \pi \lambda^{s} (1-\lambda)^{-2} \int_{\lambda}^{1} (1-l) dl
  \int_{0}^{\infty} v^{-\frac{1}{2}} (1+v)^{s-2} (l^2+v)^{-s}\\
 &\times \cos 2k \arctan \left[\frac{(1-l) \sqrt{v}}{v+l} \right] dv.
\endaligned$$
Let $x=\arctan \frac{(1-l) \sqrt{v}}{v+l}$, then
$$\cos 2x=\frac{1-\tan^2 x}{1+\tan^2 x}=\frac{(v+l)^2-(1-l)^2 v}
  {(v+l)^2+(1-l)^2 v}.$$
We have
$$\cos 2kx=\cos \left(k \arccos \frac{(v+l)^2-(1-l)^2 v}
  {(v+l)^2+(1-l)^2 v}\right).$$
It is known that the Tchebichef polynomial
$$T_n(x)=\cos(n \arccos x)=F(-n, n; \frac{1}{2}; \frac{1-x}{2}).$$
Hence
$$\cos 2kx=F(-k, k; \frac{1}{2}; \frac{(1-l)^2 v}{v^2+(1+l^2)v+l^2}).$$
Therefore,
$$\aligned
  \phi(\frac{1}{2} \rho, \frac{1}{2} \rho^{\prime}; k, s)
=&8 \pi \lambda^s (1-\lambda)^{-2} \int_{\lambda}^{1} (1-l) dl
  \int_{0}^{\infty} v^{-\frac{1}{2}} (1+v)^{s-2} (l^2+v)^{-s}\\
 &\times F(-k, k; \frac{1}{2}; \frac{(1-l)^2 v}{(1+v)(l^2+v)}) dv.
\endaligned$$
We know that
$$F(-k, k; \frac{1}{2}; x)=\sum_{n=0}^{|k|} \frac{(-k)_{n} (k)_{n}}
  {n! (\frac{1}{2})_{n}} x^n.$$
Thus
$$\aligned
  \phi(\frac{1}{2} \rho, \frac{1}{2} \rho^{\prime}; k, s)
=&8 \pi \lambda^{s} (1-\lambda)^{-2} \sum_{n=0}^{|k|}
  \frac{(-k)_{n} (k)_{n}}{n! (\frac{1}{2})_{n}}
  \int_{\lambda}^{1} (1-l)^{2n+1} l^{-2(s+n)} dl\\
 &\times \int_{0}^{\infty} v^{n-\frac{1}{2}}
  (1+v)^{s-n-2} (1+v l^{-2})^{-s-n} dv. 
\endaligned$$
By the following formula
$$F(a, b; c; 1-z)=\frac{\Gamma(c)}{\Gamma(b) \Gamma(c-b)}
  \int_{0}^{\infty} s^{b-1} (1+s)^{a-c} (1+sz)^{-a} ds$$
for $\text{Re}(c)>\text{Re}(b)>0$, $|\arg(z)|<\pi$, we have
$$\aligned
 &\int_{0}^{\infty} v^{n-\frac{1}{2}} (1+v)^{s-n-2}
  (1+v l^{-2})^{-s-n} dv\\
=&\frac{\Gamma(n+\frac{1}{2}) \Gamma(n+\frac{3}{2})}
  {\Gamma(2n+2)} F(s+n, \frac{1}{2}+n; 2n+2; 1-l^{-2}).
\endaligned$$
Hence,
$$\aligned
  \phi(\frac{1}{2} \rho, \frac{1}{2} \rho^{\prime}; k, s)
=&8 \pi \lambda^{s} (1-\lambda)^{-2} \sum_{n=0}^{|k|}
  \frac{(-k)_{n} (k)_{n}}{n! (\frac{1}{2})_{n}}
  \frac{\Gamma(n+\frac{1}{2}) \Gamma(n+\frac{3}{2})}
  {\Gamma(2n+2)}\\  
 &\times \int_{\lambda}^{1} (1-l)^{2n+1} l^{-2(s+n)}
  F(s+n, \frac{1}{2}+n; 2n+2; 1-l^{-2}) dl.  
\endaligned$$
Note that
$$\frac{1}{n! (\frac{1}{2})_{n}} \frac{\Gamma(n+\frac{1}{2})
  \Gamma(n+\frac{3}{2})}{\Gamma(2n+2)}=\frac{\pi}{2^{2n+1} 
  (n!)^{2}}.$$
We have
$$\aligned
  \phi(\frac{1}{2} \rho, \frac{1}{2} \rho^{\prime}; k, s)
=&4 \pi^{2} \lambda^{s} (1-\lambda)^{-2} \sum_{n=0}^{|k|}
  \frac{(-k)_{n} (k)_{n}}{2^{2n} (n!)^{2}} \int_{\lambda}^{1}
  (1-l)^{2n+1} l^{-2(s+n)}\\
 &\times F(s+n, \frac{1}{2}+n; 2n+2; 1-l^{-2}) dl.   
\endaligned$$

(1) $z=1-l^{-2}$, $\frac{z}{z-1}=1-l^2$.

  By $F(a, b; c; z)=(1-z)^{-a} F(a, c-b; c; \frac{z}{z-1})$, we have
$$F(s+n, \frac{1}{2}+n; 2n+2; 1-l^{-2})=l^{2(s+n)}
  F(s+n, n+\frac{3}{2}; 2n+2; 1-l^{2}).$$   

(2) $z=1-l^2$, $1-z=l^2$.

  By
$$\aligned
 &F(a, b; c; z)=\frac{\Gamma(c) \Gamma(c-a-b)}{\Gamma(c-a) \Gamma(c-b)}
  F(a, b; a+b-c+1; 1-z)\\
 &+\frac{\Gamma(c) \Gamma(a+b-c)}{\Gamma(a) \Gamma(b)} (1-z)^{c-a-b}
  F(c-a, c-b; c-a-b+1; 1-z),
\endaligned$$  
where $|\arg(1-z)|<\pi$, we have
$$\aligned
 &F(s+n, n+\frac{3}{2}; 2n+2; 1-l^{2})\\
=&\frac{\Gamma(2n+2) \Gamma(\frac{1}{2}-s)}{\Gamma(n-s+2) \Gamma(n+
  \frac{1}{2})} F(s+n, n+\frac{3}{2}; s+\frac{1}{2}; l^2)\\
+&\frac{\Gamma(2n+2) \Gamma(s-\frac{1}{2})}{\Gamma(s+n) \Gamma(n+
  \frac{3}{2})} l^{1-2s} F(n-s+2, n+\frac{1}{2}; \frac{3}{2}-s; l^2).
\endaligned$$
Therefore,
$$\phi(\frac{1}{2} \rho, \frac{1}{2} \rho^{\prime}; k, s)
 =4 \pi^2 \lambda^{s} (1-\lambda)^{-2} \sum_{n=0}^{|k|}
  \frac{(-k)_{n} (k)_{n}}{2^{2n} (n!)^2} (F_1+F_2),$$
where
$$F_1=\frac{(2n+1)! \Gamma(\frac{1}{2}-s)}{\Gamma(n-s+2) 
  \Gamma(n+\frac{1}{2})} \int_{\lambda}^{1} (1-l)^{2n+1} 
  F(s+n, n+\frac{3}{2}; s+\frac{1}{2}; l^2) dl,$$
$$F_2=\frac{(2n+1)! \Gamma(s-\frac{1}{2})}{\Gamma(s+n) 
  \Gamma(n+\frac{3}{2})} \int_{\lambda}^{1} (1-l)^{2n+1} 
  l^{1-2s} F(n-s+2, n+\frac{1}{2}; \frac{3}{2}-s; l^2) dl.$$

  It is known that the incomplete beta function is defined by 
(see \cite{Er}, p.87)
$$B_{x}(p, q)=\int_{0}^{x} t^{p-1} (1-t)^{q-1} dt.$$
Set 
$I_{x}(p, q)=\frac{B_{x}(p, q)}{B_{1}(p, q)}$.
Note that
$B_{1}(p, q)=\frac{\Gamma(p) \Gamma(q)}{\Gamma(p+q)}$.

  Hence,
$$\aligned
 &\int_{\lambda}^{1} (1-l)^{2n+1} F(s+n, n+\frac{3}{2}; s+
  \frac{1}{2}; l^2) dl\\
=&\sum_{m=0}^{\infty} \frac{(s+n)_{m} (n+\frac{3}{2})_{m}}{m!
  (s+\frac{1}{2})_{m}} [B_{1}(2m+1, 2n+2)-B_{\lambda}(2m+1, 2n+2)].  
\endaligned$$
Similarly,
$$\aligned
 &\int_{\lambda}^{1} (1-l)^{2n+1} l^{1-2s} F(n-s+2, n+\frac{1}{2};
  \frac{3}{2}-s; l^2) dl\\
=&\sum_{m=0}^{\infty} \frac{(2-s+n)_{m} (n+\frac{1}{2})_{m}}{m!
  (\frac{3}{2}-s)_{m}} [B_{1}(2m-2s+2, 2n+2)-B_{\lambda}(2m-2s+2,
  2n+2)].
\endaligned$$
We have
$$\aligned
 &\sum_{m=0}^{\infty} \frac{(s+n)_{m} (n+\frac{3}{2})_{m}}
  {m! (s+\frac{1}{2})_{m}} B_{1}(2m+1, 2n+2)
=\sum_{m=0}^{\infty} \frac{(s+n)_{m}}{(s+\frac{1}{2})_{m}}
  \frac{(2m-1)!! (2n)!!}{(2m+2n+2)!!}\\
=&\frac{1}{2n+2} \sum_{m=0}^{\infty} \frac{(s+n)_{m} 
  (\frac{1}{2})_{m}}{(s+\frac{1}{2})_{m} (n+2)_{m}}.
\endaligned\tag 4.26$$
By 
$$(n+2)_{m}=\frac{(n+m+1)!}{(n+1)!}, \quad
  (\frac{1}{2})_{m}=\frac{(-\frac{1}{2}-n)_{m+n+1}}
  {(-\frac{1}{2}-n)_{n+1}},$$
$$(s+\frac{1}{2})_{m}=\frac{(s-n-\frac{1}{2})_{m+n+1}}
  {(s-n-\frac{1}{2})_{n+1}}, \quad
  (s+n)_{m}=\frac{(s-1)_{m+n+1}}{(s-1)_{n+1}},$$
(4.26) is equal to
$$\aligned
 &\frac{1}{2n+2} \frac{(n+1)! (s-n-\frac{1}{2})_{n+1}}
  {(s-1)_{n+1} (-\frac{1}{2}-n)_{n+1}} \sum_{m=0}^{\infty}
  \frac{(s-1)_{m+n+1} (-\frac{1}{2}-n)_{m+n+1}}{(s-n-\frac{1}
  {2})_{m+n+1} (m+n+1)!}\\
=&\frac{1}{2n+2} \frac{(n+1)! (s-n-\frac{1}{2})_{n+1}}
  {(s-1)_{n+1} (-\frac{1}{2}-n)_{n+1}}[F(s-1, -\frac{1}{2}-n;
  s-n-\frac{1}{2}; 1)\\
 &-\sum_{r=0}^{n} \frac{(s-1)_{r} (-\frac{1}{2}-n)_{r}}
  {(s-n-\frac{1}{2})_{r} r!}].
\endaligned\tag 4.27$$
Similarly,
$$\aligned
 &\sum_{m=0}^{\infty} \frac{(2-s+n)_{m} (n+\frac{1}{2})_{m}}
  {m! (\frac{3}{2}-s)_{m}} B_{1}(2m-2s+2, 2n+2)\\
=&\frac{(2n+1)! \Gamma(2-2s)}{\Gamma(4-2s+2n)} \sum_{m=0}^{\infty}
  \frac{(1-s)_{m} (n+\frac{1}{2})_{m}}{(\frac{5}{2}-s+n)_{m} m!}\\  
=&\frac{(2n+1)! \Gamma(2-2s)}{\Gamma(4-2s+2n)}
  F(1-s, n+\frac{1}{2}; \frac{5}{2}-s+n; 1).
\endaligned\tag 4.28$$
Therefore,
$$\aligned
 &\sum_{m=0}^{\infty} \frac{(s+n)_{m} (n+\frac{3}{2})_{m}}
  {m! (s+\frac{1}{2})_{m}} B_{\lambda}(2m+1, 2n+2)\\
=&\sum_{m=0}^{\infty} \frac{(s+n)_{m} (n+\frac{3}{2})_{m}}
  {m! (s+\frac{1}{2})_{m}} B_{1}(2m+1, 2n+2) I_{\lambda}(2m+1, 2n+2)\\
=&\frac{1}{2n+2} \sum_{m=0}^{\infty} \frac{(s+n)_{m} (\frac{1}{2})_{m}}
  {(s+\frac{1}{2})_{m} (n+2)_{m}} I_{\lambda}(2m+1, 2n+2).
\endaligned$$
Now, we need the following lemma:

{\smc Lemma 4.11}. {\it The explicit expression of $I_{x}(p, q)$: For
$p \in {\Bbb C}$, $\text{Re}(p)>0$ and $q \in {\Bbb Z}$ with $q>0$, 
the following formula holds:
$$I_{x}(p, q+1)=(p)_{q+1} \sum_{r=0}^{q} \frac{(-1)^{r}}{r!
  (q-r)!} \frac{x^{p+r}}{p+r}.\tag 4.29$$}

{\it Proof}.
$$\aligned
 &I_{x}(p, q+1)
 =\frac{(p)_{q+1}}{q!} B(p, q+1) I_{x}(p, q+1)\\
=&\frac{(p)_{q+1}}{q!} B_{x}(p, q+1)
 =\frac{(p)_{q+1}}{q!} \int_{0}^{x} t^{p-1} (1-t)^{q} dt\\
=&\frac{(p)_{q+1}}{q!} \sum_{r=0}^{q} (-1)^{r} C_{q}^{r}
  \frac{x^{p+r}}{p+r}
 =(p)_{q+1} \sum_{r=0}^{q} \frac{(-1)^{r}}{r! (q-r)!}
  \frac{x^{p+r}}{p+r}.  
\endaligned$$ 
$\qquad \qquad \qquad \qquad \qquad \qquad \qquad \qquad \qquad
 \qquad \qquad \qquad \qquad \qquad \qquad \qquad \qquad \qquad
 \quad \boxed{}$

  By Lemma 4.11, we have
$$\aligned
 &\frac{1}{2n+2} \sum_{m=0}^{\infty} \frac{(s+n)_{m} (\frac{1}{2})_{m}} 
  {(s+\frac{1}{2})_{m} (n+2)_{m}} I_{\lambda}(2m+1, 2n+2)\\
=&\frac{1}{2n+2} \sum_{m=0}^{\infty} \frac{(s+n)_{m} (\frac{1}{2})_{m}}
  {(s+\frac{1}{2})_{m} (n+2)_{m}} (2m+1)_{2n+2} \sum_{r=0}^{2n+1}
  \frac{(-1)^{r}}{r! (2n+1-r)!} \frac{\lambda^{2m+r+1}}{2m+r+1}\\
=&\frac{1}{2n+2} \sum_{r=0}^{2n+1} \frac{(-1)^{r}}{r! (2n+1-r)!}
  \sum_{m=0}^{\infty} \frac{(s+n)_{m} (\frac{1}{2})_{m}}{(s+\frac{1}
  {2})_{m} (n+2)_{m}}\\
 &\times \frac{(2m+1)(2m+2) \cdots (2m+2n+2)}{2m+r+1}
  \lambda^{2m+r+1}.
\endaligned$$   

  If $n=0$, then $r=0, 1$.
\flushpar
(1) $r=0$,
$$\sum_{m=0}^{\infty} \frac{(s)_{m} (\frac{1}{2})_{m}}
  {(s+\frac{1}{2})_{m} (2)_{m}} (2m+2) \lambda^{2m+1}
 =2 \lambda F(s, \frac{1}{2}; s+\frac{1}{2}; \lambda^2).$$
(2) $r=1$,
$$\sum_{m=0}^{\infty} \frac{(s)_{m} (\frac{1}{2})_{m}}
  {(s+\frac{1}{2})_{m} (2)_{m}} (2m+1) \lambda^{2m+2}
 =\frac{2s-1}{s-1} [F(s-1, \frac{1}{2}; s-\frac{1}{2};
  \lambda^2)-1].$$
On the other hand, when $n=0$, (4.27) is equal to
$$\frac{2s-1}{2s-2}[1-F(s-1, -\frac{1}{2}; s-\frac{1}{2}; 1)].$$
Hence,
$$\aligned
 &\int_{\lambda}^{1} (1-l)^{2n+1} F(s+n, n+\frac{3}{2}; 
  s+\frac{1}{2}; l^2) dl_{n=0}\\
=&-\lambda F(s, \frac{1}{2}; s+\frac{1}{2}; \lambda^2)+
  \frac{2s-1}{2s-2} F(s-1, \frac{1}{2}; s-\frac{1}{2}; \lambda^2)
  -\frac{2s-1}{2s-2} F(s-1, -\frac{1}{2}; s-\frac{1}{2}; 1).
\endaligned$$

  If $n=1$, then $r=0, 1, 2, 3$.
\flushpar
(1) $r=0$,
$$\aligned
 &\sum_{m=0}^{\infty} \frac{(s+1)_{m} (\frac{1}{2})_{m}}
  {(s+\frac{1}{2})_{m} (3)_{m}} (2m+2)(2m+3)(2m+4) \lambda^{2m+1}\\
=&8 \lambda F(s+1, \frac{3}{2}; s+\frac{1}{2}; \lambda^2)+16 \lambda
  F(s+1, \frac{1}{2}; s+\frac{1}{2}; \lambda^2).  
\endaligned$$
(2) $r=1$,
$$\sum_{m=0}^{\infty} \frac{(s+1)_{m} (\frac{1}{2})_{m}}
  {(s+\frac{1}{2})_{m} (3)_{m}} (2m+1)(2m+3)(2m+4) \lambda^{2m+2}
 =8 \frac{s-\frac{1}{2}}{s}[F(s, \frac{3}{2}; s-\frac{1}{2}; 
  \lambda^2)-1].$$
(3) $r=2$,
$$\sum_{m=0}^{\infty} \frac{(s+1)_{m} (\frac{1}{2})_{m}}
  {(s+\frac{1}{2})_{m} (3)_{m}} (2m+1)(2m+2)(2m+4) \lambda^{2m+3}
 =8 \lambda^{3} F(s+1, \frac{3}{2}; s+\frac{1}{2}; \lambda^2).$$
(4) $r=3$,
$$\aligned
 &\sum_{m=0}^{\infty} \frac{(s+1)_{m} (\frac{1}{2})_{m}}
  {(s+\frac{1}{2})_{m} (3)_{m}} (2m+1)(2m+2)(2m+3) \lambda^{2m+4}\\
=&8 \frac{s-\frac{1}{2}}{s} \lambda^{2} [F(s, \frac{3}{2}; 
  s-\frac{1}{2}; \lambda^2)-1]-16 \frac{(s-\frac{1}{2})(s-\frac{3}{2})}  
  {s(s-1)}[F(s-1, \frac{1}{2}; s-\frac{3}{2}; \lambda^{2})\\
 &-\frac{s-1}{2s-3} \lambda^2-1].
\endaligned$$

  On the other hand, when $n=1$, (4.27) is equal to
$$\frac{2}{3} \frac{(s-\frac{1}{2})(s-\frac{3}{2})}{s(s-1)}
  F(s-1, -\frac{3}{2}; s-\frac{3}{2}; 1)-\frac{2}{3} \frac{(s-
  \frac{1}{2})(s-\frac{3}{2})}{s(s-1)}+\frac{s-\frac{1}{2}}{s}.$$  
Therefore,
$$\aligned
 &\int_{\lambda}^{1} (1-l)^{2n+1} F(s+n, n+\frac{3}{2}; s+
  \frac{1}{2}; l^2) dl_{n=1}\\
=&\frac{2}{3} \frac{(s-\frac{1}{2})(s-\frac{3}{2})}{s(s-1)}
  F(s-1, -\frac{3}{2}; s-\frac{3}{2}; 1)-\frac{2}{3} \frac{(s-
  \frac{1}{2})(s-\frac{3}{2})}{s(s-1)} F(s-1, \frac{1}{2}; s-
  \frac{3}{2}; \lambda^{2})\\
 &-\frac{2}{3} \lambda F(s+1, \frac{1}{2}; s+\frac{1}{2}; \lambda^{2})  
  +(\frac{1}{3} \lambda^{2}+1) \frac{s-\frac{1}{2}}{s} F(s, 
  \frac{3}{2}; s-\frac{1}{2}; \lambda^{2})\\
 &-(\lambda^{3}+\frac{1}{3} \lambda) F(s+1, \frac{3}{2}; s+\frac{1}{2};
  \lambda^{2}). 
\endaligned$$

  Now, we need the following relation of Gauss between contiguous
functions (see \cite{Er}, p.103, (36)):  
$$(c-a-b)F-(c-a)F(a-1)+b(1-z)F(b+1)=0,$$
where $F$ denotes $F(a, b; c; z)$ and $F(a \pm 1)$, $F(b \pm 1)$, and
$F(c \pm 1)$ stands for $F(a \pm 1, b; c; z)$, $F(a, b \pm 1; c; z)$,
and $F(a, b; c \pm 1; z)$ respectively.

  We have
$$F(s-1, \frac{1}{2}; s-\frac{3}{2}; \lambda^{2})=\frac{1}{2}
  F(s-2, \frac{1}{2}; s-\frac{3}{2}; \lambda^{2})+\frac{1}{2}
  (1-\lambda^{2}) F(s-1, \frac{3}{2}; s-\frac{3}{2}; \lambda^{2}),$$
$$F(s+1, \frac{1}{2}; s+\frac{1}{2}; \lambda^{2})=\frac{1}{2}
  F(s, \frac{1}{2}; s+\frac{1}{2}; \lambda^{2})+\frac{1}{2} (1-
  \lambda^{2}) F(s+1, \frac{3}{2}; s+\frac{1}{2}; \lambda^{2}).$$

  By \cite{Er}, p.104, (46),
$$F(a, b; c; 1)=\frac{\Gamma(c) \Gamma(c-a-b)}{\Gamma(c-a) \Gamma(c-b)},
  \quad c \neq 0, -1, -2, \cdots, \text{Re}(c)>\text{Re}(a+b),$$
we have
$$F(s-1, -\frac{3}{2}; s-\frac{3}{2}; 1)=\frac{\Gamma(s-\frac{3}{2})}
  {-2 \sqrt{\pi} \Gamma(s)}, \quad
  F(s-1, -\frac{1}{2}; s-\frac{1}{2}; 1)=\frac{\Gamma(s-\frac{1}{2})}
  {\sqrt{\pi} \Gamma(s)}.$$
Therefore,
$$\aligned
 &\int_{\lambda}^{1} (1-l)^{2n+1} F(s+n, n+\frac{3}{2}; s+\frac{1}{2};
  l^{2}) dl_{n=0}\\
=&-\lambda F(s, \frac{1}{2}; s+\frac{1}{2}; \lambda^{2})+\frac{s-
  \frac{1}{2}}{s-1} F(s-1, \frac{1}{2}; s-\frac{1}{2}; \lambda^{2})  
  -\frac{1}{\sqrt{\pi}} \frac{\Gamma(s+\frac{1}{2})}{(s-1) \Gamma(s)},
\endaligned$$
and
$$\aligned
 &\int_{\lambda}^{1} (1-l)^{2n+1} F(s+n, n+\frac{3}{2}; s+\frac{1}{2};
  l^{2}) dl_{n=1}\\
=&-\frac{1}{3} \frac{\Gamma(s+\frac{1}{2})}{\sqrt{\pi} (s-1) 
  \Gamma(s+1)}-\frac{1}{3} \frac{(s-\frac{1}{2})(s-\frac{3}{2})}
  {s(s-1)} F(s-2, \frac{1}{2}; s-\frac{3}{2}; \lambda^{2})\\
 &-\frac{1}{3} \frac{(s-\frac{1}{2})(s-\frac{3}{2})}{s(s-1)}
  (1-\lambda^{2}) F(s-1, \frac{3}{2}; s-\frac{3}{2}; \lambda^{2})
  +(\frac{1}{3} \lambda^{2}+1) \frac{s-\frac{1}{2}}{s}
  F(s, \frac{3}{2}; s-\frac{1}{2}; \lambda^{2})\\
 &-\frac{1}{3} \lambda F(s, \frac{1}{2}; s+\frac{1}{2}; \lambda^{2})
  -\frac{2}{3} \lambda(\lambda^{2}+1) F(s+1, \frac{3}{2}; s+\frac{1}
  {2}; \lambda^{2}). 
\endaligned$$

  By \cite{Er}, p.111, (5),
$$F(a, b; 2b; \frac{4z}{(1+z)^{2}})=(1+z)^{2a} F(a, a+\frac{1}{2}-b;
  b+\frac{1}{2}; z^{2}),$$
we have
$$\aligned
 &F(s, \frac{1}{2}; s+\frac{1}{2}; \lambda^{2})=(1-\lambda)^{-2s}
  F(s, s; 2s; -u^{-1}),\\
 &F(s-1, \frac{1}{2}; s-\frac{1}{2}; \lambda^{2})=(1-\lambda)^{-2(s-1)}
  F(s-1, s-1; 2s-2; -u^{-1}),\\
 &F(s-2, \frac{1}{2}; s-\frac{3}{2}; \lambda^{2})=(1-\lambda)^{-2(s-2)}
  F(s-2, s-2; 2s-4; -u^{-1}),\\
 &F(s-1, \frac{3}{2}; s-\frac{3}{2}; \lambda^{2})=(1-\lambda)^{-2(s-1)}
  F(s-1, s-2; 2s-4; -u^{-1}),\\
 &F(s, \frac{3}{2}; s-\frac{1}{2}; \lambda^{2})=(1-\lambda)^{-2s}
  F(s, s-1; 2s-2; -u^{-1}),\\
 &F(s, \frac{1}{2}; s+\frac{1}{2}; \lambda^{2})=(1-\lambda)^{-2s}
  F(s, s; 2s; -u^{-1}),\\
 &F(s+1, \frac{3}{2}; s+\frac{1}{2}; \lambda^{2})=(1-\lambda)^{-2(s+1)}
  F(s+1, s; 2s; -u^{-1}).
\endaligned$$  
where $u=\frac{(\lambda-1)^{2}}{4 \lambda}$.
Hence,
$$\aligned
 &\int_{\lambda}^{1} (1-l)^{2n+1} F(s+n, n+\frac{3}{2}; s+\frac{1}{2};
  l^{2}) dl_{n=0}=-\lambda (1-\lambda)^{-2s} F(s, s; 2s; -u^{-1})\\
+&\frac{s-\frac{1}{2}}{s-1} (1-\lambda)^{-2(s-1)} 
  F(s-1, s-1; 2s-2; -u^{-1})-\frac{1}{\sqrt{\pi}} 
  \frac{1}{s-1} \frac{\Gamma(s+\frac{1}{2})}{\Gamma(s)}.
\endaligned$$
$$\aligned
 &\int_{\lambda}^{1} (1-l)^{2n+1} F(s+n, n+\frac{3}{2}; s+\frac{1}{2};
  l^{2}) dl_{n=1}=-\frac{1}{3} \frac{1}{\sqrt{\pi}} \frac{1}{s-1} 
  \frac{\Gamma(s+\frac{1}{2})}{\Gamma(s+1)}\\
 &-\frac{1}{3} \frac{(s-\frac{1}{2})(s-\frac{3}{2})}{s(s-1)} 
  (1-\lambda)^{-2(s-2)} F(s-2, s-2; 2s-4; -u^{-1})\\  
 &-\frac{1}{3} \frac{(s-\frac{1}{2})(s-\frac{3}{2})}{s(s-1)}
  (1-\lambda^{2})(1-\lambda)^{-2(s-1)} F(s-1, s-2; 2s-4; -u^{-1})\\
 &+(\frac{1}{3} \lambda^{2}+1) \frac{s-\frac{1}{2}}{s} (1-\lambda)^{-2s}
  F(s, s-1; 2s-2; -u^{-1})\\
 &-\frac{1}{3} \lambda (1-\lambda)^{-2s} F(s, s; 2s; -u^{-1})
  -\frac{2}{3} \lambda (\lambda^{2}+1) (1-\lambda)^{-2(s+1)}
  F(s+1, s; 2s; -u^{-1}). 
\endaligned$$

  Similarly,
$$\aligned
 &\sum_{m=0}^{\infty} \frac{(2-s+n)_{m} (n+\frac{1}{2})_{m}}
  {m! (\frac{3}{2}-s)_{m}} B_{\lambda}(2m+2-2s, 2n+2)\\
=&\sum_{m=0}^{\infty} \frac{(2-s+n)_{m} (n+\frac{1}{2})_{m}}
  {m! (\frac{3}{2}-s)_{m}} B_{1}(2m+2-2s, 2n+2)
  I_{\lambda}(2m+2-2s, 2n+2)\\
=&\frac{(2n+1)! \Gamma(2-2s)}{\Gamma(4-2s+2n)}   
  \sum_{m=0}^{\infty} \frac{(1-s)_{m} (n+\frac{1}{2})_{m}}
  {(\frac{5}{2}-s+n)_{m} m!} I_{\lambda}(2m+2-2s, 2n+2).
\endaligned$$
By Lemma 4.11, we have
$$I_{\lambda}(2m+2-2s, 2n+2)=(2m+2-2s)_{2n+2}
  \sum_{r=0}^{2n+1} \frac{(-1)^{r}}{r! (2n+1-r)!}
  \frac{\lambda^{2m+2-2s+r}}{2m+2-2s+r}.$$
Thus,
$$\aligned
 &\sum_{m=0}^{\infty} \frac{(1-s)_{m} (n+\frac{1}{2})_{m}}
  {(\frac{5}{2}-s+n)_{m} m!} I_{\lambda}(2m+2-2s, 2n+2)\\
=&\sum_{r=0}^{2n+1} \frac{(-1)^{r}}{r! (2n+1-r)!}
  \sum_{m=0}^{\infty} \frac{(1-s)_{m} (n+\frac{1}{2})_{m}}
  {(\frac{5}{2}-s+n)_{m} m!}\\
 &\times \frac{(2m+2-2s)(2m+3-2s) \cdots (2m+2n+3-2s)}{2m+2-2s+r} 
  \lambda^{2m+2-2s+r}.
\endaligned$$

  If $n=0$, then $r=0, 1$. 
\flushpar
(1) $r=0$,
$$\sum_{m=0}^{\infty} \frac{(1-s)_{m} (\frac{1}{2})_{m}}
  {(\frac{5}{2}-s)_{m} m!} (2m+3-2s) \lambda^{2m+2-2s}
 =(3-2s) \lambda^{2-2s} F(1-s, \frac{1}{2}; \frac{3}{2}-s;
  \lambda^{2}).$$
(2) $r=1$,
$$\sum_{m=0}^{\infty} \frac{(1-s)_{m} (\frac{1}{2})_{m}}
  {(\frac{5}{2}-s)_{m} m!} (2m+2-2s) \lambda^{2m+3-2s}
 =2(1-s) \lambda^{3-2s} F(2-s, \frac{1}{2}; \frac{5}{2}-s;
  \lambda^{2}).$$

  On the other hand, when $n=0$, (4.28) is equal to
$$\frac{1}{(3-2s)(2-2s)} F(1-s, \frac{1}{2}; \frac{5}{2}-s; 1).$$  
Thus,
$$\aligned
 &\int_{\lambda}^{1} (1-l)^{2n+1} l^{1-2s} F(n-s+2, n+\frac{1}{2};
  \frac{3}{2}-s; l^{2}) dl_{n=0}\\
=&\frac{1}{(3-2s)(2-2s)} F(1-s, \frac{1}{2}; \frac{5}{2}-s; 1)
  -\frac{\lambda^{2-2s}}{2-2s} F(1-s, \frac{1}{2}; \frac{3}{2}
  -s; \lambda^{2})\\
 &+\frac{\lambda^{3-2s}}{3-2s} F(2-s, \frac{1}{2};
  \frac{5}{2}-s; \lambda^{2}).
\endaligned$$

  If $n=1$, then $r=0, 1, 2, 3$. 
\flushpar
(1) $r=0$,
$$\aligned
 &\sum_{m=0}^{\infty} \frac{(1-s)_{m} (\frac{3}{2})_{m}}
  {(\frac{7}{2}-s)_{m} m!} (2m+3-2s)(2m+4-2s)(2m+5-2s)
  \lambda^{2m+2-2s}\\
=&8(1-s)(\frac{3}{2}-s)(\frac{5}{2}-s) \lambda^{2-2s}
  F(2-s, \frac{3}{2}; \frac{3}{2}-s; \lambda^{2})\\
 &+8(\frac{3}{2}-s)(\frac{5}{2}-s) \lambda^{2-2s}
  F(1-s, \frac{3}{2}; \frac{3}{2}-s; \lambda^{2}).  
\endaligned$$
(2) $r=1$,
$$\aligned
 &\sum_{m=0}^{\infty} \frac{(1-s)_{m} (\frac{3}{2})_{m}}
  {(\frac{7}{2}-s)_{m} m!} (2m+2-2s)(2m+4-2s)(2m+5-2s)
  \lambda^{2m+3-2s}\\
=&8(1-s)(2-s)(\frac{5}{2}-s) \lambda^{3-2s}
  F(3-s, \frac{3}{2}; \frac{5}{2}-s; \lambda^{2}).
\endaligned$$
(3) $r=2$,
$$\aligned
 &\sum_{m=0}^{\infty} \frac{(1-s)_{m} (\frac{3}{2})_{m}}
  {(\frac{7}{2}-s)_{m} m!} (2m+2-2s)(2m+3-2s)(2m+5-2s)
  \lambda^{2m+4-2s}\\
=&8(1-s)(\frac{3}{2}-s)(\frac{5}{2}-s) \lambda^{4-2s}
  F(2-s, \frac{3}{2}; \frac{3}{2}-s; \lambda^{2}).
\endaligned$$
(3) $r=3$,
$$\aligned
 &\sum_{m=0}^{\infty} \frac{(1-s)_{m} (\frac{3}{2})_{m}}
  {(\frac{7}{2}-s)_{m} m!} (2m+2-2s)(2m+3-2s)(2m+4-2s)
  \lambda^{2m+5-2s}\\
=&8(1-s)(2-s)(\frac{5}{2}-s) \lambda^{5-2s} F(3-s, \frac{3}{2};
  \frac{5}{2}-s; \lambda^{2})\\
 &-8(1-s)(2-s) \lambda^{5-2s} F(3-s, \frac{3}{2}; \frac{7}{2}-s; 
  \lambda^{2}).
\endaligned$$

  On the other hand, when $n=1$, (4.28) is equal to
$$\frac{6 \Gamma(2-2s)}{\Gamma(6-2s)} F(1-s, \frac{3}{2};
  \frac{7}{2}-s; 1).$$  
Thus, 
$$\aligned
 &\int_{\lambda}^{1} (1-l)^{2n+1} l^{1-2s}
  F(n-s+2, n+\frac{1}{2}; \frac{3}{2}-s; l^{2})_{n=1}\\
=&\frac{6 \Gamma(2-2s)}{\Gamma(6-2s)} F(1-s, \frac{3}{2};
  \frac{7}{2}-s; 1)  
  -\frac{1}{2} \frac{1}{(2-s)} \lambda^{2-2s} (1+3 \lambda^{2})
  F(2-s, \frac{3}{2}; \frac{3}{2}-s; \lambda^{2})\\
 &+\frac{1}{2} \frac{1}{(\frac{3}{2}-s)} \lambda^{3-2s}
  (3+\lambda^{2}) F(3-s, \frac{3}{2}; \frac{5}{2}-s; \lambda^{2})
  -\frac{1}{2} \frac{1}{(1-s)(2-s)} \lambda^{2-2s}\\
 &\times F(1-s, \frac{3}{2}; \frac{3}{2}-s; \lambda^{2})
  -\frac{1}{2} \frac{1}{(\frac{3}{2}-s)(\frac{5}{2}-s)}
  \lambda^{5-2s} F(3-s, \frac{3}{2}; \frac{7}{2}-s; \lambda^{2}).   
\endaligned$$

  Now, we need the following relation of Gauss between contiguous
functions (see \cite{Er}, p.103, (33)):
$$(c-a-b) F+a(1-z) F(a+1)-(c-b) F(b-1)=0.$$  
We have
$$F(1-s, \frac{3}{2}; \frac{3}{2}-s; \lambda^{2})
 =(1-s)(1-\lambda^{2}) F(2-s, \frac{3}{2}; \frac{3}{2}-s; \lambda^{2})
 +s F(1-s, \frac{1}{2}; \frac{3}{2}-s; \lambda^{2}).$$
$$F(3-s, \frac{3}{2}; \frac{7}{2}-s; \lambda^{2})
 =(3-s)(1-\lambda^{2}) F(4-s, \frac{3}{2}; \frac{7}{2}-s; \lambda^{2})
 -(2-s) F(3-s, \frac{1}{2}; \frac{7}{2}-s; \lambda^{2}).$$
By \cite{Er}, p.104, (46),
$$F(a, b; c; 1)=\frac{\Gamma(c) \Gamma(c-a-b)}{\Gamma(c-a) \Gamma(c-b)},
  \quad c \neq 0, -1, -2, \cdots, \quad \text{Re}(c)>\text{Re}(a+b),$$
we have
$$F(1-s, \frac{1}{2}; \frac{5}{2}-s; 1)
 =2 \frac{\Gamma(\frac{5}{2}-s)}{\sqrt{\pi} \Gamma(2-s)}, \quad
  F(1-s, \frac{3}{2}; \frac{7}{2}-s; 1)
 =\frac{4}{3} \frac{\Gamma(\frac{7}{2}-s)}{\sqrt{\pi} \Gamma(2-s)}.$$
Therefore,
$$\aligned
 &\int_{\lambda}^{1} (1-l)^{2n+1} l^{1-2s} F(n-s+2, n+\frac{1}{2};
  \frac{3}{2}-s; l^{2}) dl_{n=0}\\
=&\frac{1}{2-2s} \frac{\Gamma(\frac{3}{2}-s)}{\sqrt{\pi} \Gamma(2-s)}
  -\frac{\lambda^{2-2s}}{2-2s} F(1-s, \frac{1}{2}; \frac{3}{2}-s; 
  \lambda^{2})
  +\frac{\lambda^{3-2s}}{3-2s} F(2-s, \frac{1}{2}; \frac{5}{2}-s;
  \lambda^{2}), 
\endaligned$$
and
$$\aligned
 &\int_{\lambda}^{1} (1-l)^{2n+1} l^{1-2s} F(n-s+2, n+\frac{1}{2};
  \frac{3}{2}-s; l^{2}) dl_{n=1}\\
=&\frac{1}{2-2s} \frac{\Gamma(\frac{3}{2}-s)}{\sqrt{\pi} \Gamma(3-s)}
  -\frac{1}{2-s} (1+\lambda^{2}) \lambda^{2-2s} F(2-s, \frac{3}{2};
  \frac{3}{2}-s; \lambda^{2})\\
 &+\frac{1}{2} \frac{1}{(\frac{3}{2}-s)} (3+\lambda^{2}) \lambda^{3-2s}
  F(3-s, \frac{3}{2}; \frac{5}{2}-s; \lambda^{2})\\
 &-\frac{1}{2} \frac{s}{(1-s)(2-s)} \lambda^{2-2s} F(1-s, \frac{1}{2};
  \frac{3}{2}-s; \lambda^{2})\\
 &-\frac{1}{2} \frac{(3-s)}{(\frac{3}{2}-s)(\frac{5}{2}-s)}
  (1-\lambda^{2}) \lambda^{5-2s} F(4-s, \frac{3}{2}; \frac{7}{2}-s;
  \lambda^{2})\\
 &+\frac{1}{2} \frac{(2-s)}{(\frac{3}{2}-s)(\frac{5}{2}-s)}
  \lambda^{5-2s} F(3-s, \frac{1}{2}; \frac{7}{2}-s; \lambda^{2}).    
\endaligned$$

  By \cite{Er}, p.111, (5),
$$F(a, b; 2b; \frac{4z}{(1+z)^{2}})=(1+z)^{2a}
  F(a, a+\frac{1}{2}-b; b+\frac{1}{2}; z^{2}),$$  
we have
$$\aligned
 &F(1-s, \frac{1}{2}; \frac{3}{2}-s; \lambda^{2})
 =(1-\lambda)^{-2(1-s)} F(1-s, 1-s; 2-2s; -u^{-1}),\\
 &F(2-s, \frac{1}{2}; \frac{5}{2}-s; \lambda^{2})
 =(1-\lambda)^{-2(2-s)} F(2-s, 2-s; 4-2s; -u^{-1}),\\
 &F(2-s, \frac{3}{2}; \frac{3}{2}-s; \lambda^{2})
 =(1-\lambda)^{-2(2-s)} F(2-s, 1-s; 2-2s; -u^{-1}),\\
 &F(3-s, \frac{3}{2}; \frac{5}{2}-s; \lambda^{2})
 =(1-\lambda)^{-2(3-s)} F(3-s, 2-s; 4-2s; -u^{-1}),\\
 &F(1-s, \frac{1}{2}; \frac{3}{2}-s; \lambda^{2})
 =(1-\lambda)^{-2(1-s)} F(1-s, 1-s; 2-2s; -u^{-1}),\\
 &F(4-s, \frac{3}{2}; \frac{7}{2}-s; \lambda^{2})
 =(1-\lambda)^{-2(4-s)} F(4-s, 3-s; 6-2s; -u^{-1}),\\
 &F(3-s, \frac{1}{2}; \frac{7}{2}-s; \lambda^{2})
 =(1-\lambda)^{-2(3-s)} F(3-s, 3-s; 6-2s; -u^{-1}).
\endaligned$$ 
Therefore,
$$\aligned
 &\int_{\lambda}^{1} (1-l)^{2n+1} l^{1-2s} F(n-s+2, n+\frac{1}{2};
  \frac{3}{2}-s; l^{2}) dl_{n=0}\\
=&\frac{1}{2-2s} \frac{\Gamma(\frac{3}{2}-s)}{\sqrt{\pi} \Gamma(2-s)}
 -\frac{1}{2-2s} \lambda^{2-2s} (1-\lambda)^{-2(1-s)} F(1-s, 1-s;
  2-2s; -u^{-1})\\
 &+\frac{1}{3-2s} \lambda^{3-2s} (1-\lambda)^{-2(2-s)} F(2-s, 2-s;
  4-2s; -u^{-1}). 
\endaligned$$ 
$$\aligned
 &\int_{\lambda}^{1} (1-l)^{2n+1} l^{1-2s} F(n-s+2, n+\frac{1}{2};
  \frac{3}{2}-s; l^{2}) dl_{n=1}\\
=&\frac{1}{2-2s} \frac{\Gamma(\frac{3}{2}-s)}{\sqrt{\pi} \Gamma(3-s)}
  -\frac{1}{2-s} (1+\lambda^{2}) \lambda^{2-2s} (1-\lambda)^{-2(2-s)}
  F(2-s, 1-s; 2-2s; -u^{-1})\\
 &+\frac{1}{2} \frac{1}{(\frac{3}{2}-s)} (3+\lambda^{2}) \lambda^{3-2s}
  (1-\lambda)^{-2(3-s)} F(3-s, 2-s; 4-2s; -u^{-1})\\
 &-\frac{1}{2} \frac{s}{(1-s)(2-s)} \lambda^{2-2s} 
  (1-\lambda)^{-2(1-s)} F(1-s, 1-s; 2-2s; -u^{-1})\\
 &-\frac{1}{2} \frac{(3-s)}{(\frac{3}{2}-s)(\frac{5}{2}-s)}
  (1-\lambda^{2}) \lambda^{5-2s} (1-\lambda)^{-2(4-s)}
  F(4-s, 3-s; 6-2s; -u^{-1})\\
 &+\frac{1}{2} \frac{(2-s)}{(\frac{3}{2}-s)(\frac{5}{2}-s)}
  \lambda^{5-2s} (1-\lambda)^{-2(3-s)} F(3-s, 3-s; 6-2s; -u^{-1}).
\endaligned$$

  Since $u=\frac{(\lambda-1)^{2}}{4 \lambda}$, we have
$$\aligned
 &\frac{\Gamma(\frac{1}{2}-s)}{\sqrt{\pi} \Gamma(2-s)}
  \lambda^{s} (1-\lambda)^{-2} \int_{\lambda}^{1} (1-l)^{2n+1}
  F(s+n, n+\frac{3}{2}; s+\frac{1}{2}; l^{2}) dl_{n=0}\\
=&\frac{\Gamma(\frac{1}{2}-s)}{\sqrt{\pi} \Gamma(2-s)}
  [-(4u)^{-(s+1)} F(s, s; 2s; -u^{-1})+\frac{s-\frac{1}{2}}{s-1} 
  (4u)^{-s}\\
 &\times F(s-1, s-1; 2s-2; -u^{-1})-\frac{1}{s-1} 
  \frac{\Gamma(s+\frac{1}{2})}{\sqrt{\pi} \Gamma(s)}
  \lambda^{s} (1-\lambda)^{-2}],
\endaligned$$
and
$$\aligned
 &\frac{2 \Gamma(s-\frac{1}{2})}{\sqrt{\pi} \Gamma(s)} \lambda^{s}
  (1-\lambda)^{-2} \int_{\lambda}^{1} (1-l)^{2n+1} l^{1-2s}
  F(n-s+2, n+\frac{1}{2}; \frac{3}{2}-s; l^{2}) dl_{n=0}\\
=&\frac{\Gamma(s-\frac{1}{2})}{\sqrt{\pi} \Gamma(s)}
  [\frac{1}{1-s} \frac{\Gamma(\frac{3}{2}-s)}{\sqrt{\pi} \Gamma(2-s)}
  \lambda^{s} (1-\lambda)^{-2}-\frac{1}{1-s} (4u)^{-(2-s)}\\
 &\times F(1-s, 1-s; 2-2s; -u^{-1})+\frac{1}{\frac{3}{2}-s} 
  (4u)^{-(3-s)} F(2-s, 2-s; 4-2s; -u^{-1})].
\endaligned$$
Note that
$$\aligned
 &\lambda^{s} (1-\lambda)^{-2} \left[ \frac{\Gamma(\frac{1}{2}-s)}
  {\sqrt{\pi} \Gamma(2-s)} \frac{(-1)}{\sqrt{\pi}} \frac{1}{s-1}
  \frac{\Gamma(s+\frac{1}{2})}{\Gamma(s)}+\frac{\Gamma(s-\frac{1}{2})}
  {\sqrt{\pi} \Gamma(s)} \frac{1}{\sqrt{\pi}} \frac{1}{1-s}
  \frac{\Gamma(\frac{3}{2}-s)}{\Gamma(2-s)} \right]\\
=&\frac{\lambda^{s} (1-\lambda)^{-2}}{\pi (1-s) \Gamma(s) \Gamma(2-s)}
  [\Gamma(\frac{1}{2}-s) \Gamma(s+\frac{1}{2})+\Gamma(s-\frac{1}{2})
  \Gamma(\frac{3}{2}-s)]\\
=&0.
\endaligned$$

  Now, we need the following formulas:
$$\gamma F-\beta z F(\beta+1, \gamma+1)-\gamma F(\alpha-1)=0.$$  
$$(\gamma-1) F(\gamma-1)-\alpha F(\alpha+1)-(\gamma-\alpha-1) F=0.$$
We have
$$(2s-1) F(s, s-1; 2s-1; z)-(s-1) z F(s, s; 2s; z)-
  (2s-1) F(s-1, s-1; 2s-1; z)=0,$$
and
$$(2s-2) F(s-1, s-1; 2s-2; z)-(s-1) F(s, s-1; 2s-1; z)
  -(s-1) F(s-1, s-1; 2s-1; z)=0.$$
Thus,
$$\frac{4s-2}{s-1} F(s, s-1; 2s-1; z)=z F(s, s; 2s; z)+
  \frac{4s-2}{s-1} F(s-1, s-1; 2s-2; z).$$
Hence, we have
$$\aligned
 &\frac{s-\frac{1}{2}}{s-1} F(s, s-1; 2s-1; -u^{-1})\\
=&-(4u)^{-1} F(s, s; 2s; -u^{-1})+\frac{s-\frac{1}{2}}{s-1}
  F(s-1, s-1; 2s-2; -u^{-1}).
\endaligned$$  
Similarly,
$$\aligned
 &\frac{\frac{3}{2}-s}{1-s} F(2-s, 1-s; 3-2s; -u^{-1})\\
=&-(4u)^{-1} F(2-s, 2-s; 4-2s; -u^{-1})+\frac{\frac{3}{2}-s}
  {1-s} F(1-s, 1-s; 2-2s; -u^{-1}).
\endaligned$$
Therefore,
$$\aligned
  \frac{1}{4 \pi^{2}} \phi(\frac{1}{2} \rho, \frac{1}{2} 
  \rho^{\prime}; 0, s)
=&\frac{1}{1-s} \frac{\Gamma(\frac{3}{2}-s)}{\sqrt{\pi} \Gamma(2-s)}
  (4u)^{-s} F(s, s-1; 2s-1; -u^{-1})\\ 
 &+\frac{1}{s-1} \frac{\Gamma(s-\frac{1}{2})}{\sqrt{\pi} \Gamma(s)}
  (4u)^{s-2} F(2-s, 1-s; 3-2s; -u^{-1}).
\endaligned$$

  For $n=1$, note that
$$\frac{3! \Gamma(\frac{1}{2}-s)}{\Gamma(3-s) \Gamma(\frac{3}{2})}
 =\frac{12 \Gamma(\frac{1}{2}-s)}{\sqrt{\pi} \Gamma(3-s)},$$  
we have
$$\aligned
 &\frac{12 \Gamma(\frac{1}{2}-s)}{\sqrt{\pi} \Gamma(3-s)}
  \lambda^{s} (1-\lambda)^{-2} \int_{\lambda}^{1} (1-l)^{2n+1}
  F(s+n, n+\frac{3}{2}; s+\frac{1}{2}; l^{2}) dl_{n=1}\\ 
=&-4 \frac{\Gamma(\frac{1}{2}-s)}{\sqrt{\pi} \Gamma(3-s)}
  \frac{1}{s-1} \frac{\Gamma(s+\frac{1}{2})}{\sqrt{\pi} 
  \Gamma(s+1)} \lambda^{s} (1-\lambda)^{-2}\\
 &-4 \frac{\Gamma(\frac{5}{2}-s)}{s(s-1) \sqrt{\pi} \Gamma(3-s)}
  \lambda (4u)^{-(s-1)} F(s-2, s-2; 2s-4; -u^{-1})\\
 &-4 \frac{\Gamma(\frac{5}{2}-s)}{s(s-1) \sqrt{\pi} \Gamma(3-s)}
  (1-\lambda^{2}) (4u)^{-s} F(s-1, s-2; 2s-4; -u^{-1})\\
 &-4 \frac{\Gamma(\frac{3}{2}-s)}{s \sqrt{\pi} \Gamma(3-s)}
  (\lambda+\frac{3}{\lambda}) (4u)^{-(s+1)} F(s, s-1; 2s-2; -u^{-1})\\
 &-4 \frac{\Gamma(\frac{1}{2}-s)}{\sqrt{\pi} \Gamma(3-s)}    
  (4u)^{-(s+1)} F(s, s; 2s; -u^{-1})\\
 &-8 \frac{\Gamma(\frac{1}{2}-s)}{\sqrt{\pi} \Gamma(3-s)}
  (\lambda+\frac{1}{\lambda}) (4u)^{-(s+2)} F(s+1, s; 2s; -u^{-1}).
\endaligned$$

  Similarly,
$$\frac{3! \Gamma(s-\frac{1}{2})}{\Gamma(s+1) \Gamma(\frac{5}{2})}
 =\frac{8 \Gamma(s-\frac{1}{2})}{\sqrt{\pi} \Gamma(s+1)},$$
and
$$\aligned
 &\frac{8 \Gamma(s-\frac{1}{2})}{\sqrt{\pi} \Gamma(s+1)}
  \lambda^{s} (1-\lambda)^{-2} \int_{\lambda}^{1} (1-l)^{2n+1}
  l^{1-2s} F(n-s+2, n+\frac{1}{2}; \frac{3}{2}-s; l^{2}) dl_{n=1}\\
=&\frac{8 \Gamma(s-\frac{1}{2})}{\sqrt{\pi} \Gamma(s+1)}
  \times [\frac{1}{2} \frac{1}{(1-s)} \frac{\Gamma(\frac{3}{2}-s)}
  {\sqrt{\pi} \Gamma(3-s)} \lambda^{s} (1-\lambda)^{-2}\\
 &-\frac{1}{2-s} (\lambda+\frac{1}{\lambda}) (4u)^{-(3-s)} 
  F(2-s, 1-s; 2-2s; -u^{-1})\\
 &+\frac{1}{2} \frac{1}{(\frac{3}{2}-s)} (\lambda+\frac{3}{\lambda})
  (4u)^{-(4-s)} F(3-s, 2-s; 4-2s; -u^{-1})\\
 &-\frac{1}{2} \frac{s}{(1-s)(2-s)} (4u)^{-(2-s)} 
  F(1-s, 1-s; 2-2s; -u^{-1})\\
 &-\frac{1}{2} \frac{(3-s)}{(\frac{3}{2}-s)(\frac{5}{2}-s)}
  (1-\lambda^{2}) (4u)^{-(5-s)} F(4-s, 3-s; 6-2s; -u^{-1})\\
 &+\frac{1}{2} \frac{(2-s)}{(\frac{3}{2}-s)(\frac{5}{2}-s)}
  \lambda (4u)^{-(4-s)} F(3-s, 3-s; 6-2s; -u^{-1})].
\endaligned$$

  Note that
$$\aligned
 &\left[\frac{-4 \Gamma(\frac{1}{2}-s)}{\sqrt{\pi} \Gamma(3-s)}
  \frac{1}{s-1} \frac{\Gamma(s+\frac{1}{2})}{\sqrt{\pi} \Gamma(s+1)}
  +\frac{-4 \Gamma(s-\frac{1}{2})}{\sqrt{\pi} \Gamma(s+1)} 
  \frac{1}{s-1} \frac{\Gamma(\frac{3}{2}-s)}{\sqrt{\pi} 
  \Gamma(3-s)} \right] \lambda^{s} (1-\lambda)^{-2}\\
=&\frac{-4 \lambda^{s} (1-\lambda)^{-2}}{\pi (s-1) \Gamma(s+1)
  \Gamma(3-s)} \left[ \Gamma(\frac{1}{2}-s) \Gamma(s+\frac{1}{2})+
  \Gamma(s-\frac{1}{2}) \Gamma(\frac{3}{2}-s) \right]\\
=&0.
\endaligned$$  

  For $k=1$, if $n=0$, $\frac{(-1)_{n} (1)_{n}}{2^{2n} (n!)^{2}}=1$.
If $n=1$, $\frac{(-1)_{n} (1)_{n}}{2^{2n} (n!)^{2}}=-\frac{1}{4}$.
Let us add the terms of $n=0$ and $n=1$, we have
$$\frac{\Gamma(\frac{1}{2}-s)}{\sqrt{\pi} \Gamma(2-s)} (4u)^{-s} S_{I}
+\frac{\Gamma(s-\frac{1}{2})}{\sqrt{\pi} \Gamma(s)} (4u)^{s-3} S_{II},$$
where
$$\aligned
  S_{I}
=&\frac{\frac{1}{2}-s}{1-s} F(s, s-1; 2s-1; -u^{-1})  
  +\frac{(\frac{3}{2}-s)(\frac{1}{2}-s)}{s(s-1)(2-s)} 4u \lambda
  F(s-2, s-2; 2s-4; -u^{-1})\\
 &+\frac{(\frac{3}{2}-s)(\frac{1}{2}-s)}{s(s-1)(2-s)} (1-\lambda^{2})
  F(s-1, s-2; 2s-4; -u^{-1})\\
 &+\frac{\frac{1}{2}-s}{s(2-s)} (\lambda+\frac{3}{\lambda})
  \frac{1}{4u} F(s, s-1; 2s-2; -u^{-1})\\
 &+\frac{1}{2-s} \frac{1}{4u} F(s, s; 2s; -u^{-1})
  +\frac{1}{2-s} 2(\lambda+\frac{1}{\lambda}) \frac{1}{(4u)^{2}}
  F(s+1, s; 2s; -u^{-1}), 
\endaligned$$
and
$$\aligned
  S_{II}
=&\frac{1}{s-1} 4u F(2-s, 1-s; 3-2s; -u^{-1})\\
 &+\frac{2}{s(2-s)} (\lambda+\frac{1}{\lambda})
  F(2-s, 1-s; 2-2s; -u^{-1})\\
 &-\frac{1}{s(\frac{3}{2}-s)} (\lambda+\frac{3}{\lambda}) 
  (4u)^{-1} F(3-s, 2-s; 4-2s; -u^{-1})\\
 &+\frac{1}{(1-s)(2-s)} 4u F(1-s, 1-s; 2-2s; -u^{-1})\\
 &+\frac{(3-s)}{s(\frac{3}{2}-s)(\frac{5}{2}-s)}
  (1-\lambda^{2}) (4u)^{-2} F(4-s, 3-s; 6-2s; -u^{-1})\\
 &-\frac{(2-s)}{s(\frac{3}{2}-s)(\frac{5}{2}-s)} \lambda
  (4u)^{-1} F(3-s, 3-s; 6-2s; -u^{-1}). 
\endaligned$$

  Since $u=\frac{(\lambda-1)^{2}}{4 \lambda}$, we have
$$\lambda=\frac{1}{(\sqrt{u+1}+\sqrt{u})^{2}}.$$
Hence,
$$\lambda+\frac{1}{\lambda}=4u+2, \quad
  \lambda+\frac{3}{\lambda}=4(2u+1+\sqrt{u+1} \sqrt{u}),$$
$$1-\lambda^{2}=\lambda(\frac{1}{\lambda}-\lambda)
 =\frac{4 \sqrt{u} \sqrt{u+1}}{(\sqrt{u+1}+\sqrt{u})^{2}}.$$
Thus, we have
$$\aligned
  S_{I}
=&\frac{s-\frac{1}{2}}{s-1} F(s, s-1; 2s-1; -u^{-1})\\
 &+\frac{(\frac{3}{2}-s)(\frac{1}{2}-s)}{s(s-1)(2-s)}
  \frac{4}{(\sqrt{1+\frac{1}{u}}+1)^{2}} F(s-2, s-2; 2s-4; -u^{-1})\\ 
 &+\frac{(\frac{3}{2}-s)(\frac{1}{2}-s)}{s(s-1)(2-s)}
  \frac{4 \sqrt{1+\frac{1}{u}}}{(\sqrt{1+\frac{1}{u}}+1)^{2}}
  F(s-1, s-2; 2s-4; -u^{-1})\\  
 &+\frac{(\frac{1}{2}-s)}{s(2-s)} (2+\frac{1}{u}+\sqrt{1+\frac{1}{u}})
  F(s, s-1; 2s-2; -u^{-1})\\
 &+\frac{1}{2-s} \frac{1}{4u} F(s, s; 2s; -u^{-1})
  +\frac{1}{2-s} (\frac{1}{2u}+\frac{1}{4 u^{2}})
  F(s+1, s; 2s; -u^{-1}).  
\endaligned$$
$$\aligned
  S_{II}
=&\frac{1}{s-1} 4u F(2-s, 1-s; 3-2s; -u^{-1})\\
 &+\frac{1}{s(2-s)} 4(2u+1) F(2-s, 1-s; 2-2s; -u^{-1})\\
 &-\frac{1}{s(\frac{3}{2}-s)} (2+\frac{1}{u}+\sqrt{1+\frac{1}{u}})
  F(3-s, 2-s; 4-2s; -u^{-1})\\
 &+\frac{1}{(1-s)(2-s)} 4u F(1-s, 1-s; 2-2s; -u^{-1})\\
 &+\frac{(3-s)}{s(\frac{3}{2}-s)(\frac{5}{2}-s)}
  \frac{\sqrt{1+\frac{1}{u}}}{(1+\sqrt{1+\frac{1}{u}})^{2}}
  \frac{1}{4 u^2} F(4-s, 3-s; 6-2s; -u^{-1})\\
 &-\frac{(2-s)}{s(\frac{3}{2}-s)(\frac{5}{2}-s)}
  \frac{1}{(1+\sqrt{1+\frac{1}{u}})^{2}} \frac{1}{4 u^2}
  F(3-s, 3-s; 6-2s; -u^{-1}).
\endaligned$$
Set $z=-u^{-1}$ and $x=\sqrt{1-z}$, then
$$\frac{u}{u+1}=\frac{1}{x^{2}}, \quad 
  \frac{u^2}{u+1}=\frac{1}{x^2 (x^2-1)}.$$
$$\aligned
  S_{I}
=&\frac{s-\frac{1}{2}}{s-1} F(s, s-1; 2s-1; 1-x^{2})\\
 &+\frac{(\frac{3}{2}-s)(\frac{1}{2}-s)}{s(s-1)(2-s)}
  \frac{4}{(x+1)^{2}} F(s-2, s-2; 2s-4; 1-x^{2})\\
 &+\frac{(\frac{3}{2}-s)(\frac{1}{2}-s)}{s(s-1)(2-s)}
  \frac{4x}{(x+1)^{2}} F(s-1, s-2; 2s-4; 1-x^{2})\\
 &+\frac{(\frac{1}{2}-s)}{s(2-s)} (x^{2}+x+1)
  F(s, s-1; 2s-2; 1-x^{2})\\
 &+\frac{1}{2-s} \frac{1}{4} (x^{2}-1) F(s, s; 2s; 1-x^{2})
  +\frac{1}{2-s} \frac{1}{4} (x^{4}-1) F(s+1, s; 2s; 1-x^{2}).  
\endaligned$$
$$\aligned
  S_{II}
=&\frac{1}{s-1} \frac{4}{x^2-1} F(2-s, 1-s; 3-2s; 1-x^2)\\
 &+\frac{1}{s(2-s)} \frac{4(x^2+1)}{x^2-1} F(2-s, 1-s; 2-2s; 1-x^2)\\
 &-\frac{1}{s(\frac{3}{2}-s)} (x^2+x+1) F(3-s, 2-s; 4-2s; 1-x^2)\\
 &+\frac{1}{(1-s)(2-s)} \frac{4}{x^2-1} F(1-s, 1-s; 2-2s; 1-x^2)\\
 &+\frac{(3-s)}{s(\frac{3}{2}-s)(\frac{5}{2}-s)} \frac{1}{4}
  x(x-1)^{2} F(4-s, 3-s; 6-2s; 1-x^2)\\
 &-\frac{(2-s)}{s(\frac{3}{2}-s)(\frac{5}{2}-s)}
  \frac{1}{4} (x-1)^{2} F(3-s, 3-s; 6-2s; 1-x^{2}). 
\endaligned$$

  At first, we calculate $S_{I}$. 
  It is known that (see \cite{Er}, p.103, (29))
$$(c-b) F(b-1)+(2b-c-bz+az) F+b(z-1) F(b+1)=0.$$
Set $a=s$, $b=s$, $c=2s$ and $z=1-x^{2}$, then
$$F(s, s+1; 2s; 1-x^{2})=\frac{1}{x^{2}} F(s, s-1; 2s; 1-x^{2}),$$
i.e.,
$$F(s+1, s; 2s; 1-x^{2})=\frac{1}{x^{2}} F(s, s-1; 2s; 1-x^{2}).$$
Set $a=s-1$, $b=s-1$, $c=2s-2$ and $z=1-x^{2}$, then
$$F(s-1, s; 2s-2; 1-x^{2})=\frac{1}{x^{2}} F(s-1, s-2; 2s-2; 1-x^{2}),$$
i.e.,
$$F(s, s-1; 2s-2; 1-x^{2})=\frac{1}{x^2} F(s-1, s-2; 2s-2; 1-x^2).$$

  In the formula (see \cite{Er}, p.103, (34))
$$c[a-(c-b)z] F-ac(1-z) F(a+1)+(c-a)(c-b)z F(c+1)=0,$$
set $a=s-2$, $b=s-2$, $c=2s-4$ and $z=1-x^{2}$, then
$$\aligned
 &F(s-2, s-2; 2s-4; 1-x^{2})\\
=&F(s-1, s-2; 2s-4; 1-x^{2})+\frac{x^{2}-1}{2 x^{2}} 
  F(s-2, s-2; 2s-3; 1-x^{2}).
\endaligned$$
Now, we have
$$\aligned
  S_{I}
=&\frac{s-\frac{1}{2}}{s-1} F(s, s-1; 2s-1; 1-x^{2})\\
 &+\frac{(\frac{3}{2}-s)(\frac{1}{2}-s)}{s(s-1)(2-s)}
  \frac{4}{x+1} F(s-1, s-2; 2s-4; 1-x^{2})\\
 &+\frac{(\frac{3}{2}-s)(\frac{1}{2}-s)}{s(s-1)(2-s)}
  \frac{2(x-1)}{x^{2} (x+1)} F(s-2, s-2; 2s-3; 1-x^{2})\\
 &+\frac{(\frac{1}{2}-s)}{s(2-s)} \frac{x^{2}+x+1}{x^{2}}
  F(s-1, s-2; 2s-2; 1-x^{2})\\
 &+\frac{1}{2-s} \frac{1}{4} (x^{2}-1) F(s, s; 2s; 1-x^{2})\\
 &+\frac{1}{2-s} \frac{1}{4} \frac{x^{4}-1}{x^{2}}
  F(s, s-1; 2s; 1-x^{2}). 
\endaligned$$

  It is known that
$$\gamma F-\beta z F(\beta+1, \gamma+1)-\gamma F(\alpha-1)=0.$$
Set $\alpha=s$, $\beta=s-1$, $\gamma=2s-1$ and $z=1-x^{2}$, then
$$F(s, s; 2s; 1-x^{2})=\frac{1}{1-x^{2}} \frac{2s-1}{s-1}
 [F(s, s-1; 2s-1; 1-x^{2})-F(s-1, s-1; 2s-1; 1-x^{2})].$$
Set $\alpha=s-1$, $\beta=s-1$, $\gamma=2s-1$ and $z=1-x^{2}$, then
$$\aligned
 &F(s, s-1; 2s; 1-x^{2})=F(s-1, s; 2s; 1-x^{2})\\
=&\frac{1}{1-x^{2}} \frac{2s-1}{s-1}[F(s-1, s-1; 2s-1; 1-x^{2})
  -F(s-2, s-1; 2s-1; 1-x^{2})]. 
\endaligned$$
Thus, we have
$$\aligned
  S_{I}
=&\frac{(\frac{1}{2}-s)(s-\frac{3}{2})}{(2-s)(s-1)}
  F(s, s-1; 2s-1; 1-x^{2})\\
 &+\frac{(\frac{3}{2}-s)(\frac{1}{2}-s)}{s(s-1)(2-s)}
  \frac{4}{x+1} F(s-1, s-2; 2s-4; 1-x^{2})\\
 &+\frac{(\frac{3}{2}-s)(\frac{1}{2}-s)}{s(s-1)(2-s)}
  \frac{2(x-1)}{x^{2} (x+1)} F(s-2, s-2; 2s-3; 1-x^{2})\\
 &+\frac{(\frac{1}{2}-s)}{s(2-s)} \frac{x^2+x+1}{x^2} 
  F(s-1, s-2; 2s-2; 1-x^2)\\
 &+\frac{(\frac{1}{2}-s)}{(2-s)(s-1)} \frac{1}{2 x^2}
  F(s-1, s-1; 2s-1; 1-x^2)\\
 &-\frac{(\frac{1}{2}-s)}{(2-s)(s-1)} \frac{x^2+1}{2 x^2}
  F(s-1, s-2; 2s-1; 1-x^2). 
\endaligned$$

  It is known that (see \cite{Er}, p.103, (37))
$$(b-a)(1-z)F-(c-a)F(a-1)+(c-b)F(b-1)=0.$$
Set $a=s$, $b=s-1$, $c=2s-1$ and $z=1-x^2$, then
$$F(s, s-1; 2s-1; 1-x^2)=\frac{1-s}{x^2} F(s-1, s-1; 2s-1; 1-x^2)
 +\frac{s}{x^2} F(s, s-2; 2s-1; 1-x^2).$$

  By \cite{Er}, p.103, (30),
$$c(c-1)(z-1)F(c-1)+c[c-1-(2c-a-b-1)z] F+(c-a)(c-b)zF(c+1)=0,$$
set $a=s-1$, $b=s-2$, $c=2s-3$ and $z=1-x^2$, then
$$\aligned
 &F(s-1, s-2; 2s-4; 1-x^2)\\
=&F(s-1, s-2; 2s-3; 1-x^2)+\frac{(s-1)}{4(s-\frac{3}{2})} 
  \frac{1-x^2}{x^2} F(s-1, s-2; 2s-2; 1-x^2).
\endaligned$$
Hence,
$$\aligned
  S_{I}
=&\frac{s(\frac{1}{2}-s)(\frac{3}{2}-s)}{(2-s)(1-s)}
  \frac{1}{x^2} F(s, s-2; 2s-1; 1-x^2)\\
 &+\frac{(\frac{3}{2}-s)(\frac{1}{2}-s)}{s(s-1)(2-s)}
  \frac{4}{x+1} F(s-1, s-2; 2s-3; 1-x^2)\\
 &+\frac{(\frac{3}{2}-s)(\frac{1}{2}-s)}{s(s-1)(2-s)}
  \frac{2(x-1)}{x^2 (x+1)} F(s-2, s-2; 2s-3; 1-x^2)\\
 &+\frac{(\frac{1}{2}-s)}{s(2-s)} \frac{x+2}{x} 
  F(s-1, s-2; 2s-2; 1-x^2)\\
 &-\frac{(\frac{1}{2}-s)^{2}}{s-1} \frac{1}{x^2}
  F(s-1, s-1; 2s-1; 1-x^2)\\
 &-\frac{(\frac{1}{2}-s)}{(2-s)(s-1)} \frac{x^2+1}{2 x^2}
  F(s-1, s-2; 2s-1; 1-x^2).  
\endaligned$$

  It is known that (see \cite{Er}, p.103, (35))
$$(c-a-1)F+aF(a+1)-(c-1)F(c-1)=0.$$
Set $a=s-2$, $b=s-1$, $c=2s-1$ and $z=1-x^2$, then
$$\aligned
 &F(s-1, s-1; 2s-1; 1-x^2)\\
=&\frac{2(s-1)}{s-2} F(s-2, s-1; 2s-2; 1-x^2)
  -\frac{s}{s-2} F(s-2, s-1; 2s-1; 1-x^2).
\endaligned$$
Set $a=s-2$, $b=s-2$, $c=2s-2$ and $z=1-x^2$, then
$$\aligned
 &F(s-2, s-2; 2s-3; 1-x^2)\\
=&\frac{(s-1)}{2(s-\frac{3}{2})} F(s-2, s-2; 2s-2; 1-x^2)+
  \frac{(s-2)}{2(s-\frac{3}{2})} F(s-1, s-2; 2s-2; 1-x^2).
\endaligned$$
Set $a=s-1$, $b=s-2$, $c=2s-1$ and $z=1-x^2$, then
$$F(s, s-2; 2s-1; 1-x^2)=2 F(s-1, s-2; 2s-2; 1-x^2)-
  F(s-1, s-2; 2s-1; 1-x^2).$$
Hence,
$$\aligned
  S_{I}
=&\left[\frac{-s(\frac{1}{2}-s)}{(s-1)(s-2)} \frac{1}{x^2}+
  \frac{(\frac{1}{2}-s)}{(s-1)(s-2)} \frac{x^2+1}{2 x^2}\right]
  F(s-1, s-2; 2s-1; 1-x^2)\\
 &+\left[\frac{(\frac{1}{2}-s)}{(s-1)(s-2)} \frac{1}{x^2}+
  \frac{(\frac{1}{2}-s)}{s(s-1)} \frac{(x-1)}{x^2 (x+1)}-
  \frac{(\frac{1}{2}-s)}{s(s-2)} \frac{x+2}{x}\right]\\
 &\times F(s-1, s-2; 2s-2; 1-x^2)\\
 &+\frac{(\frac{1}{2}-s)(\frac{3}{2}-s)}{s(s-1)(2-s)}
  \frac{4}{x+1} F(s-1, s-2; 2s-3; 1-x^2)\\
 &+\frac{(\frac{1}{2}-s)}{s(s-2)} \frac{(x-1)}{x^2 (x+1)}
  F(s-2, s-2; 2s-2; 1-x^2). 
\endaligned$$

  We know that (see \cite{Er}, p.103, (35))
$$(c-a-1)F+aF(a+1)-(c-1)F(c-1)=0.$$
Set $a=s-2$, $b=s-2$, $c=2s-1$ and $z=1-x^2$, then
$$\aligned
 &F(s-2, s-2; 2s-2; 1-x^2)\\
=&\frac{s}{2(s-1)} F(s-2, s-2; 2s-1; 1-x^2)
 +\frac{s-2}{2(s-1)} F(s-1, s-2; 2s-1; 1-x^2).
\endaligned$$
For $F(s-1, s-2; 2s-3; 1-x^2)=F(s-2, s-1; 2s-3; 1-x^2)$,
set $a=s-2$, $b=s-1$, $c=2s-2$ and $z=1-x^2$, then
$$\aligned
 &F(s-1, s-2; 2s-3; 1-x^2)\\
=&\frac{s-1}{2(s-\frac{3}{2})} F(s-1, s-2; 2s-2; 1-x^2)+
  \frac{s-2}{2(s-\frac{3}{2})} F(s-1, s-1; 2s-2; 1-x^2).
\endaligned$$
Therefore, we have
$$\aligned
  S_{I}
=&\left[\frac{-s(\frac{1}{2}-s)}{(s-1)(s-2)} \frac{1}{x^2}+
  \frac{(\frac{1}{2}-s)}{(s-1)(s-2)} \frac{x^2+1}{2 x^2}+
  \frac{(\frac{1}{2}-s)}{s(s-1)} \frac{(x-1)}{2 x^2 (x+1)}\right]\\
 &\times F(s-1, s-2; 2s-1; 1-x^2)\\
 &+\left[\frac{(\frac{1}{2}-s)}{(s-1)(s-2)} \frac{1}{x^2}
  +\frac{(\frac{1}{2}-s)}{s(s-1)} \frac{(x-1)}{x^2 (x+1)}
  -\frac{(\frac{1}{2}-s)}{s(s-2)} \frac{(x+2)}{x}
  +\frac{(\frac{1}{2}-s)}{s(s-2)} \frac{2}{x+1}\right]\\
 &\times F(s-1, s-2; 2s-2; 1-x^2)\\
 &+\frac{(\frac{1}{2}-s)}{s(s-1)} \frac{2}{x+1}
  F(s-1, s-1; 2s-2; 1-x^2)\\
 &+\frac{(\frac{1}{2}-s)}{(s-1)(s-2)} \frac{(x-1)}{2 x^2 (x+1)}
  F(s-2, s-2; 2s-1; 1-x^2). 
\endaligned$$

  Since (see \cite{Er}, p.103, (38))
$$c(1-z)F-cF(a-1)+(c-b)zF(c+1)=0.$$
Set $a=s-1$, $b=s-1$, $c=2s-2$ and $z=1-x^2$, then
$$\aligned
 &F(s-1, s-1; 2s-2; 1-x^2)\\
=&\frac{1}{x^2} F(s-1, s-2; 2s-2; 1-x^2)
 +\frac{x^2-1}{2 x^2} F(s-1, s-1; 2s-1; 1-x^2).
\endaligned$$
Thus,
$$\aligned
  S_{I}
=&\left[\frac{-s(\frac{1}{2}-s)}{(s-1)(s-2)} \frac{1}{x^2}+
  \frac{(\frac{1}{2}-s)}{(s-1)(s-2)} \frac{(x^2+1)}{2 x^2}+
  \frac{(\frac{1}{2}-s)}{s(s-1)} \frac{(x-1)}{2 x^2 (x+1)}\right]\\
 &\times F(s-1, s-2; 2s-1; 1-x^2)\\
 &+\frac{(\frac{1}{2}-s)}{s(s-2)} \frac{-(x-1)(x^2+2x+2)}{x^2 (x+1)}
  F(s-1, s-2; 2s-2; 1-x^2)\\
 &+\frac{(\frac{1}{2}-s)}{s(s-1)} \frac{(x-1)}{x^2}
  F(s-1, s-1; 2s-1; 1-x^2)\\
 &+\frac{(\frac{1}{2}-s)}{(s-1)(s-2)} \frac{(x-1)}{2 x^2 (x+1)}
  F(s-2, s-2; 2s-1; 1-x^2). 
\endaligned$$

  It is known that (see \cite{Er}, p.103, (35))
$$(c-a-1)F+aF(a+1)-(c-1)F(c-1)=0.$$
Set $a=s-2$, $b=s-1$, $c=2s-1$ and $z=1-x^2$, then
$$\aligned
 &F(s-1, s-2; 2s-2; 1-x^2)\\
=&\frac{s}{2(s-1)} F(s-1, s-2; 2s-1; 1-x^2)+
  \frac{(s-2)}{2(s-1)} F(s-1, s-1; 2s-1; 1-x^2).
\endaligned$$
Hence,
$$\aligned
  S_{I}
=&[\frac{-s(\frac{1}{2}-s)}{(s-1)(s-2)} \frac{1}{x^2}+
  \frac{(\frac{1}{2}-s)}{(s-1)(s-2)} \frac{(x^2+1)}{2 x^2}+
  \frac{(\frac{1}{2}-s)}{s(s-1)} \frac{(x-1)}{2 x^2 (x+1)}\\
 &+\frac{(\frac{1}{2}-s)}{(s-1)(s-2)} \frac{-(x-1)(x^2+2x+2)}
  {2 x^2 (x+1)}] F(s-1, s-2; 2s-1; 1-x^2)\\
 &+\frac{(\frac{1}{2}-s)}{s(s-1)} \frac{-(x-1)}{2(x+1)}
  F(s-1, s-1; 2s-1; 1-x^2)\\
 &+\frac{(\frac{1}{2}-s)}{(s-1)(s-2)} \frac{(x-1)}{2 x^2 (x+1)}
  F(s-2, s-2; 2s-1; 1-x^2). 
\endaligned\tag 4.30$$

  In the following formula (see \cite{Er}, p.103, (36))
$$(c-a-b)F-(c-a)F(a-1)+b(1-z)F(b+1)=0,$$
set $a=s-1$, $b=s-2$, $c=2s-1$ and $z=1-x^2$, then
$$\aligned
 &2 F(s-1, s-2; 2s-1; 1-x^2)\\
=&-(s-2) x^2 F(s-1, s-1, 2s-1; 1-x^2)+s F(s-2, s-2; 2s-1; 1-x^2).
\endaligned$$
Now, the sum of the latter two terms in (4.30) is equal to
$$\aligned
 &\frac{(\frac{1}{2}-s)}{s(s-1)(s-2)} \frac{(x-1)}{2 x^2 (x+1)}
  [-(s-2) x^2 F(s-1, s-1; 2s-1; 1-x^2)\\
 &+s F(s-2, s-2; 2s-1; 1-x^2)]\\
=&\frac{(\frac{1}{2}-s)}{s(s-1)(s-2)} \frac{(x-1)}{x^2 (x+1)}   
  F(s-1, s-2; 2s-1; 1-x^2).
\endaligned$$
Thus,
$$\aligned
  S_{I}
=&[\frac{-s(\frac{1}{2}-s)}{(s-1)(s-2)} \frac{1}{x^2}
  +\frac{(\frac{1}{2}-s)}{(s-1)(s-2)} \frac{(x^2+1)}{2 x^2}
  +\frac{(\frac{1}{2}-s)}{s(s-1)} \frac{(x-1)}{2 x^2 (x+1)}\\
 &+\frac{(\frac{1}{2}-s)}{s(s-1)(s-2)} \frac{(x-1)}{x^2 (x+1)}
  +\frac{(\frac{1}{2}-s)}{(s-1)(s-2)} \frac{-(x-1)(x^2+2x+2)}
  {2 x^2 (x+1)}]\\
 &\times F(s-1, s-2; 2s-1; 1-x^2).
\endaligned$$
Note that
$$\frac{x^2+1}{2 x^2}-\frac{(x-1)(x^2+2x+2)}{2 x^2 (x+1)}
 =\frac{x+3}{2 x^2 (x+1)},$$  
and
$$\frac{(\frac{1}{2}-s)}{s(s-1)} \frac{(x-1)}{2 x^2 (x+1)}+
  \frac{(\frac{1}{2}-s)}{s(s-1)(s-2)} \frac{(x-1)}{x^2 (x+1)}
 =\frac{(\frac{1}{2}-s)}{(s-1)(s-2)} \frac{(x-1)}{2 x^2 (x+1)}.$$
We have
$$\aligned
  S_{I}
=&\left[\frac{-s(\frac{1}{2}-s)}{(s-1)(s-2)} \frac{1}{x^2}+
  \frac{(\frac{1}{2}-s)}{(s-1)(s-2)} \frac{(x+3)}{2 x^2 (x+1)}+
  \frac{(\frac{1}{2}-s)}{(s-1)(s-2)} \frac{(x-1)}{2 x^2 (x+1)}\right]\\
 &\times F(s-1, s-2; 2s-1; 1-x^2)\\
=&\frac{\frac{1}{2}-s}{2-s} \frac{1}{x^2} F(s-1, s-2; 2s-1; 1-x^2)\\
=&\frac{\frac{1}{2}-s}{2-s} \frac{u}{u+1} F(s-1, s-2; 2s-1; -u^{-1}).
\endaligned$$
Therefore, we have
$$\frac{\Gamma(\frac{1}{2}-s)}{\sqrt{\pi} \Gamma(2-s)} (4u)^{-s} S_{I}
 =\frac{\Gamma(\frac{3}{2}-s)}{\sqrt{\pi} \Gamma(3-s)} 4^{-s}
  (u+1)^{-1} u^{-(s-1)} F(s-1, s-2; 2s-1; -u^{-1}).$$

  Now, we compute $S_{II}$. It is known that (see \cite{Er}, 
p.103, (28))
$$(c-a)F(a-1)+(2a-c-az+bz)F+a(z-1)F(a+1)=0.$$
Set $a=3-s$, $b=3-s$, $c=6-2s$ and $z=1-x^2$, then
$$F(4-s, 3-s; 6-2s; 1-x^2)=\frac{1}{x^2} F(2-s, 3-s; 6-2s; 1-x^2).$$
By \cite{Er}, p.103, (39), 
$$[a-1-(c-b-1)z]F+(c-a)F(a-1)-(c-1)(1-z)F(c-1)=0,$$
set $a=3-s$, $b=3-s$, $c=6-2s$ and $z=1-x^2$, then
$$\aligned
 &F(3-s, 3-s; 6-2s; 1-x^2)\\ 
=&\frac{2(\frac{5}{2}-s)}{2-s} F(3-s, 3-s; 5-2s; 1-x^2)
 -\frac{3-s}{2-s} \frac{1}{x^2} F(2-s, 3-s; 6-2s; 1-x^2).
\endaligned$$
Thus,
$$\aligned
  S_{II}
=&\frac{1}{s-1} \frac{4}{x^2-1} F(2-s, 1-s; 3-2s; 1-x^2)\\
 &+\frac{1}{s(2-s)} \frac{4(x^2+1)}{x^2-1} F(2-s, 1-s; 2-2s; 1-x^2)\\
 &-\frac{1}{s(\frac{3}{2}-s)} (x^2+x+1) F(3-s, 2-s; 4-2s; 1-x^2)\\
 &+\frac{1}{(1-s)(2-s)} \frac{4}{x^2-1} F(1-s, 1-s; 2-2s; 1-x^2)\\
 &+\frac{(3-s)}{s(\frac{3}{2}-s)(\frac{5}{2}-s)} \frac{(x-1)^{2} 
  (x+1)}{4 x^2} F(3-s, 2-s; 6-2s; 1-x^2)\\
 &-\frac{1}{s(\frac{3}{2}-s)} \frac{1}{2} (x-1)^{2}
  F(3-s, 3-s; 5-2s; 1-x^2).  
\endaligned$$
 
  We know that (see \cite{Er}, p.103, (30)) 
$$c(c-1)(z-1)F(c-1)+c[c-1-(2c-a-b-1)z]F+(c-a)(c-b)zF(c+1)=0.$$
Set $a=3-s$, $b=2-s$, $c=5-2s$ and $z=1-x^2$, then
$$\aligned
  F(3-s, 2-s; 6-2s; 1-x^2)
=&\frac{-4(\frac{5}{2}-s)}{3-s} \frac{x^2}{x^2-1}
  F(3-s, 2-s; 4-2s; 1-x^2)\\
 &+\frac{4(\frac{5}{2}-s)}{3-s} \frac{x^2}{x^2-1} 
  F(3-s, 2-s; 5-2s; 1-x^2).
\endaligned$$
It is known that (see \cite{Er}, p.103, (35))
$$(c-a-1)F+aF(a+1)-(c-1)F(c-1)=0.$$
Set $a=2-s$, $b=3-s$, $c=5-2s$ and $z=1-x^2$, then
$$F(3-s, 3-s; 5-2s; 1-x^2)=2 F(3-s, 2-s; 4-2s; 1-x^2)
 -F(3-s, 2-s; 5-2s; 1-x^2).$$
Hence,
$$\aligned
  S_{II}
=&\frac{1}{s-1} \frac{4}{x^2-1} F(2-s, 1-s; 3-2s; 1-x^2)\\
 &+\frac{1}{s(2-s)} \frac{4(x^2+1)}{x^2-1} F(2-s, 1-s; 2-2s; 1-x^2)\\
 &-\frac{1}{s(\frac{3}{2}-s)} (2 x^2+1) F(3-s, 2-s; 4-2s; 1-x^2)\\
 &+\frac{1}{(1-s)(2-s)} \frac{4}{x^2-1} F(1-s, 1-s; 2-2s; 1-x^2)\\
 &+\frac{1}{s(\frac{3}{2}-s)} \frac{1}{2}(x^2-1) 
  F(3-s, 2-s; 5-2s; 1-x^2).
\endaligned$$

  By \cite{Er}, p.103, (39)
$$[a-1-(c-b-1)z]F+(c-a)F(a-1)-(c-1)(1-z)F(c-1)=0,$$
set $a=3-s$, $b=2-s$, $c=5-2s$ and $z=1-x^2$, then
$$F(3-s, 2-s; 5-2s; 1-x^2)=-\frac{1}{x^2} F(2-s, 2-s; 5-2s; 1-x^2)
 +2 F(3-s, 2-s; 4-2s; 1-x^2).$$
We have
$$\aligned
  S_{II}
=&\frac{1}{s-1} \frac{4}{x^2-1} F(2-s, 1-s; 3-2s; 1-x^2)\\
 &+\frac{1}{s(2-s)} \frac{4(x^2+1)}{x^2-1} F(2-s, 1-s; 2-2s; 1-x^2)\\
 &-\frac{1}{s(\frac{3}{2}-s)} (x^2+2) F(3-s, 2-s; 4-2s; 1-x^2)\\
 &+\frac{1}{(1-s)(2-s)} \frac{4}{x^2-1} F(1-s, 1-s; 2-2s; 1-x^2)\\
 &-\frac{1}{s(\frac{3}{2}-s)} \frac{x^2-1}{2 x^2}
  F(2-s, 2-s; 5-2s; 1-x^2).
\endaligned$$

  It is known that (see \cite{Er}, p.103, (28))
$$(c-a)F(a-1)+(2a-c-az+bz)F+a(z-1)F(a+1)=0.$$
Set $a=2-s$, $b=2-s$, $c=4-2s$ and $z=1-x^2$, then
$$F(3-s, 2-s; 4-2s; 1-x^2)=\frac{1}{x^2} F(2-s, 1-s; 4-2s; 1-x^2).$$
By
$$\gamma F-\beta z F(\beta+1, \gamma+1)-\gamma F(\alpha-1)=0,$$
set $\alpha=2-s$, $\beta=1-s$, $\gamma=4-2s$ and $z=1-x^2$, then
$$\aligned
 &F(2-s, 2-s; 5-2s; 1-x^2)\\
=&\frac{2(2-s)}{1-s} \frac{1}{x^2-1}[F(1-s, 1-s; 4-2s; 1-x^2)-
  F(2-s, 1-s; 4-2s; 1-x^2)].
\endaligned$$
So,
$$\aligned
  S_{II}
=&\frac{1}{s-1} \frac{4}{x^2-1} F(2-s, 1-s; 3-2s; 1-x^2)\\
 &+\frac{1}{s(2-s)} \frac{4(x^2+1)}{x^2-1} F(2-s, 1-s; 2-2s; 1-x^2)\\
 &+\left[\frac{1}{(1-s)(\frac{3}{2}-s)} \frac{1}{x^2}-
  \frac{1}{s(\frac{3}{2}-s)}\right] F(2-s, 1-s; 4-2s; 1-x^2)\\
 &+\frac{1}{(1-s)(2-s)} \frac{4}{x^2-1} F(1-s, 1-s; 2-2s; 1-x^2)\\
 &-\frac{(2-s)}{s(1-s)(\frac{3}{2}-s)} \frac{1}{x^2}
  F(1-s, 1-s; 4-2s; 1-x^2).
\endaligned$$

  We know that (see \cite{Er}, p.103, (38))
$$c(1-z)F-cF(a-1)+(c-b)zF(c+1)=0.$$
Set $a=2-s$, $b=1-s$, $c=2-2s$ and $z=1-x^2$, then
$$\aligned
 &F(2-s, 1-s; 2-2s; 1-x^2)\\
=&\frac{1}{x^2} F(1-s, 1-s; 2-2s; 1-x^2)+\frac{x^2-1}{2 x^2}
  F(2-s, 1-s; 3-2s; 1-x^2).
\endaligned$$
Thus,
$$\aligned
  S_{II}
=&\left[\frac{1}{s-1} \frac{4}{x^2-1}+\frac{1}{s(2-s)} 
  \frac{2(x^2+1)}{x^2} \right] F(2-s, 1-s; 3-2s; 1-x^2)\\
 &+\left[\frac{1}{(1-s)(2-s)} \frac{4}{x^2-1}+\frac{1}{s(2-s)} 
  \frac{4(x^2+1)}{x^2(x^2-1)} \right] F(1-s, 1-s; 2-2s; 1-x^2)\\
 &+\left[\frac{1}{(1-s)(\frac{3}{2}-s)} \frac{1}{x^2}-
  \frac{1}{s(\frac{3}{2}-s)} \right] F(2-s, 1-s; 4-2s; 1-x^2)\\
 &-\frac{(2-s)}{s(1-s)(\frac{3}{2}-s)} \frac{1}{x^2}
  F(1-s, 1-s; 4-2s; 1-x^2).
\endaligned$$

  It is known that 
$$\gamma F-\beta z F(\beta+1, \gamma+1)-\gamma F(\alpha-1)=0.$$
Set $\alpha=1-s$, $\beta=-s$, $\gamma=3-2s$ and $z=1-x^2$, then
$$\aligned
 &F(1-s, 1-s; 4-2s; 1-x^2)\\ 
=&\frac{2(\frac{3}{2}-s)}{s} \frac{1}{x^2-1}
  [F(1-s, -s; 3-2s; 1-x^2)-F(-s, -s; 3-2s; 1-x^2)].
\endaligned$$
Set $\alpha=2-s$, $\beta=-s$, $\gamma=3-2s$ and $z=1-x^2$, then
$$\aligned
 &F(2-s, 1-s; 4-2s; 1-x^2)\\
=&\frac{2(\frac{3}{2}-s)}{s} \frac{1}{x^2-1}
 [F(2-s, -s; 3-2s; 1-x^2)-F(1-s, -s; 3-2s; 1-x^2)]. 
\endaligned$$
Therefore,
$$\aligned
  S_{II}
=&\left[\frac{1}{s-1} \frac{4}{x^2-1}+\frac{1}{s(2-s)}
  \frac{2(x^2+1)}{x^2}\right] F(2-s, 1-s; 3-2s; 1-x^2)\\
 &+\left[\frac{1}{(1-s)(2-s)} \frac{4}{x^2-1}+\frac{1}{s(2-s)}
  \frac{4(x^2+1)}{x^2(x^2-1)}\right] F(1-s, 1-s; 2-2s; 1-x^2)\\
 &+\left[\frac{1}{s(1-s)} \frac{2}{x^2(x^2-1)}-\frac{1}{s^2}
  \frac{2}{x^2-1}\right] F(2-s, -s; 3-2s; 1-x^2)\\  
 &+\frac{(2-s)}{s^2(1-s)} \frac{2}{x^2(x^2-1)} F(-s, -s; 3-2s; 1-x^2)\\
 &+\left[-\frac{1}{s^2(1-s)} \frac{4}{x^2(x^2-1)}+
  \frac{1}{s^2} \frac{2}{x^2-1}\right] F(1-s, -s; 3-2s; 1-x^2).
\endaligned$$

  It is known that (see \cite{Er}, p.103, (28))
$$(c-a)F(a-1)+(2a-c-az+bz)F+a(z-1)F(a+1)=0.$$
Set $a=1-s$, $b=1-s$, $c=3-2s$ and $z=1-x^2$, then
$$\aligned
 &F(2-s, 1-s; 3-2s; 1-x^2)\\
=&\frac{2-s}{1-s} \frac{1}{x^2} F(1-s, -s; 3-2s; 1-x^2)
 -\frac{1}{1-s} \frac{1}{x^2} F(1-s, 1-s; 3-2s; 1-x^2).
\endaligned$$
Set $a=1-s$, $b=-s$, $c=3-2s$ and $z=1-x^2$, then
$$\aligned
 &F(2-s, -s; 3-2s; 1-x^2)\\ 
=&\frac{2-s}{1-s} \frac{1}{x^2} F(-s, -s; 3-2s; 1-x^2)
 +\frac{1}{1-s} \frac{x^2-2}{x^2} F(1-s, -s; 3-2s; 1-x^2).
\endaligned$$

  Note that
$$\aligned
 &\frac{1}{(1-s)(2-s)} \frac{4}{x^2-1}+\frac{1}{s(2-s)}
  \frac{4(x^2+1)}{x^2(x^2-1)}\\
=&\frac{1}{s(1-s)(2-s)} \frac{4}{x^2-1}+\frac{1}{s(2-s)} 
  \frac{4}{x^2(x^2-1)},
\endaligned$$
and
$$-\frac{1}{s^2(1-s)} \frac{4}{x^2(x^2-1)}+\frac{1}{s^2}
  \frac{2}{x^2-1}-\frac{1}{s^2(1-s)} \frac{2(x^2-2)}{x^2(x^2-1)}
 =\frac{-1}{s(1-s)} \frac{2}{x^2-1}.$$
We have
$$\aligned
  S_{II}
=&[-\frac{1}{s(1-s)} \frac{2}{x^2-1}-\frac{(2-s)}{(1-s)^2}
  \frac{4}{x^2(x^2-1)}+\frac{1}{s(1-s)} \frac{2(x^2+1)}{x^4}\\
 &+\frac{1}{s(1-s)^2} \frac{2(x^2-2)}{x^4(x^2-1)}]
  F(1-s, -s; 3-2s; 1-x^2)\\
 &+[\frac{1}{(1-s)^2} \frac{4}{x^2(x^2-1)}-\frac{1}{s(1-s)(2-s)}
  \frac{2(x^2+1)}{x^4}] F(1-s, 1-s; 3-2s; 1-x^2)\\
 &+[\frac{1}{s(2-s)} \frac{4}{x^2(x^2-1)}+\frac{1}{s(1-s)(2-s)}
  \frac{4}{x^2-1}] F(1-s, 1-s; 2-2s; 1-x^2)\\
 &+\frac{(2-s)}{s(1-s)^2} \frac{2}{x^4(x^2-1)}
  F(-s, -s; 3-2s; 1-x^2). 
\endaligned$$

  It is known that (see \cite{Er}, p.103, (38))
$$c(1-z)F-cF(a-1)+(c-b)zF(c+1)=0.$$
Set $a=1-s$, $b=1-s$, $c=2-2s$ and $z=1-x^2$, then
$$\aligned
 &F(1-s, 1-s; 2-2s; 1-x^2)\\ 
=&\frac{1}{x^2} F(-s, 1-s; 2-2s; 1-x^2)
 +\frac{x^2-1}{2 x^2} F(1-s, 1-s; 3-2s; 1-x^2).
\endaligned$$
Note that
$$\aligned
 &\frac{1}{(1-s)^2} \frac{4}{x^2(x^2-1)}-\frac{1}{s(1-s)(2-s)}
  \frac{2(x^2+1)}{x^4}+\frac{1}{s(2-s)} \frac{2}{x^4}+
  \frac{1}{s(1-s)(2-s)} \frac{2}{x^2}\\
=&\frac{1}{(1-s)^2} \frac{4}{x^2(x^2-1)}-
  \frac{1}{(1-s)(2-s)} \frac{2}{x^4}.
\endaligned$$
Thus,
$$\aligned
  S_{II}
=&[-\frac{1}{s(1-s)} \frac{2}{x^2-1}-\frac{(2-s)}{(1-s)^2}
  \frac{4}{x^2(x^2-1)}+\frac{1}{s(1-s)} \frac{2(x^2+1)}{x^4}\\
 &+\frac{1}{s(1-s)^2} \frac{2(x^2-2)}{x^4(x^2-1)}]
  F(1-s, -s; 3-2s; 1-x^2)\\
 &+[\frac{1}{(1-s)^2} \frac{4}{x^2(x^2-1)}-\frac{1}{(1-s)(2-s)}
  \frac{2}{x^4}] F(1-s, 1-s; 3-2s; 1-x^2)\\
 &+[\frac{1}{s(2-s)} \frac{4}{x^4(x^2-1)}+\frac{1}{s(1-s)(2-s)}
  \frac{4}{x^2(x^2-1)}] F(-s, 1-s; 2-2s; 1-x^2)\\
 &+\frac{(2-s)}{s(1-s)^2} \frac{2}{x^4(x^2-1)}
  F(-s, -s; 3-2s; 1-x^2).  
\endaligned$$

  We know that (see \cite{Er}, p.103, (35))
$$(c-a-1)F+aF(a+1)-(c-1)F(c-1)=0.$$
Set $a=-s$, $b=1-s$, $c=3-2s$ and $z=1-x^2$, then
$$\aligned
 &F(-s, 1-s; 2-2s; 1-x^2)\\
=&\frac{2-s}{2(1-s)} F(-s, 1-s; 3-2s; 1-x^2)
 -\frac{s}{2(1-s)} F(1-s, 1-s; 3-2s; 1-x^2). 
\endaligned$$
Note that
$$\aligned
 &\frac{1}{(1-s)^2} \frac{4}{x^2(x^2-1)}-\frac{1}{(1-s)(2-s)}
  \frac{2}{x^4}-\frac{1}{(1-s)(2-s)} \frac{2}{x^4(x^2-1)}\\
 &-\frac{1}{(1-s)^2 (2-s)} \frac{2}{x^2(x^2-1)}
  =\frac{1}{(1-s)^2} \frac{2}{x^2(x^2-1)}.  
\endaligned$$
Therefore,
$$\aligned
  S_{II}
=&[-\frac{1}{s(1-s)} \frac{2}{x^2-1}-\frac{(2-s)}{(1-s)^2}
  \frac{4}{x^2(x^2-1)}+\frac{1}{s(1-s)} \frac{2(x^2+1)}{x^4}\\
 &+\frac{1}{s(1-s)^2} \frac{2(x^2-2)}{x^4(x^2-1)}+\frac{1}{s(1-s)}
  \frac{2}{x^4(x^2-1)}+\frac{1}{s(1-s)^2} \frac{2}{x^2(x^2-1)}]\\
 &\times F(1-s, -s; 3-2s; 1-x^2)\\
 &+\frac{1}{(1-s)^2} \frac{2}{x^2(x^2-1)} F(1-s, 1-s; 3-2s; 1-x^2)\\
 &+\frac{(2-s)}{s(1-s)^2} \frac{2}{x^4(x^2-1)} F(-s, -s; 3-2s; 1-x^2)
\endaligned$$
Note that
$$-\frac{1}{s(1-s)} \frac{2}{x^2-1}+\frac{1}{s(1-s)}
  \frac{2}{x^4(x^2-1)}+\frac{1}{s(1-s)} \frac{2(x^2+1)}{x^4}=0,$$
and
$$\aligned
 &\frac{1}{s(1-s)^2} \frac{2(x^2-2)}{x^4(x^2-1)}+\frac{1}{s(1-s)^2}
  \frac{2}{x^2(x^2-1)}\\
=&\frac{-1}{s(1-s)^2} \frac{4}{x^4(x^2-1)}+
  \frac{1}{s(1-s)^2} \frac{4}{x^2(x^2-1)}.
\endaligned$$
We have
$$\aligned
  S_{II}
=&\left[-\frac{(2-s)}{(1-s)^2} \frac{4}{x^2(x^2-1)}-
  \frac{1}{s(1-s)^2} \frac{4}{x^4(x^2-1)}+\frac{1}{s(1-s)^2}
  \frac{4}{x^2(x^2-1)}\right]\\
 &\times F(1-s, -s; 3-2s; 1-x^2)\\
 &+\frac{1}{(1-s)^2} \frac{2}{x^2(x^2-1)} F(1-s, 1-s; 3-2s; 1-x^2)\\
 &+\frac{2-s}{s(1-s)^2} \frac{2}{x^4(x^2-1)} F(-s, -s; 3-2s; 1-x^2).
\endaligned$$

  In the formula (see \cite{Er}, p.103, (36))
$$(c-a-b)F-(c-a)F(a-1)+b(1-z)F(b+1)=0,$$
set $a=1-s$, $b=-s$, $c=3-2s$ and $z=1-x^2$, then
$$\aligned
 &F(1-s, 1-s; 3-2s; 1-x^2)+\frac{2-s}{s} \frac{1}{x^2}
  F(-s, -s; 3-2s; 1-x^2)\\
=&\frac{1}{s} \frac{2}{x^2} F(1-s, -s; 3-2s; 1-x^2).
\endaligned$$
Now, we have
$$\aligned
 &\frac{1}{(1-s)^2} \frac{2}{x^2(x^2-1)} F(1-s, 1-s; 3-2s; 1-x^2)\\ 
 &+\frac{(2-s)}{s(1-s)^2} \frac{2}{x^4(x^2-1)} F(-s, -s; 3-2s; 1-x^2)\\
=&\frac{1}{(1-s)^2} \frac{2}{x^2(x^2-1)} \frac{1}{s} \frac{2}{x^2}
  F(1-s, -s; 3-2s; 1-x^2)\\
=&\frac{1}{s(1-s)^2} \frac{4}{x^4(x^2-1)} F(1-s, -s; 3-2s; 1-x^2).
\endaligned$$
Hence,
$$\aligned
  S_{II}
=&[-\frac{(2-s)}{(1-s)^2} \frac{4}{x^2(x^2-1)}-\frac{1}{s(1-s)^2}
  \frac{4}{x^4(x^2-1)}+\frac{1}{s(1-s)^2} \frac{4}{x^2(x^2-1)}\\
 &+\frac{1}{s(1-s)^2} \frac{4}{x^4(x^2-1)}] F(1-s, -s; 3-2s; 1-x^2)\\
=&\frac{1}{s} \frac{4}{x^2(x^2-1)} F(1-s, -s; 3-2s; 1-x^2)\\
=&\frac{4}{s} \frac{u^2}{u+1} F(1-s, -s; 3-2s; -u^{-1}).
\endaligned$$
Therefore, we have
$$\frac{\Gamma(s-\frac{1}{2})}{\sqrt{\pi} \Gamma(s)} (4u)^{s-3} S_{II}
 =\frac{\Gamma(s-\frac{1}{2})}{\sqrt{\pi} \Gamma(s+1)} 4^{s-2} u^{s-1}
  (u+1)^{-1} F(1-s, -s; 3-2s; -u^{-1}).$$
Consequently,
$$\phi(Z, Z^{\prime}; k, s)=K(Z, Z^{\prime}; k, s)+
                            K(Z, Z^{\prime}; k, 2-s),$$
where $k=\pm 1$.							
\flushpar
$\qquad \qquad \qquad \qquad \qquad \qquad \qquad \qquad \qquad
 \qquad \qquad \qquad \qquad \qquad \qquad \qquad \qquad \qquad
 \quad \boxed{}$

{\smc Corollary 4.12}. {\it The functional equation of Eisenstein 
series of weight $k$ for the trivial group on ${\frak S}_2$: If 
$\text{Re}(s)>1$, then the following identity holds:
$$\aligned
 &\int_{\partial {\frak S}_2} P_{k}(Z, W; s) 
  \rho(W, W^{\prime})^{s+k-2} \rho(\overline{W}, 
  \overline{W^{\prime}})^{s-k-2} dm(W)\\
=&\frac{\pi^{\frac{3}{2}} \Gamma(s-\frac{1}{2})}
  {(|k|+s-1) \Gamma(s)} 4^{s-1} P_{k}(Z, W^{\prime}; 2-s),
\endaligned\tag 4.31$$
for $k=0, \pm 1$, where $W^{\prime} \in \partial {\frak S}_2$.}

{\it Proof}. The integral on the left-hand side is absolutely
convergent. We set
$$P_{-k}(Z^{\prime}, W; 2-s)=\left[\frac{\rho(Z^{\prime}, W)}
  {\rho(\overline{Z^{\prime}}, \overline{W})}\right]^{-k}
  \frac{\rho(Z^{\prime})^{2-s}}{|\rho(Z^{\prime}, W)|^{2(2-s)}}$$
in Theorem 4.10. Next, we multiply both sides of the formula by 
$\rho(Z^{\prime})^{s-2}$ and take the limit as
$Z^{\prime} \to W^{\prime} \in \partial {\frak S}_2$. This completes
the proof.
\flushpar
$\qquad \qquad \qquad \qquad \qquad \qquad \qquad \qquad \qquad
 \qquad \qquad \qquad \qquad \qquad \qquad \qquad \qquad \qquad
 \quad \boxed{}$
 
  We define the S-matrix of weight $k$ as follows:
$$S(W, W^{\prime}; k, s)
:=\sum_{\gamma \in \Gamma} j(\gamma, W)^{k-s}
  j(\overline{\gamma}, \overline{W})^{-k-s} 
  \rho(\gamma(W), W^{\prime})^{-k-s} 
  \rho(\overline{\gamma(W)}, \overline{W^{\prime}})^{k-s}
  \tag 4.32$$  
for $W, W^{\prime} \in \Omega(\Gamma)$ and 
$\text{Re}(s)>\delta(\Gamma)$.

{\smc Theorem 4.13}. {\it Assume that $\Gamma$ is convex
cocompact and $\delta(\Gamma)<1$. Then the following 
functional equation for Eisenstein series of weight $k$ 
holds:
$$\aligned
 &\int_{\Gamma \backslash \Omega(\Gamma)} E(Z, W; k, s) 
  S(W, W^{\prime}; -k, 2-s) dm(W)\\
=&\frac{\pi^{\frac{3}{2}} \Gamma(s-\frac{1}{2})}
  {(|k|+s-1) \Gamma(s)} 4^{s-1} E(Z, W^{\prime}; k, 2-s),
\endaligned\tag 4.33$$
for $1<\text{Re}(s)<2-\delta(\Gamma)$ and $k=0, \pm 1$.}

{\it Proof}. By Corollary 4.12,
$$\aligned
 &\int_{\partial {\frak S}_2} E(Z, W; k, s) 
  \rho(W, W^{\prime})^{s+k-2} \rho(\overline{W}, 
  \overline{W^{\prime}})^{s-k-2} dm(W)\\
=&\frac{\pi^{\frac{3}{2}} \Gamma(s-\frac{1}{2})}
  {(|k|+s-1) \Gamma(s)} 4^{s-1} E(Z, W^{\prime}; k, 2-s).
\endaligned\tag 4.34$$
On the other hand,
$$\aligned
  \int_{\partial {\frak S}_2} f(W) dW
=&\int_{\Gamma \backslash \Omega(\Gamma)} \sum_{\gamma \in \Gamma} 
  f(\gamma(W)) dm(\gamma(W))\\
=&\int_{\Gamma \backslash \Omega(\Gamma)} \sum_{\gamma \in \Gamma} 
  |J(\gamma(W))| f(\gamma(W)) dm(W).
\endaligned$$

  By Lemma 4.9, (3), the left hand side of (4.34) is equal to
$$\aligned
 &\int_{\Gamma \backslash \Omega(\Gamma)} \sum_{\gamma \in \Gamma}
  |J(\gamma(W))| E(Z, \gamma(W); k, s) \rho(\gamma(W), W^{\prime})^
  {s+k-2} \rho(\overline{\gamma(W)}, \overline{W^{\prime}})^{s-k-2} 
  dm(W)\\
=&\int_{\Gamma \backslash \Omega(\Gamma)} \sum_{\gamma \in \Gamma}
  |j(\gamma, W)|^{-4} j(\gamma, W)^{s-k} j(\overline{\gamma},
  \overline{W})^{s+k} E(Z, W; k, s)\\
 &\times \rho(\gamma(W), W^{\prime})^{s+k-2} 
  \rho(\overline{\gamma(W)}, \overline{W^{\prime}})^{s-k-2} dm(W)\\
=&\int_{\Gamma \backslash \Omega(\Gamma)} E(Z, W; k, s) 
  S(W, W^{\prime}; -k, 2-s) dm(W).
\endaligned$$  
$\qquad \qquad \qquad \qquad \qquad \qquad \qquad \qquad \qquad
 \qquad \qquad \qquad \qquad \qquad \qquad \qquad \qquad \qquad
 \quad \boxed{}$
 
  By Theorem 4.10, we get the main theorem of this paper.

{\smc Theorem 4.14}. {\it Assume that $\Gamma$ is convex cocompact
and $\delta(\Gamma)<1$. Then the following product formula holds:
$$\aligned
 &\int_{\Gamma \backslash \Omega(\Gamma)} E(Z, W; k, s)
  E(Z^{\prime}, W; -k, 2-s) dm(W)\\
=&G(Z, Z^{\prime}; k, s)+G(Z, Z^{\prime}; k, 2-s)\\
=&G(Z^{\prime}, Z; -k, s)+G(Z^{\prime}, Z; -k, 2-s),
\endaligned\tag 4.35$$
for $\delta(\Gamma)<\text{Re}(s)<2-\delta(\Gamma)$
and $k=0$ or $k=\pm 1$.}

\vskip 0.5 cm
\centerline{\bf Appendix. Product formulas on $SL(2, {\Bbb R})$}
\vskip 0.5 cm

  In this appendix, by the same method as in the above argument
we give the product formulas on $SL(2, {\Bbb R})$. For simplicity, 
we omit some details.

  Let
$$\Delta_{k}=-(z-\overline{z})^{2} \frac{\partial^2}{\partial z
  \partial \overline{z}}-k(z-\overline{z})(\frac{\partial}
  {\partial z}+\frac{\partial}{\partial \overline{z}}).$$  
Set $H_{k}(z, w)=(\overline{z}-w)^{k}(z-\overline{w})^{-k}$ and 
$P(z, \zeta)=\frac{\text{Im}(z)}{|\overline{z}-\zeta|^{2}}$, 
where $z, w \in {\Bbb H}$ and $\zeta \in \partial {\Bbb H}={\Bbb R}$. 
The Poisson kernel of weight $k$ is defined as follows:
$$P_{k}(z, \zeta; s)=(\overline{z}-\zeta)^{k} (z-\zeta)^{-k}
  P(z, \zeta)^{s}.$$
$$E(z, \zeta; k, s):=\sum_{\gamma \in \Gamma} 
  P_{k}(\gamma(z), \zeta; k, s), \quad \text{for} \quad 
  \text{Re}(s)>\delta(\Gamma).$$
Here $\Gamma$ is a discrete subgroup of $SL(2, {\Bbb R})$.
The point-pair invariant
$$u(z, z^{\prime}):=\frac{|z-z^{\prime}|^{2}}{4 \text{Im}(z) 
  \text{Im}(z^{\prime})}, \quad
  \sigma(z, z^{\prime}):=\frac{|\overline{z}-z^{\prime}|^{2}}
  {4 \text{Im}(z) \text{Im}(z^{\prime})}.$$
In fact, $u=\sigma-1$.

  If $f=H_{k}(z, w) \Phi(\sigma(z, w))$, then
$$\Delta_{k} f=H_{k}(z, w)\left[\sigma(\sigma-1) 
  \Phi^{\prime \prime}(\sigma)+(2 \sigma-1) 
  \Phi^{\prime}(\sigma)+\frac{k^2}{\sigma} \Phi(\sigma)\right].$$    
By $\Delta_{k}f=s(s-1) f$, we have
$$\sigma(\sigma-1) \Phi^{\prime \prime}(\sigma)+(2 \sigma-1)
  \Phi^{\prime}(\sigma)+\left[\frac{k^2}{\sigma}-s(s-1)\right]
  \Phi(\sigma)=0.$$
A solution is
$$\aligned
  \Phi(\sigma)
=&P\left\{\matrix
    0 & 1 & \infty\\
 -|k| & 0 & s     \\
  |k| & 0 & 1-s 
  \endmatrix; \sigma \right\}\\ 
=&\sigma^{-|k|}(1-\sigma)^{|k|-s} F(s-|k|, s-|k|; 2s; 
  -\frac{1}{\sigma-1}). 
\endaligned$$
Set
$$\aligned
  K(z, z^{\prime}; k, s)
=&\frac{\sqrt{\pi} \Gamma(\frac{1}{2}-s)}{\Gamma(1-s)}
  4^{-s} \frac{s}{s-|k|} H_{k}(z, z^{\prime})\\
 &\times u^{|k|-s}(1+u)^{-|k|} F(s-|k|, s-|k|; 2s; -u^{-1}),
\endaligned$$
where $u=u(z, z^{\prime})$. The automorphic Green function of
weight $k$ is given by
$$G(z, z^{\prime}; k, s)=\sum_{\gamma \in \Gamma}
  K(z, \gamma(z^{\prime}); k, s), \quad \text{for} \quad
  \text{Re}(s)>\delta(\Gamma).$$
Let
$$\phi(z, z^{\prime}; k, s)=\int_{\Bbb R} P_{k}(z, \zeta; s)
  P_{-k}(z^{\prime}, \zeta; 1-s) d \zeta.$$
Then $\phi$ is a point-pair invariant covariant with respect to the
weight $k$. 

  Put $z=iy$, $z^{\prime}=i y^{\prime}$, $y>0$, $y^{\prime}>0$. Then
$$\aligned
 &\phi(iy, iy^{\prime}; k, s)\\
=&y^{s} {y^{\prime}}^{1-s} \int_{-\infty}^{\infty}
  \left(\frac{\zeta+iy}{\zeta-iy}\right)^{k}
  \left(\frac{\zeta+i y^{\prime}}{\zeta-i y^{\prime}}\right)^{-k}
  |\zeta-iy|^{-2s} |\zeta-i y^{\prime}|^{-2(1-s)} d \zeta\\
=&y^{s} {y^{\prime}}^{1-s} \int_{-\infty}^{\infty} 
  \exp(2ki \arctan \frac{y}{\zeta}-2ki \arctan 
  \frac{y^{\prime}}{\zeta}) \frac{1}{(\zeta^2+y^2)^{s}
  ({y^{\prime}}^{2}+\zeta^2)^{1-s}} d \zeta.
\endaligned$$
By $\tan(\arctan \frac{y}{\zeta}-\arctan \frac{y^{\prime}}{\zeta})
=\frac{\zeta(y-y^{\prime})}{\zeta^2+y y^{\prime}}$,
$$\cos(2 \arctan \left[\frac{(y-y^{\prime}) \zeta}{y y^{\prime}+\zeta^2} 
  \right])=\frac{(y y^{\prime}+\zeta^2)^{2}-(y-y^{\prime})^{2} \zeta^2}
  {(y y^{\prime}+\zeta^2)^{2}+(y-y^{\prime})^{2} \zeta^{2}},$$
and
$$\aligned
 &\cos(2k \arctan \left[\frac{(y-y^{\prime}) \zeta}{y y^{\prime}+
  \zeta^2} \right])\\
=&\cos(k \arccos \frac{(y y^{\prime}+\zeta^2)^{2}-(y-y^{\prime})^{2}
  \zeta^2}{(y y^{\prime}+\zeta^2)^{2}+(y-y^{\prime})^{2} \zeta^{2}})\\ 
=&F(-k, k; \frac{1}{2}; \frac{1}{2}\left[1-\frac{(y y^{\prime}+
 \zeta^2)^{2}-(y-y^{\prime})^{2} \zeta^2}{(y y^{\prime}+\zeta^2)^2+
 (y-y^{\prime})^2 \zeta^2} \right])\\
=&F(-k, k; \frac{1}{2}; \frac{(y-y^{\prime})^{2} \zeta^{2}}
  {(y^2+\zeta^2)({y^{\prime}}^2+\zeta^2)}), 
\endaligned$$
we have
$$\aligned
 \phi(iy, iy^{\prime}; k, s)
=&2 y^{s} {y^{\prime}}^{1-s} \sum_{n=0}^{|k|}
  \frac{(-k)_{n} (k)_{n}}{n! (\frac{1}{2})_{n}} (y-y^{\prime})^{2n}\\
 &\times \int_{0}^{\infty} \zeta^{2n} (\zeta^2+y^2)^{-s-n}
  (\zeta^2+{y^{\prime}}^2)^{s-n-1} d \zeta.
\endaligned$$
Set $l=\frac{y^{\prime}}{y}$, then $u=\frac{(l-1)^{2}}{4l}$.
By the same method as the above argument for $GU(2, 1)$, we have
$$\aligned
 &\phi(iy, iy^{\prime}; k, s)\\
=&\sum_{n=0}^{|k|} \frac{(-k)_{n} (k)_{n}}{n! (\frac{1}{2})_{n}}
  \Gamma(n+\frac{1}{2})[\frac{\Gamma(\frac{1}{2}-s)}{\Gamma(n+1-s)}
  (4u)^{-s} F(n+s, s; 2s; -u^{-1})\\
 &+\frac{\Gamma(s-\frac{1}{2})}{\Gamma(n+s)} (4u)^{-(1-s)}
  F(n+1-s, 1-s; 2-2s; -u^{-1})].
\endaligned$$
\flushpar
(1) $k=0$, 
$$\phi(iy, iy^{\prime}; k, s)=\frac{\sqrt{\pi} \Gamma(\frac{1}{2}-s)}
  {\Gamma(1-s)} (4u)^{-s} F(s, s; 2s; -u^{-1})+(s \mapsto 1-s).$$
(2) $k=\pm 1$,
$$S_{I}:=\frac{\sqrt{\pi} \Gamma(\frac{1}{2}-s)}{\Gamma(1-s)}
  (4u)^{-s}[F(s, s; 2s; -u^{-1})+\frac{1}{s-1} 
  F(s+1, s; 2s; -u^{-1})].$$
By \cite{Er}, p.103, (28),
$$(c-a)F(a-1)+(2a-c-az+bz)F+a(z-1)F(a+1)=0,$$
set $a=s$, $b=s$, $c=2s$ and $z=-u^{-1}$, then
$$F(s+1, s; 2s; -u^{-1})=\frac{u}{u+1} F(s-1, s; 2s; -u^{-1}).$$
By \cite{Er}, p.103, (36),
$$(c-a-b)F-(c-a)F(a-1)+b(1-z)F(b+1)=0,$$
set $a=s$, $b=s-1$, $c=2s$ and $z=-u^{-1}$, then
$$\aligned
 &\frac{1}{s-1} \frac{u}{u+1} F(s, s-1; 2s; -u^{-1})+
  F(s, s; 2s; -u^{-1})\\
=&\frac{s}{s-1} \frac{u}{u+1} F(s-1, s-1; 2s; -u^{-1}).
\endaligned$$
Thus,
$$S_{I}=\frac{\sqrt{\pi} \Gamma(\frac{1}{2}-s)}{\Gamma(1-s)}
  (4u)^{-s} \frac{s}{s-1} \frac{u}{u+1} F(s-1, s-1; 2s; -u^{-1}).$$
Consequently,
$$\phi(iy, iy^{\prime}; k, s)
 =\frac{\sqrt{\pi} \Gamma(\frac{1}{2}-s)}{\Gamma(1-s)}
  4^{-s} \frac{s}{s-1} u^{1-s} (u+1)^{-1} F(s-1, s-1; 2s; -u^{-1})
  +(s \mapsto 1-s).$$
Thus, we have the following theorem:

{\smc Theorem}. {\it Assume that $\Gamma$ is convex cocompact
and $\delta(\Gamma)<\frac{1}{2}$. Then the following product 
formula on $SL(2, {\Bbb R})$ holds:
$$\int_{\Gamma \backslash \Omega(\Gamma)} E(z, \zeta; k, s)
  E(z^{\prime}, \zeta; -k, 1-s) dm(\zeta)
 =G(z, z^{\prime}; k, s)+G(z, z^{\prime}; k, 1-s),$$
for $\delta(\Gamma)<\text{Re}(s)<1-\delta(\Gamma)$ and $k=0, \pm 1$.}

\vskip 2.0 cm
{\smc Department of Mathematics, Peking University}

{\smc Beijing 100871, P. R. China}

{\it E-mail address}: yanglei\@sxx0.math.pku.edu.cn
\vskip 1.5 cm
\Refs

\item{[Ap]} {\smc B. Apanasov}, Deformations and stability in complex
	        hyperbolic geometry, MSRI preprint, No. 1997-111. 

\item{[BB]} {\smc W. L. Baily and A. Borel}, Compactification of
            arithmetic quotients of bounded symmetric domains,
			Ann. of Math. {\bf 84} (1966), 442-528. 

\item{[BO]} {\smc U. Bunke and M. Olbrich}, The spectrum of
             Kleinian manifolds, math.DG/9810146.

\item{[Co]} {\smc K. Corlette}, Hausdorff dimensions of limit sets
             I, Invent. Math. {\bf 102} (1990), 521-542.

\item{[CoI]} {\smc K. Corlette and A. Iozzi}, Limit sets of
              discrete groups of isometries of exotic hyperbolic
			  spaces, Trans. AMS {\bf 351} (1999), 1507-1530.

\item{[El]} {\smc J. Elstrodt}, Die Resolvente zum Eigenwertproblem
            der automorphen Formen in der hyperbolischen Ebene, I,
			Math. Ann. {\bf 203} (1973), 295-330; II, Math. Z. 
			{\bf 132} (1973), 99-134; III, Math. Ann. {\bf 208} 
			(1974), 99-132.

\item{[ElGM]} {\smc J. Elstrodt, F. Grunewald and J. Mennicke}, {\it
            Groups Acting on Hyperbolic Space, Harmonic Analysis and
			Number Theory}, Springer-Verlag, 1998. 

\item{[Er]} {\smc A. Erd\'{e}lyi, W. Magnus, F. Oberhettinger and
            F. G. Tricomi}, {\it Higher Transcendental Functions},
            Vol. I, McGraw-Hill Book Company, Inc., 1953.

\item{[Fa]} {\smc J. D. Fay}, Fourier coefficients of the resolvent
            for a Fuchsian group, J. Reine Angew. Math. {\bf 294}
			(1977), 143-203.

\item{[FHP]} {\smc R. Froese, P. Hislop and P. Perry}, The Laplace
            operator on hyperbolic three manifolds with cusps of
			non-maximal rank, Invent. Math. {\bf 106} (1991),
			295-333.

\item{[Gir]} {\smc G. Giraud}, Sur certaines fonctions automorphes
           de deux variables, Ann. Ec. Norm., {\bf (3)}, {\bf 38}
           (1921), 43-164.

\item{[Go]} {\smc W. M. Goldman}, {\it Complex Hyperbolic 
            Geometry}, Oxford Mathematical Monographs,
			Oxford, 1999. 
		   
\item{[H]} {\smc D. Hejhal}, {\it The Selberg Trace Formula for
           $PSL(2, {\Bbb R})$}, Vol. {\bf 1}, Lecture Notes in Math.
		   {\bf 548}, Springer-Verlag, 1976; Vol. {\bf 2}, 
		   Lecture Notes in Math. {\bf 1001}, Springer-Verlag, 1983.
		   
\item{[He]} {\smc J. C. Hemperly}, The parabolic contribution to the
             number of linearly independent automorphic forms on a
			 certain bounded domain, Amer. J. Math. {\bf 94} (1972),
			 1078-1100.
			 
\item{[HP]} {\smc S. Hersonsky and F. Paulin}, On the volumes
            of complex hyperbolic manifolds, Duke Math. J.
            {\bf 84} (1996), 719-737.

\item{[Ho1]} {\smc R. P. Holzapfel}, {\it Geometry and Arithmetic
            Around Euler Partial Differential Equations}, Dt. Verl.
            d. Wiss., Berlin/Reidel, Dordrecht 1986.

\item{[Ho2]} {\smc R. P. Holzapfel}, {\it The Ball and Some Hilbert
            Problems}, Lectures in Mathematics, ETH Z\"{u}rich,
            Birkh\"{a}user, 1995.
	
\item{[K]} {\smc M. Kashiwara, A. Kowata, K. Minemura, K. Okamoto, 
            T. Oshima and M. Tanaka}, Eigenfunctions of invariant 
			differential operators on a symmetric space, 
			Ann. of Math. 107(1978), 1-39. 

\item{[Ko]} {\smc A. Kor\'{a}nyi}, The Poisson integral for
             generalized half-planes and bounded symmetric domains,
			 Ann. of Math. {\bf 82} (1965), 332-350.
			 
\item{[KR]} {\smc A. Kor\'{a}nyi and H. M. Reimann}, Quasiconformal
             mappings on the Heisenberg group, Invent. Math. {\bf 80}
			 (1985), 309-338.
			 
\item{[KW]} {\smc A. Kor\'{a}nyi and J. A. Wolf}, Realization of
             hermitian symmetric spaces as generalized half-planes,
			 Ann. of Math. {\bf 81} (1965), 265-288.

\item{[Ku]} {\smc T. Kubota}, {\it Elementary Theory of Eisenstein 
            Series}. Wiley, New York, 1973.        

\item{[LR]} {\smc R. P. Langlands and D. Ramakrishnan}, {\it The
            Zeta Functions of Picard Modular Surfaces}, Les
            Publications CRM, Universit\'{e} de Montr\'{e}al,
            1992.

\item{[M]} {\smc H. Maass}, Die Differentialgleichungen in der 
            Theorie der elliptischen Modulfunktionen, Math. Ann.
			{\bf 125} (1953), 235-263. 

\item{[Ma1]} {\smc N. Mandouvalos}, Spectral theory and Eisenstein
            series for Kleinian groups, Proc. London Math. Soc.
			(3) {\bf 57} (1988), 209-238.

\item{[Ma2]} {\smc N. Mandouvalos}, {\it Scattering Operator, 
            Eisenstein Series, Inner Product Formula and 
			``Maass-Selberg'' Relations for Kleinian Groups}, 
			Memoirs of A.M.S. {\bf 78},	No. {\bf 400} (1989).			

\item{[P]} {\smc J. R. Parker}, On the volumes of cusped, complex
           hyperbolic manifolds and orbifolds, Duke Math. J.
           {\bf 94} (1998), 433-464.

\item{[Pa1]} {\smc S. J. Patterson}, The limit set of a Fuchsian
              group, Acta Math. {\bf 136} (1976), 241-273.

\item{[Pa2]} {\smc S. J. Patterson}, The Laplacian operator on a
              Riemann surface, I, Comp. Math. {\bf 31} (1975),
			  83-107; II, Comp. Math. {\bf 32} (1976), 71-119,
			  III, Comp. Math. {\bf 33} (1976), 227-259.

\item{[Pa3]} {\smc S. J. Patterson}, The Selberg zeta function of a
           Kleinian group. In {\it Number Theory, Trace Formulas,
		   and Discrete Groups}: Symposium in honor of Atle Selberg,
		   Olso, Norway (1987), 409-442.

\item{[Pe1]} {\smc P. A. Perry}, The Laplace operator on a hyperbolic
           manifold. II. Eisenstein series and the scattering matrix,
		   J. Reine Angew. Math. {\bf 398} (1989), 67-91.
		    
\item{[Pe2]} {\smc P. A. Perry}, The Selberg zeta function and a local
           trace formula for Kleinian groups, J. Reine Angew. Math.
		   {\bf 410} (1990), 116-152.

\item{[Pic1]} {\smc E. Picard}, Sur une extension aux fonctions de
           deux variables du probl\`{e}me de Riemann relatif aux
           fonctions hypergeom\'{e}triques s\'{e}ries
           hypergeom\'{e}triques, Ann. Ec. Norm. {\bf 10} (1881),
           305-322.

\item{[Pic2]} {\smc E. Picard}, Sur les fonctions hyperfuchsiennes
           provenant des s\'{e}ries hypergeom\'{e}triques de deux
           variables, Ann. Ec. Norm. {\bf III 2} (1885), 357-384.

\item{[Ro]} {\smc W. Roelcke}, Das Eigenwertproblem der automorphen
             Formen in der hyperbolischen Ebene, I, Math. Ann. 
			 {\bf 167} (1966), 292-337; II, Math. Ann. {\bf 168} 
			 (1967), 261-324.

\item{[Se]} {\smc A. Selberg}, Harmonic analysis and discontinuous
            groups in weakly symmetric Riemannian spaces with 
			applications to Dirichlet series, J. Indian Math. Soc.
			B. {\bf 20} (1956), 47-87.
			
\item{[Sh1]} {\smc G. Shimura}, The arithmetic of automorphic forms
             with respect to a unitary group, Ann. of Math. 
			 {\bf 107} (1978), 569-605.
			 
\item{[Sh2]} {\smc G. Shimura}, L-functions and eigenvalue problems.
             In {\it Algebraic Analysis, Geometry, and Number Theory},
             (J. I. Igusa, ed.), Johns Hopkins University Press,
			 Baltimore, 1989, 341-396. 

\item{[Su]} {\smc D. Sullivan}, The density at infinity of a discrete
             group of hyperbolic motions, Publ. Math. I.H.E.S. 
			 {\bf 50} (1979), 171-209.

\endRefs
\end{document}